\documentclass[10pt]{scrartcl}
\usepackage[
bookmarks,
  bookmarksopen=true,
  bookmarksnumbered=true,
  pdfusetitle,
  pdfcreator={},
  colorlinks,
  linkcolor=black,
  urlcolor=black,
  citecolor=black,
  plainpages=false,
  ]{hyperref}
	\usepackage[
  left=2.5cm,
  right=2.5cm,
  top=2cm,
]{geometry}

\usepackage{amsfonts}
\usepackage{amsmath}
\usepackage{amssymb}
\usepackage{amsthm}
\usepackage{algorithm}
\usepackage{algorithmic}
\usepackage{graphicx}
\usepackage{color}
\usepackage{subcaption}
\usepackage{xfrac}
\usepackage{csvsimple} 

\usepackage{mathabx}

\usepackage{tikz}
\usepackage{pgfplots}\pgfplotsset{compat=1.9}
\usepackage{tikz-3dplot}
\usetikzlibrary{plotmarks}
\usetikzlibrary{cd}
\usetikzlibrary{decorations.pathreplacing}

\def\addlegendimage{\csname pgfplots@addlegendimage\endcsname}

\usepackage{prettyref}
\usepackage{mathrsfs}
\usepackage[utf8]{inputenc}
\usepackage{floatflt} 
\usepackage{float} 
\usepackage{wrapfig}

\usepackage{chngcntr}
\counterwithin{figure}{section}
\counterwithin{table}{section}

\renewcommand{\vec}[1]{\ensuremath \mathbf{\boldsymbol{#1}}}
\renewcommand{\r}[1]{\ensuremath \mathring{#1}}

\renewcommand{\phi}{\varphi}
\renewcommand{\epsilon}{\varepsilon}
\renewcommand{\rho}{\varrho}
\newcommand{\bB}{\mathbf{\boldsymbol{B}}}

\newcommand{\Rho}{\mathrm{R}}
\newcommand{\N}{\ensuremath{\mathbb{N}}}
\newcommand{\T}{\ensuremath{\mathbb{T}}}

\newcommand{\I}{\ensuremath{\mathcal{I}}}
\newcommand{\J}{\ensuremath{\mathcal{J}}}

\newcommand{\Z}{\ensuremath{\mathbb{Z}}}

\newcommand{\R}{\ensuremath{\mathbb{R}}}
\newcommand{\C}{\ensuremath{\mathbb{C}}}
\newcommand{\E}{\ensuremath{\mathbb{E}}}
\newcommand{\Eps}{\ensuremath{\mathcal{E}}}

\newcommand{\per}{\mathrm{per}}
\newcommand{\erf}{\mathrm{erf}}
\newcommand{\MISE}{\text{MISE}}
\newcommand{\AMISE}{\text{AMISE}}

\renewcommand{\P}{\ensuremath{\mathbb{P}}}

\newcommand{\dx}{\mathrm{d}}
\newcommand{\e}{\mathrm{e}}
\newcommand{\im}{\mathrm{i}}
\newcommand{\mix}{\mathrm{mix}}

\DeclareMathOperator*{\supp}{supp}

\DeclareMathOperator*{\argmin}{argmin}

\DeclareMathOperator*{\lin}{span}

\DeclareMathOperator{\sgn}{sgn}
\renewcommand{\d}{\,\mathrm{d}}
\newcommand{\Dx}{\mathrm{D}}

\newcommand{\norm}[1]{\left\lVert \smash{#1} \right\rVert}

\newcommand\multiset[2]
{\mathchoice{\left(\kern-0.5em{\binom{#1}{#2}}\kern-0.5em\right)}
	{\bigl(\kern-0.3em{\binom{#1}{#2}}\kern-0.3em\bigr)}
	{\bigl(\kern-0.3em{\binom{#1}{#2}}\kern-0.3em\bigr)}
	{\bigl(\kern-0.3em{\binom{#1}{#2}}\kern-0.3em\bigr)}}
	
\newcommand{\X}{\ensuremath{\mathcal{X}}}
\newcommand{\Y}{\ensuremath{\mathcal{Y}}}

\makeatletter
\def\moverlay{\mathpalette\mov@rlay}
\def\mov@rlay#1#2{\leavevmode\vtop{\baselineskip\z@skip \lineskiplimit-\maxdimen
   \ialign{\hfil$\m@th#1##$\hfil\cr#2\crcr}}}
\newcommand{\charfusion}[3][\mathord]{
    #1{\ifx#1\mathop\vphantom{#2}\fi
        \mathpalette\mov@rlay{#2\cr#3}
      }
    \ifx#1\mathop\expandafter\displaylimits\fi}
\makeatother
\newcommand{\cupdot}{\charfusion[\mathbin]{\cup}{\cdot}}
\newcommand{\bigcupdot}{\charfusion[\mathop]{\bigcup}{\cdot}}

\newtheorem{theorem}{Theorem}[section]
\newtheorem{lemma}[theorem]{Lemma}

\newtheorem{corollary}[theorem]{Corollary}
\newtheorem{proposition}[theorem]{Proposition}
\newtheorem{problem}[theorem]{Problem}

\theoremstyle{definition}
\newtheorem{definition}[theorem]{Definition}
\newtheorem{example}[theorem]{Example}
\newtheorem{remark}[theorem]{Remark}

\newenvironment{Remark}[1][noisnotdefined]{ \ifthenelse{\equal{#1}{noisnotdefined}}{\begin{remark}}{\begin{remark}[#1]}\normalfont\rmfamily}{\bend\end{remark}}
\newenvironment{Example}{ \begin{example}\normalfont\rmfamily}{\bend\end{example}}
\newenvironment{Definition}{ \begin{definition}\normalfont\rmfamily}{\end{definition}}

\numberwithin{equation}{section}

\newcommand{\bend}{\hspace*{0ex} \hfill \hbox{\vrule height
    1.5ex\vbox{\hrule width 1.4ex \vskip 1.4ex\hrule  width 1.4ex}\vrule
    height 1.5ex}}

\long\def\symbolfootnote[#1]#2{\begingroup
\def\thefootnote{\fnsymbol{footnote}}\footnote[#1]{#2}\endgroup}

\newcounter{todocounter}
\newcommand{\todo}[2][noisnotdefined]{
 \marginpar{\fcolorbox{black}{yellow}{\footnotesize\textbf{todo}}
 \ifthenelse{\equal{#1}{noisnotdefined}}{}{\textcolor{black}{\newline\tiny #1}}}
 \textbf{\ifthenelse{\equal{#2}{.}}
   {\fcolorbox{red}{white}{\textcolor{red}{$\maltese$}}}{{\textcolor{red}{#2}}}}
 \refstepcounter{todocounter}}

\newrefformat{def}{Definition \ref{#1}}
\newrefformat{sub}{Subsection \ref{#1}}
\newrefformat{prop}{Proposition \ref{#1}}
\newrefformat{fig}{Figure \ref{#1}}
\newrefformat{rem}{Remark \ref{#1}}
\newrefformat{cor}{Corollary \ref{#1}}

\definecolor{tuccolor}{rgb}{0, 0.371, 0.3125}

\newcommand{\changed}[1]{\textcolor{black}{#1}} 

\title{Variable Transformations in combination with Wavelets and ANOVA for high-dimensional approximation}

\date{\today}
\date{}
\author{\textbf{Daniel Potts}\\Faculty of Mathematics\\Chemnitz University of Technology\\09107 Chemnitz \\
Germany\\
potts@math.tu-chemnitz.de\\
ORCID: 0000-0003-3651-4364
\and
\textbf{Laura Weidensager\thanks{corresponding author}}\\Faculty of Mathematics\\Chemnitz University of Technology\\09107 Chemnitz \\
Germany\\
laura.weidensager@math.tu-chemnitz.de\\
ORCID: 0000-0001-7988-1485
}

\begin{document}
\maketitle
\begin{abstract}
\subsection*{Abstract}
We use hyperbolic wavelet regression for the fast reconstruction of high-dimensional functions having only low dimensional variable interactions. 
Compactly supported periodic Chui-Wang wavelets are used for the tensorized
hyperbolic wavelet basis on the torus. With a variable transformation we are able 
to transform the approximation rates and fast algorithms from the torus to other domains. 
We perform and analyze scattered-data approximation for smooth but arbitrary density functions by using a least squares
method. The corresponding system matrix is sparse due to the compact support of the wavelets, which leads to a significant acceleration of the matrix
vector multiplication. 
For non-periodic functions we propose a new extension method. A proper choice of the extension parameter together with the piece-wise polynomial Chui-Wang wavelets extends the functions appropriately. In every case we are able to bound the approximation error 
with high probability. Additionally, if the function
has low effective dimension (i.e.\,only interactions of few variables), we qualitatively determine the variable interactions and omit ANOVA terms with low
variance in a second step in order to decrease the approximation error. This allows us to
suggest an adapted model for the approximation. Numerical results show the
efficiency of the proposed method.
\end{abstract}
 \textit{Keywords:} Variable transformations, least squares approximation, random sampling, wavelets, ANOVA decomposition
\newpage 
\section{Introduction}
The distribution of data points is a key component in machine learning. In many applications a target variable has to be predicted from given high-dimensional 
samples. We want to reconstruct an underlying function $f$ to give an interpretable approximation algorithm, which allows a prediction of the target 
variable for new samples. We consider the domains $\Omega\in \{\T^d,\R^d,[0,1]^d\}$ and also tensor products of these cases. We consider the setting of reconstructing a $d$-dimensional function $f\colon \Omega\rightarrow \C$ from discrete samples on the set of nodes $\{\vec y_1,\ldots,\vec y_M\}\subset \Omega$, which are distributed to the continuous density $\rho:\Omega\rightarrow \R_{+}$. One main aim is to also deal with an unknown density. 
Besides the natural question of finding a good approximation for $f$, we want to consider the question
of interpretability, i.e.\,analyzing the importance of the input variables and variable interactions of the function. \\

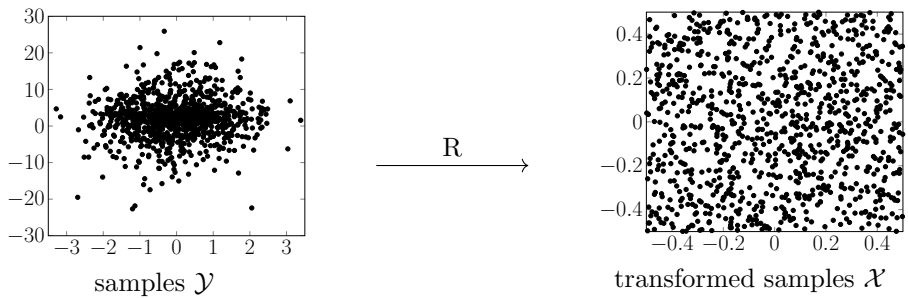
\begin{figure}[htb]
\centering
\begin{minipage}{0.3\textwidth}
\centering
\begin{tikzpicture}[scale = 0.4,every node/.append style={font=\huge}]
    \begin{axis}[scale only axis,
					ymin=-30,
						ymax = 30,
						xmin = -3.5,
						xmax = 3.5,
         ]        
          \addplot[black,only marks] table [x index=0,y index=1,col sep=comma] {Y_12.csv};

    \end{axis}
\end{tikzpicture}\\
samples $\Y$
\end{minipage}\begin{minipage}{0.19\textwidth}
\centering

\begin{tikzpicture}
\draw[->] (0,0) -- node[above]{$\Rho$} (2, 0)
 ;
\end{tikzpicture}
\end{minipage}\begin{minipage}{0.3\textwidth}
\centering
\begin{tikzpicture}[scale=0.4,every node/.append style={font=\huge}]
    \begin{axis}[scale only axis,
						ymin=-0.5,
						ymax = 0.5,
						xmin = -0.5,
						xmax = 0.5,						
         ]        
          \addplot[black,only marks] table [x index=0,y index=1,col sep=comma] {X_12.csv};

    \end{axis}
\end{tikzpicture}\\
transformed samples $\X$
\end{minipage}
\caption{Transformation of the samples in the two-dimensional case.}
\label{fig:samples}
\end{figure}

\textbf{Motivation}\\
The starting point of our considerations is the question whether it is possible to transform the good approximation results and the related fast algorithms for periodic functions on the torus $\T^d$ to the domain $\Omega$. To investigate the scattered data problem on the torus we are engaged with the sample set $\X$, the corresponding function values $\vec f = (f(\vec x))_{\vec x\in \X}$ and we constructed a recovery operator $S_I^\X$ in \cite{LiPoUl21}. This operator computes  
a best least squares fit 
\begin{equation}\label{eq:S_IX}
S_I^\X f = \sum_{\vec k\in I}a_{\vec k} \psi_{\vec k},
\end{equation} 
in the finite dimensional subspace spanned by semi-orthogonal wavelets $\psi_{\vec k} \colon \T^d \rightarrow \R$ with indices in the hyperbolic cross type set $I$. 
Assuming i.i.d.\,uniformly samples $\X\subset\T^d$ we showed in~\cite[Corollary~3.22.]{LiPoUl21}: 
Let $m$ be the order of vanishing moments of the wavelets and the sample size have logarithmic oversampling, i.e. $|\X|\gtrsim r |I|\log |I|$. For $1/2<s<m$ there is a constant $C(r,d,s)>0$ such that for fixed Besov norm $\|f\|_{\bB^s_{2,\infty}(\T^d)} \leq 1$
$$\P\Big(\norm{f-S_I^\X f}_{L_2(\T^d)}\leq C(r,d,s) \,\frac{(\log |\X|)^{(d-1)(s+1/2)+s}}{|\X|^s}\Big) \geq 1-2|\X|^{-r}\,,$$
see Appendix~\ref{sec:A1} for the definition of the Besov space.
Also in~\cite{PoSc19a} uniformly i.i.d.\,samples on the torus in combination with different basis functions perform well and there are possibilities for dealing with the curse of dimensionality. \changed{The results in~\cite{LiPoUl21} for the periodic case} serve as basis for this paper. 
In many practical applications, we have to take the data set as it is and have no uniform samples available.  
For that reason, we study here the case where the given sample points $\Y\subset\Omega$ are sampled from an arbitrary 
(but possibly unknown) density $\rho(\vec y)$. In Figure~\ref{fig:samples} we illustrate some random two-dimensional samples with respect to $\rho$. We can not guarantee good approximation rates and stability if we would use these samples directly.\\

To this end, we investigate transformations $\Rho$ of the given samples. 
The main result of this paper is the transformed approximation operator $S^\Y_I f$. Besides the interesting results for $\Omega \in \{\T^d,\R^d\}$ we use a new extension technique for the non-periodic case $\Omega = [0,1]^d$. Furthermore we present a detailed error analysis. \\

\textbf{The approach}\\
Our main approach is to transform given samples $\Y\subset \Omega$ to the d-dimensional torus $\T^d$ by $\X = \Rho(\Y)$, \changed{using the idea of inverse transform sampling: Let $F$ be the cumulative distribution function of a distribution $\rho$ and $\X\sim \mathcal U([0,1])$. Then the random variable $F^{-1}(\X) \sim \rho$ is distributed according to the distribution $\rho$. Based on this, we give possibilities for constructing a transformation $\R$ in~\eqref{eq:Rho}, which transforms the samples $\Y\subset \Omega^d$ to $\X\subset \T^d$ on the torus.}
In Figure~\ref{fig:samples} we show an illustration of what our constructed transformation $\Rho$ does with the samples.  
In order to investigate the scattered data problem on $\Omega$, we then use the recovery operator $S_I^\X$ from \eqref{eq:S_IX} on the torus. This operator minimizes the $\ell_2$-loss function
$$\sum_{\vec y\in \Y}|f\circ \Rho^{-1}(\vec y)-S_I^\X \left(f\circ \Rho^{-1}(\vec y)\right)|^2$$ 
by using an iterative LSQR algorithm. To transform the approximation back, we have to apply the transformation $\Rho$. We give some explicit densities and the corresponding transformations in Example~\ref{ex:rhos} for $\Omega = \R^d$ and in Example~\ref{ex:beta} for $\Omega = [0,1]^d$.

Our procedure coincides with transforming the function $f$ to the function $f\circ \Rho^{-1}$, which is a function on the torus.
In approximation theory it is known that 
the error decay gains from smoothness of the function. We will introduce weighted function spaces of mixed Sobolev smoothness $H^m_{\mix}(\Omega,\rho)$ in Definition~\ref{def:Hm_Omega} and even for non-periodic functions in Definition~\ref{def:Hm_cube}. For the more general function class of mixed Besov spaces we also introduce in Definition~\ref{def:besov_on_Omega} weighted spaces $\bB_{2,\infty}^s(\Omega,\rho)$ on $\Omega$. All our function space definitions rely on the definition of the periodic spaces. The relevant facts about Sobolev and Besov spaces of mixed smoothness on $\T^d$ have been collected in the appendix~\ref{sec:A1}.\\

Since we take the position that we can 
only learn the function where we have sample points, our aim is to find an approximation to the function $f$, which minimizes the $L_2$-error 
with respect to the density $\rho$. Indeed, it shows that for functions in the defined weighted function spaces $H^m_{\mix}(\Omega,\rho)$ or $\bB_{2,\infty}^s(\Omega,\rho)$ we receive in Theorem~\ref{thm:decay_trafo} the same approximation rates for the $L_2(\Omega, \rho)$-error as in the periodic setting, i.e.\,we provide that for $1/2<s<m$, $r>1$ and logarithmic oversampling $|\Y|\gtrsim r |I|\log |I|$ there is a constant $C(r,d,s)>0$ such that for fixed $\|f\|_{\bB^s_{2,\infty}(\Omega,\rho)} \leq 1$
$$\P\Big(\norm{f-(S_n^\X (f\circ \Rho^{-1}) )\circ \Rho}_{L_2(\Omega,\rho)}\leq C(r,d,s) \,\frac{(\log |\Y|)^{(d-1)(s+1/2)+s}}{|\Y|^s}\Big) \geq 1-2|\Y|^{-r}\,.$$

Wavelets have many applications in signal processing. Most commonly they are used in compression,
edge detection, noise reduction, and other signal enhancements. The broad practicality of wavelets 
is mainly due to the localization properties of wavelets in time and frequency, so that many
signals can be sparsely represented. Hence, the hyperbolic wavelet regression is a reasonable choice for  our purposes. \\

The approximation of non-periodic functions is more challenging than the periodic setting because of the boundary behaviour. 
For wavelet approximation one possibility is to construct boundary wavelets, see~\cite{DaKuUr99, Jia09}. We avoid these complicated 
constructions by extending the function, similar to Fourier extension,~\cite{AdHu20,Boyd10,Hu10}. Especially the Chui-Wang 
wavelets provide an opportunity for letting the approximation extend the function itself, so that we do not have to construct 
the extension explicitly. We suggest to choose the transformation with a fixed but small parameter $\eta$
$$\Rho(y) = \eta +(1-\eta)\int_{0}^y \rho(t) \dx t -\tfrac 12.$$
More details are described in Section~\ref{sec:extension_cube}. Our new extension method with a proper chosen extension parameter $\eta$ relies 
on the compact support of the Chui-Wang wavelets. With this approach it is possible in Corollary~\ref{cor:nonper_end} to end with nearly 
the same approximation rate as in the periodic setting.\\

In some applications we usually do not know the underlying density $\rho$ and we only get the samples $\Y$. Therefore, we first estimate the underlying density 
by $\r{\rho}$ and construct the slightly different transformation $\r{\Rho}$ in Section~\ref{sec:unknown_rho}. Using an approximation operator on $\T^d$ is also in this case 
the core idea. Naturally, the approximation error depends on the quality of the density estimation. But in Theorem~\ref{thm:bound_KDE} we 
state that we expect similar approximation results as in the case where we know the density in advance. Numerical experiments confirm this result.\\

For dealing with the curse of dimensionality we introduce the \textit{analysis of variance} (ANOVA) decomposition, see \cite{CaMoOw97,LiOw06,Holtz11}, \cite[Section 3.1.6]{NW08}, which decomposes the $d$-variate function into $2^d$ ANOVA terms $f_{\vec u}$, i.e.
$$f(\vec y) = \sum_{\vec u \subseteq \{1,\ldots,d\}}f_{\vec u}(\vec y_{\vec u}).$$
Each term corresponding to $\vec u$ only depends on variables $y_i$, where $i\in \vec u$. The number of these variables is called \textit{order} of the ANOVA term. 
However, in practical applications with high-dimensional functions, often 
only the ANOVA terms of low order play a role in order to describe the function well, see~\cite{CaMoOw97,KuSlWaWo09, DePeVo10, Wu2011, PoSc21}. For a rigorous mathematical treatment of this observation we work with functions of low effective dimension, which allow for a truncation of the hyperbolic wavelet regression. The starting point of our work is~\cite{LiPoUl21}, where the usage of the ANOVA decomposition was also beneficial to approximate periodic high-dimensional functions. \\  

Mathematical modelling of complex systems often requires sensitivity analysis to determine how an output variable of
interest is influenced by individual or subsets of input variables. A global sensitivity analysis
constitutes the study of how the output uncertainty from a
mathematical model is divvied up to distinct sources of input variation in the model. We transform the classical 
sensitivity analysis from the torus to a weighted function space. The transformation helps to tune the hyperbolic wavelet regression.
Our main suggestion is Algorithm~\ref{alg:1}, which gives is a tool for approximating high-dimensional functions from given arbitrary distributed samples from independent input variables. Furthermore, it is possible to interpret the results, since we get a knowledge about which input variables and variable interactions play a role and which do not.   \\

One main advantage of our transformation approach is that we can deal with different domains in every variable direction. 
In applications it is often the case that we have a mixture of periodic, non-periodic and real-valued input variables on the larger domain $\R$. 
Our proposed Algorithm~\ref{alg:1} is also applicable in these cases.
Furthermore we can use information of the densities, i.e.\,we handle every input variable separately, 
which enables a strategy to use density information where it is available. For a typical example see Section~\ref{sec:num_high}. \\

\textbf{Related work and other approaches}\\
\changed{We will heavily use the results of~\cite{LiPoUl21}, which gives approximation bounds and fast algorithms for the periodic setting on $\T^d$. 
In this paper we want to generalize this to a more general (tensor product) domain $\Omega$. Clearly, the main idea is the inverse transform sampling. 
But beyond that we study function spaces on $\Omega$, which provide enough smoothness of the transformed function to fulfill the assumptions for 
the periodic approximation. Further new important aspects in this paper are the extension of non-periodic functions, similar to Fourier extension and the idea of combining density estimation and the transformation. }\\

A nice introduction with a detailed description of the challenges in high-dimensional approximation is given in the book~\cite{AdBrWe22}.
The change of variables was successfully used in many applications. 

In~\cite{KaUlVo19} the authors construct a least squares approximation method for the recovery of
functions from a reproducing kernel Hilbert space on $\Omega\subset \R^d$. The key is to construct the orthonormal 
basis $(\eta_k)_{k=1}^N$ in $L_2(\Omega,\rho)$, which is in general not accessible for arbitrary or unknown densities $\rho$. 
Also the considerations~\cite{CoDaLe13,CoMi17} are based on the knowledge of the basis $\eta_k$.
With our approach we construct the concatenated functions 
$$(\eta_k(y))_{k\in I}= (\psi^\per_k(\Rho(y)))_{k\in I},$$ 
which form a semi-orthogonal 
basis in $L_2(\Omega, \rho)$. It is also possible to use other basis functions on $\T^d$ instead of the wavelet functions, but in any case the benefit is that we have the basis in $L_2(\Omega,\rho)$ available, even for a very general class of density functions. Furthermore, we are able to transform the fast algorithms from $\T^d$ to the domain $\Omega$. A recent improvement was done in~\cite{CoDo21}, where the authors use a weighted least-squares algorithm with weights related to the \textit{Christoffel function} and reduced the sampling budget by canceling the logarithmic factor. But they also assume that an orthogonal basis is known. Furthermore, in contrast to this literature, we give in Theorem~\ref{thm:decay_trafo} and Corollary~\ref{cor:nonper_end} a concentration inequality for the approximation error based on the probabilistic Bernstein inequality in comparison to estimating the expected value. 

For the examples of the Chebyshev density, which is a special case of our examples, 
\cite[Section 10.3]{KaUlVo19},\cite{CoDo21} propose the Chebyshev polynomials $\eta_{\vec k}(\vec y) 	\sim \cos(\vec k \, \arccos(\vec y))$ as basis in $L_2([-1,1]^d,\rho)$ where the inner function 
coincides with our transformation. The case that the samples are normally distributed was considered in \cite{PoSc22}. This approach coincides with our transformation. In Section~\ref{sec:other} we give more details about the connection of our weighted function spaces to those in the literature.
We study the case of fixed given samples $\Y$. In contrast to that, the task of choosing sampling points was solved successfully in~\cite{nasdaladiss} transforming rank-$1$-lattices from the torus to $\R^d$ or the cube $[0,1]^d$.  \\

\textbf{Outline}\\
This paper is organized as follows. 
In Section~\ref{sec:pre} we recall an approximation operator for periodic functions, which is based on the hyperbolic wavelet regression and the well-known ANOVA
decomposition of a function on the d-dimensional torus. Section~\ref{sec:trafo} describes the main idea of our approach, namely how we construct a transformation $\Rho$. Section~\ref{sec:weighted_function_spaces} is dedicated to the introduction of weighted function spaces. We study the spaces of mixed dominating Sobolev regularity in~\ref{sec:sobolev}, mixed dominating Besov regularity in~\ref{sec:besov} and end with defining similar spaces for non-periodic functions in~\ref{sec:cube}. \\
We study in this paper two settings: First, in Section~\ref{sec:known_rho} we assume that the underlying density is known. There we show in Theorems~\ref{thm:decay_trafo} and Corollary~\ref{cor:nonper_end} that we transfer the approximation rates from the torus to our setting on $\Omega$.
Second, we investigate in Section~\ref{sec:unknown_rho} the setting where we are given only the samples $\Y$ and no density function $\rho$. Also in this case we are able to transfer the approximation results, see Theorem~\ref{thm:bound_KDE}.
Finally, Section~\ref{sec:high-dim} is dedicated to the presentation of Algorithm~\ref{alg:1}, which gives an interpretable high-dimensional approximation. Our theoretical results are supplemented by some numerical experiments in Section~\ref{sec:num_high} that demonstrate the practical efficacy of our algorithm.

 \section{Preliminaries}\label{sec:pre}
Let us introduce the general setting and notation. 
Let $f\colon \Omega \to \C$ be a 
function on a $d$-dimensional domain $\Omega$. Given are the function values $\vec f=\left(f(\vec y)\right)_{\vec y\in \Y}$ at random 
points $\Y\subset \Omega$ with $|\Y|=M$. These samples are i.i.d.\,according to the density $\rho:\Omega\rightarrow \R$, i.e.\,$\int_{\Omega} \rho(\vec y)\d \vec y =1$. 
We will assume in this paper, that $\rho(\vec y) >0$, since we otherwise omit parts of $\Omega$ where the density is equal to zero. Furthermore, we assume that the density is continuous, sufficiently smooth and integrable. We aim to approximate the function $f$. \\

\textbf{Notation}\\
Let us introduce the weighted $L_p$-norm,
\begin{equation*}
\norm{f}_{L_p(\Omega,\rho)}:=
\begin{cases}\displaystyle\left(\int_{\Omega}|f(\vec y)|^p\rho(\vec y)\dx \vec y\right)^{\sfrac 1p} &\text{ if } p<\infty,\\
\displaystyle\sup_{\vec y\in \Omega} |f(\vec y)| &\text{ if } p = \infty.
\end{cases}
\end{equation*} 
In the case where the density $\rho$ is the uniform distribution, we use the usual notations $\norm{f}_{L_2(\Omega)}$ and $\norm{f}_{L_\infty(\Omega)}$. We focus on the case $p=2$, since in this case we have the scalar product
$$\langle f, g\rangle_\rho = \int_{\Omega} f(\vec y) g(\vec y) \rho(\vec y) \d \vec y.$$
The multi-dimensional Fourier coefficients on the torus are defined by 
\begin{equation}\label{eq:Fourier}
c_{\vec k}(f)=\int_{\T^d}f(\vec x)\,\e^{-2\pi \im\langle\vec k,\vec x\rangle}\d \vec x.
\end{equation}
This allows to write every function $f\in L_2(\T^d)$ as a Fourier series
\begin{equation*}f({\vec x})=\sum_{\vec k\in \Z^d } c_{\vec k}(f)\e^{2\pi \im\langle\vec k,{\vec x}\rangle}.
\end{equation*}
 In this paper we denote by $[d]$ the set $\{1,\ldots,d\}$. We work with a transformation idea, so we will always denote the $d$-dimensional input variable of the function $f$ in $\Omega$ by $\vec y$ and the transformed values by $\vec x\in \T^d$. The subset-vector is denoted by $\vec y_{\vec u}=(y_i)_{i\in \vec u}$ for a subset $\vec u\subseteq [d]$. 
 The complement of those subsets is always with respect to $[d]$, i.e.\,$\vec u^c=[d]\backslash \vec u$.
For an index set $\vec u\subseteq [d]$ we define the \textit{order} $|\vec u|$ as the number of elements in $\vec u$.\\
We will study the cases where 
$$\Omega=\bigtimes_{i=1}^d \Omega_i, \quad  \Omega_i\in\{\T,\R,[0,1]\} \text{ for all }  i\in [d]. $$ 
Note that a general interval $[a,b]$ with $b>a$ can be transferred to the unit interval via $y\mapsto \frac{y-a}{b-a}$. Similar to the vector notation also a subset of the domain directions is denoted by
$\Omega_{\vec u} = \times_{i\in \vec u}\Omega_i$ for $\vec u\subseteq [d]$.

\subsection{Hyperbolic Wavelet Regression on the Torus}
In this section we introduce an approximation operator for periodic functions. For a more detailed description see~\cite{LiPoUl21}. 
We introduce the notation $\psi_{j,k}(x)=2^{j /2}\psi(2^jx-k)$, for $j\in \N_0, k\in \Z$ and a wavelet function $\psi$. 
We use the periodization
$$\psi_{j,k}^\per(x)=\sum_{\ell\in \Z}\psi_{j,k}(x+\ell),$$
where $\psi^\per_{-1,0}(x)$ denotes the scalar function as well as the tensorization 
\begin{equation*}
\psi^\per_{\vec j,\vec k}(\vec x)=\prod_{i=1}^d\psi^\per_{j_i,k_i}(\vec x_i),
\end{equation*}
where $\vec j=(j_i)_{i=1}^d,\, j_i\in \{-1,0,2,\ldots\}$ and $\vec k=(k_i)_{i=1}^d$ are multi-indices $\vec k\in \I_{\vec j}$. Hence, we define the sets 
\begin{equation*}
\displaystyle
\I_{\vec j}=\bigtimes\limits_{i=1}^d \begin{cases} \{0,1,\ldots 2^{j_i}-1\}&\text{ if } j_i\geq 0,\\
\{0\} &\text{ if } j_i=-1.
\end{cases}
\end{equation*}
Furthermore, we introduce the parameter $n$, which always denotes the maximal level of the used wavelets, 
i.e.\,$\J_n=\{\vec j\in \Z^d\mid \vec j\geq -\vec 1,|\vec j|_1\leq n\}$. For notation shortening we introduce the index-set 
\begin{equation*}I_n = \{(\vec j,\vec k)\mid \vec j\in \J_n, \vec k\in \I_{\vec j}\}.
\end{equation*}
To construct an approximation operator which takes given samples $\X$ we solve the overdetermined 
system $\vec A \vec a = \vec f$, where 
\begin{equation}\label{eq:A}
\vec A=(\psi_{\vec j,\vec k}^\per(\vec x))_{\vec x\in \X,(\vec j,\vec k)\in I_n}\in \C^{M\times N}
\end{equation} 
is the \textit{hyperbolic wavelet matrix} with $M>N$. We will always denote the number of 
parameters, i.e.\,the number of columns of our wavelet matrix by $N$ with $N = |I_n|$ and the number of samples by $|\X| = M$. A detailed connection between the maximal wavelet level $n$ and the number of wavelet functions $N$ can be found in~\cite[Lemma 3.11]{LiPoUl21}. We compute the coefficients 
$\vec a$ by $\vec a = \left(\vec A^*\vec A\right)^{-1}\vec A^*\vec f.$
We will do this iteratively by minimizing the norm $\norm{\vec A\vec a -\vec f}_2$.
This gives us the wavelet coefficients of an approximation $S_n^\X f$ to $f$, i.e.
\begin{equation}\label{eq:def_S_n^X}
S_n^\X f:=\sum_{\vec j\in \J_n}\sum_{\vec k\in \I_j}a_{\vec j,\vec k}\psi_{\vec j,\vec k}^\per. 
\end{equation}
A further analysis of this operator can be found in~\cite[Corollary 3.22]{LiPoUl21}.\\
The estimates there are valid for general wavelets, which are compactly supported, i.e.
\begin{equation*}
\supp \psi =[0,2m-1],
\end{equation*} 
have vanishing moments of order $m$, i.e.
\begin{equation*}\int_{-\infty}^{\infty} \psi(x)x^\beta \d x=0,\quad \beta = 0,\ldots ,m-1,
\end{equation*}
	and the periodized wavelets form a Riesz-Basis for every index $j$ with
	\begin{equation}\label{eq:Riesz}
\gamma_m \sum_{k=0}^{2^j-1} |d_{j,k}|^2\leq \left\|\sum_{k=0}^{2^j-1}d_{j,k}\psi_{j,k}^\per(x)\right\|_{L_2(\T)}^2\leq  \delta_m \sum_{k=0}^{2^j-1} |d_{j,k}|^2.
\end{equation}
Because of this semi-orthogonality we have to use the dual basis $\psi^{\per, *}_{\vec j,\vec k}$, such that every function ${f\in L_2(\T^d)}$ can be decomposed as
\begin{equation}\label{eq:decomp_f}
f =\sum_{\vec j\geq -\vec 1}\sum_{\vec k\in \I_{\vec j}}\langle f,\psi^{\per, *}_{\vec j,\vec k} \rangle \psi^\per_{\vec j,\vec k}. 
\end{equation}
Furthermore, this decomposition gives us a connection between the wavelet coefficients, 
\begin{equation*}\langle f,\psi^{\per}_{\vec j,\vec k} \rangle = \sum_{\vec k'\in \I_{\vec j}} \langle f,\psi^{\per,*}_{\vec j,\vec k} \rangle \langle\psi^{\per}_{\vec j,\vec k},\psi^{\per}_{\vec j,\vec k'}\rangle.
\end{equation*}
In~\cite{LiPoUl21} we also introduced and analyzed the projection operator onto the wavelet space,
\begin{equation}\label{eq:defPn}
P_nf := \sum_{\vec j\in \J_n } \sum_{\vec k\in \I_{\vec j}}\langle f,\psi_{\vec j,\vec k}^{\per *} \rangle \psi_{\vec j,\vec k}^\per.
\end{equation} 
The following estimates are a short summary of the results from~\cite{LiPoUl21} for periodic functions. For the definition of the function spaces see~Appendix~\ref{sec:A1}.
\begin{align}
\norm{f-P_nf}_{L_2(\T^d)}&\lesssim 2^{-mn}n^{(d-1)/2}\norm{f}_{H^m_\mix(\T^d)}\label{eq:P_n_2}\\
\norm{f-P_nf}_{L_\infty(\T^d)} &\lesssim 2^{-n(m-1/2)}n^{(d-1)} \norm{f}_{H^m_{\mix}(\T^d)},\label{eq:P_n_infty}\\
\P\left(\norm{f-S_n^\X f}_{L_2(\T^d)}^2 \right.&\lesssim\left.   2^{-2nm}n^{d-1}\norm{f}^2_{H^m_{\mix}(\T^d)}\right)\geq 1- 2\,M^{-r},\label{eq:S_n_2}\\
\P\left(\norm{f-S_n^\X f}_{L_2(\T^d)}^2 \right.&\lesssim\left.   2^{-2ns}n^{d-1}\norm{f}^2_{\bB_{2,\infty}^s(\T^d)}\right)\geq 1- 2\,M^{-r},\quad \tfrac 12 <s<m\label{eq:S_n_3}
\end{align}
where the last result holds for some $r>1$ if we have $M \gtrsim r N \log N$ and uniformly i.i.d.\,samples $\X$. 
We will focus our numerical experiments on Chui-Wang-wavelets, where we will always denote by $m$ the \textit{order of the wavelets}, 
which denotes the number of vanishing moments of the wavelets. \\

\subsection{The ANOVA decomposition}
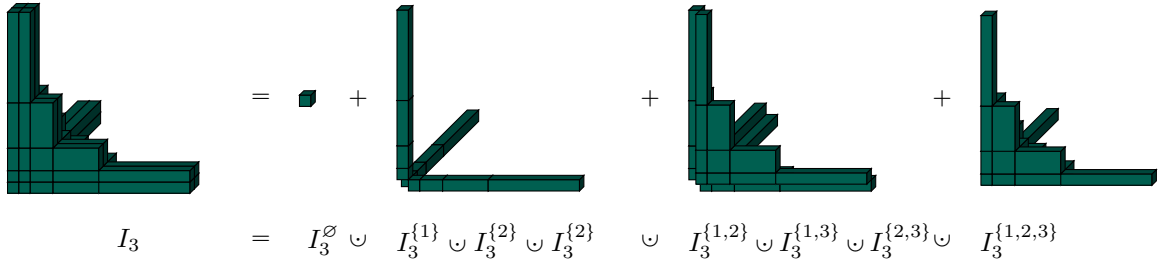
\begin{figure}[ht]
 \centering
		\begin{minipage}{.2\textwidth}
			\begin{tikzpicture}[scale=0.15]
\draw[fill=tuccolor] (-1,-1,1) -- (0,-1,1) -- (0,0,1) -- (-1,0,+1) -- cycle;
				\draw[fill=black!50!tuccolor] (0,-1,1) -- (0,-1,0) -- (0,0,0) -- (0,0,1) -- cycle;
				\draw[fill=black!25!tuccolor] (-1,0,1) -- (0,0,1) -- (0,0,0) -- (-1,0,0) -- cycle;
\foreach \y\z\color in {0/0/tuccolor,1/0/tuccolor,2/0/tuccolor,3/0/tuccolor,0/1/tuccolor,1/1/tuccolor,2/1/tuccolor,0/2/tuccolor,1/2/tuccolor,0/3/tuccolor}{
				\draw[fill=\color] (-1,2^\y-1,-2^\z+1) -- (0,2^\y-1,-2^\z+1) -- (0,2*2^\y-1,-2^\z+1) -- (-1,2*2^\y-1,-2^\z+1) -- cycle;
				\draw[fill=black!50!\color] (0,2^\y-1,-2^\z+1) -- (0,2^\y-1,-2*2^\z+1) -- (0,2*2^\y-1,-2*2^\z+1) -- (0,2*2^\y-1,-2^\z+1) -- cycle;
				\draw[fill=black!25!\color] (-1,2*2^\y-1,-2^\z+1) -- (0,2*2^\y-1,-2^\z+1) -- (0,2*2^\y-1,-2*2^\z+1) -- (-1,2*2^\y-1,-2*2^\z+1) -- cycle;
				}
\foreach \x\z\color in {0/0/tuccolor,1/0/tuccolor,2/0/tuccolor,3/0/tuccolor,0/1/tuccolor,1/1/tuccolor,2/1/tuccolor,0/2/tuccolor,1/2/tuccolor,0/3/tuccolor}{
				\draw[fill=\color] (2^\x-1,-1,-2^\z+1) -- (2*2^\x-1,-1,-2^\z+1) -- (2*2^\x-1,0,-2^\z+1) -- (2^\x-1,0,-2^\z+1) -- cycle;
				\draw[fill=black!50!\color] (2*2^\x-1,-1,-2^\z+1) -- (2*2^\x-1,-1,-2*2^\z+1) -- (2*2^\x-1,0,-2*2^\z+1) -- (2*2^\x-1,0,-2^\z+1) -- cycle;
				\draw[fill=black!25!\color] (2^\x-1,0,-2^\z+1) -- (2*2^\x-1,0,-2^\z+1) -- (2*2^\x-1,0,-2*2^\z+1) -- (2^\x-1,0,-2*2^\z+1) -- cycle;
				}
				
\foreach \y\color in {0/tuccolor,1/tuccolor,2/tuccolor,3/tuccolor}{
				\draw[fill=\color] (-1,2^\y-1,1) -- (0,2^\y-1,1) -- (0,2*2^\y-1,1) -- (-1,2*2^\y-1,1) -- cycle;
				\draw[fill=black!50!\color] (0,2^\y-1,+1) -- (0,2^\y-1,0) -- (0,2*2^\y-1,0) -- (0,2*2^\y-1,1) -- cycle;
				\draw[fill=black!25!\color] (-1,2*2^\y-1,1) -- (0,2*2^\y-1,1) -- (0,2*2^\y-1,0) -- (-1,2*2^\y-1,0) -- cycle;
				}
\foreach \x\color in {0/tuccolor,1/tuccolor,2/tuccolor,3/tuccolor}{
				\draw[fill=\color] (2^\x-1,-1,1) -- (2*2^\x-1,-1,1) -- (2*2^\x-1,0,1) -- (2^\x-1,0,1) -- cycle;
				\draw[fill=black!50!\color] (2*2^\x-1,-1,1) -- (2*2^\x-1,-1,0) -- (2*2^\x-1,0,0) -- (2*2^\x-1,0,1) -- cycle;
				\draw[fill=black!25!\color] (2^\x-1,0,1) -- (2*2^\x-1,0,1) -- (2*2^\x-1,0,0) -- (2^\x-1,0,0) -- cycle;
				}

				\foreach \x\y\z\color in {0/0/3/tuccolor,0/0/2/tuccolor,0/0/1/tuccolor,0/2/1/tuccolor,0/1/2/tuccolor,0/1/1/tuccolor,2/0/1/tuccolor,1/0/2/tuccolor,1/0/1/tuccolor,1/1/1/tuccolor,0/0/0/tuccolor,1/0/0/tuccolor,0/1/0/tuccolor,1/1/0/tuccolor,2/0/0/tuccolor,3/0/0/tuccolor,0/2/0/tuccolor,0/3/0/tuccolor,1/2/0/tuccolor,2/1/0/tuccolor}{
				\draw[fill=\color] (2^\x-1,2^\y-1,-2^\z+1) -- (2*2^\x-1,2^\y-1,-2^\z+1) -- (2*2^\x-1,2*2^\y-1,-2^\z+1) -- (2^\x-1,2*2^\y-1,-2^\z+1) -- cycle;
					\draw[fill=black!50!\color] (2*2^\x-1,2^\y-1,-2^\z+1) -- (2*2^\x-1,2^\y-1,-2*2^\z+1) -- (2*2^\x-1,2*2^\y-1,-2*2^\z+1) -- (2*2^\x-1,2*2^\y-1,-2^\z+1) -- cycle;
					\draw[fill=black!25!\color] (2^\x-1,2*2^\y-1,-2^\z+1) -- (2*2^\x-1,2*2^\y-1,-2^\z+1) -- (2*2^\x-1,2*2^\y-1,-2*2^\z+1) -- (2^\x-1,2*2^\y-1,-2*2^\z+1) -- cycle;
				}
\foreach \x\y\color in {0/0/tuccolor,1/0/tuccolor,2/0/tuccolor,3/0/tuccolor,0/1/tuccolor,1/1/tuccolor,2/1/tuccolor,0/2/tuccolor,1/2/tuccolor,0/3/tuccolor}{
				\draw[fill=\color] (2^\x-1,2^\y-1,1) -- (2*2^\x-1,2^\y-1,1) -- (2*2^\x-1,2*2^\y-1,1) -- (2^\x-1,2*2^\y-1,1) -- cycle;
				\draw[fill=black!50!\color] (2*2^\x-1,2^\y-1,1) -- (2*2^\x-1,2^\y-1,0) -- (2*2^\x-1,2*2^\y-1,-0) -- (2*2^\x-1,2*2^\y-1,1) -- cycle;
				\draw[fill=black!25!\color] (2^\x-1,2*2^\y-1,1) -- (2*2^\x-1,2*2^\y-1,1) -- (2*2^\x-1,2*2^\y-1,0) -- (2^\x-1,2*2^\y-1,0) -- cycle;
				}
			\end{tikzpicture}\end{minipage}\begin{minipage}{0.04\textwidth}
			$=$
		\end{minipage}\begin{minipage}{0.04\textwidth}
			\begin{tikzpicture}[scale=0.15]
\draw[fill=tuccolor] (-1,-1,1) -- (0,-1,1) -- (0,0,1) -- (-1,0,+1) -- cycle;
				\draw[fill=black!50!tuccolor] (0,-1,1) -- (0,-1,0) -- (0,0,0) -- (0,0,1) -- cycle;
				\draw[fill=black!25!tuccolor] (-1,0,1) -- (0,0,1) -- (0,0,0) -- (-1,0,0) -- cycle;
			\end{tikzpicture}\end{minipage}\begin{minipage}{0.04\textwidth}
			$+$
		\end{minipage}\begin{minipage}{0.2\textwidth}
			\begin{tikzpicture}[scale=0.15]
\foreach \z\color in {3/tuccolor,2/tuccolor,1/tuccolor,0/tuccolor}{
				\draw[fill=\color] (-1,-1,-2^\z+1) -- (0,-1,-2^\z+1) -- (0,0,-2^\z+1) -- (-1,0,-2^\z+1) -- cycle;
				\draw[fill=black!50!\color] (0,-1,-2^\z+1) -- (0,-1,-2*2^\z+1) -- (0,0,-2*2^\z+1) -- (0,0,-2^\z+1) -- cycle;
				\draw[fill=black!25!\color] (-1,0,-2^\z+1) -- (0,0,-2^\z+1) -- (0,0,-2*2^\z+1) -- (-1,0,-2*2^\z+1) -- cycle;
				}
\foreach \y\color in {0/tuccolor,1/tuccolor,2/tuccolor,3/tuccolor}{
				\draw[fill=\color] (-1,2^\y-1,1) -- (0,2^\y-1,1) -- (0,2*2^\y-1,1) -- (-1,2*2^\y-1,1) -- cycle;
				\draw[fill=black!50!\color] (0,2^\y-1,+1) -- (0,2^\y-1,0) -- (0,2*2^\y-1,0) -- (0,2*2^\y-1,1) -- cycle;
				\draw[fill=black!25!\color] (-1,2*2^\y-1,1) -- (0,2*2^\y-1,1) -- (0,2*2^\y-1,0) -- (-1,2*2^\y-1,0) -- cycle;
				}
\foreach \x\color in {0/tuccolor,1/tuccolor,2/tuccolor,3/tuccolor}{
				\draw[fill=\color] (2^\x-1,-1,1) -- (2*2^\x-1,-1,1) -- (2*2^\x-1,0,1) -- (2^\x-1,0,1) -- cycle;
				\draw[fill=black!50!\color] (2*2^\x-1,-1,1) -- (2*2^\x-1,-1,0) -- (2*2^\x-1,0,0) -- (2*2^\x-1,0,1) -- cycle;
				\draw[fill=black!25!\color] (2^\x-1,0,1) -- (2*2^\x-1,0,1) -- (2*2^\x-1,0,0) -- (2^\x-1,0,0) -- cycle;
				}
			\end{tikzpicture}
		\end{minipage}\begin{minipage}{0.04\textwidth}
			$+$
		\end{minipage}\begin{minipage}{0.2\textwidth}
			\begin{tikzpicture}[scale=0.15]
\foreach \y\z\color in {0/3/tuccolor,0/2/tuccolor,0/1/tuccolor,0/0/tuccolor,1/2/tuccolor,1/1/tuccolor,1/0/tuccolor,2/1/tuccolor,2/0/tuccolor,3/0/tuccolor}{
				\draw[fill=\color] (-1,2^\y-1,-2^\z+1) -- (0,2^\y-1,-2^\z+1) -- (0,2*2^\y-1,-2^\z+1) -- (-1,2*2^\y-1,-2^\z+1) -- cycle;
				\draw[fill=black!50!\color] (0,2^\y-1,-2^\z+1) -- (0,2^\y-1,-2*2^\z+1) -- (0,2*2^\y-1,-2*2^\z+1) -- (0,2*2^\y-1,-2^\z+1) -- cycle;
				\draw[fill=black!25!\color] (-1,2*2^\y-1,-2^\z+1) -- (0,2*2^\y-1,-2^\z+1) -- (0,2*2^\y-1,-2*2^\z+1) -- (-1,2*2^\y-1,-2*2^\z+1) -- cycle;
				}
\foreach \x\z\color in {0/3/tuccolor,0/2/tuccolor,1/2/tuccolor,0/1/tuccolor,1/1/tuccolor,2/1/tuccolor,0/0/tuccolor,1/0/tuccolor,2/0/tuccolor,3/0/tuccolor}{
				\draw[fill=\color] (2^\x-1,-1,-2^\z+1) -- (2*2^\x-1,-1,-2^\z+1) -- (2*2^\x-1,0,-2^\z+1) -- (2^\x-1,0,-2^\z+1) -- cycle;
				\draw[fill=black!50!\color] (2*2^\x-1,-1,-2^\z+1) -- (2*2^\x-1,-1,-2*2^\z+1) -- (2*2^\x-1,0,-2*2^\z+1) -- (2*2^\x-1,0,-2^\z+1) -- cycle;
				\draw[fill=black!25!\color] (2^\x-1,0,-2^\z+1) -- (2*2^\x-1,0,-2^\z+1) -- (2*2^\x-1,0,-2*2^\z+1) -- (2^\x-1,0,-2*2^\z+1) -- cycle;
				}
				
\foreach \x\y\color in {0/0/tuccolor,1/0/tuccolor,2/0/tuccolor,3/0/tuccolor,0/1/tuccolor,1/1/tuccolor,2/1/tuccolor,0/2/tuccolor,1/2/tuccolor,0/3/tuccolor}{
				\draw[fill=\color] (2^\x-1,2^\y-1,1) -- (2*2^\x-1,2^\y-1,1) -- (2*2^\x-1,2*2^\y-1,1) -- (2^\x-1,2*2^\y-1,1) -- cycle;
				\draw[fill=black!50!\color] (2*2^\x-1,2^\y-1,1) -- (2*2^\x-1,2^\y-1,0) -- (2*2^\x-1,2*2^\y-1,-0) -- (2*2^\x-1,2*2^\y-1,1) -- cycle;
				\draw[fill=black!25!\color] (2^\x-1,2*2^\y-1,1) -- (2*2^\x-1,2*2^\y-1,1) -- (2*2^\x-1,2*2^\y-1,0) -- (2^\x-1,2*2^\y-1,0) -- cycle;
				}
			\end{tikzpicture}\end{minipage}\begin{minipage}{0.04\textwidth}
			$+$
		\end{minipage}\begin{minipage}{0.2\textwidth}
			\begin{tikzpicture}[scale=0.15]
				\foreach \x\y\z\color in {0/0/3/tuccolor,0/0/2/tuccolor,0/0/1/tuccolor,0/2/1/tuccolor,0/1/2/tuccolor,0/1/1/tuccolor,2/0/1/tuccolor,1/0/2/tuccolor,1/0/1/tuccolor,1/1/1/tuccolor,0/0/0/tuccolor,1/0/0/tuccolor,0/1/0/tuccolor,1/1/0/tuccolor,2/0/0/tuccolor,3/0/0/tuccolor,0/2/0/tuccolor,0/3/0/tuccolor,1/2/0/tuccolor,2/1/0/tuccolor}{
				\draw[fill=\color] (2^\x-1,2^\y-1,-2^\z+1) -- (2*2^\x-1,2^\y-1,-2^\z+1) -- (2*2^\x-1,2*2^\y-1,-2^\z+1) -- (2^\x-1,2*2^\y-1,-2^\z+1) -- cycle;
				\draw[fill=black!50!\color] (2*2^\x-1,2^\y-1,-2^\z+1) -- (2*2^\x-1,2^\y-1,-2*2^\z+1) -- (2*2^\x-1,2*2^\y-1,-2*2^\z+1) -- (2*2^\x-1,2*2^\y-1,-2^\z+1) -- cycle;
				\draw[fill=black!25!\color] (2^\x-1,2*2^\y-1,-2^\z+1) -- (2*2^\x-1,2*2^\y-1,-2^\z+1) -- (2*2^\x-1,2*2^\y-1,-2*2^\z+1) -- (2^\x-1,2*2^\y-1,-2*2^\z+1) -- cycle;
				}
			\end{tikzpicture}\end{minipage}\\
		\vspace{10pt}
		\begin{minipage}{0.2\textwidth}
			\centering $I_3$
		\end{minipage}\begin{minipage}{0.04\textwidth}
			$=$
		\end{minipage}\begin{minipage}{0.04\textwidth}
			\centering
			$I_3^\varnothing$
		\end{minipage}\begin{minipage}{0.04\textwidth}
			$\cupdot$
		\end{minipage}\begin{minipage}{0.2\textwidth}
$I_3^{\{1\}}\cupdot I_3^{\{2\}} \cupdot I_3^{\{2\}}$
		\end{minipage}\begin{minipage}{0.04\textwidth}
			$\cupdot$
		\end{minipage}\begin{minipage}{0.2\textwidth}
$I_3^{\{1,2\}}\cupdot I_3^{\{1,3\}} \cupdot I_3^{\{2,3\}}$
		\end{minipage}\begin{minipage}{0.04\textwidth}
			$\cupdot$
		\end{minipage}\begin{minipage}{0.2\textwidth}
$I_3^{\{1,2,3\}}$
		\end{minipage}

\caption{Illustration of the ANOVA indices of a three-dimensional function: Every cuboid belongs to one index $\vec j$. The size represents the number of translation indices $\vec k\in \I_{\vec j}$, which gives this dyadic structure. All indices in $I_3$ are decomposed into indices belonging to the ANOVA terms with index $\vec u\subseteq [3]$. }		
\label{fig:indices_3d}	
\end{figure}

The curse of dimensionality comes into play whenever one deals with high-dimensional functions. The aim of sensitivity analysis is to describe the structure of multivariate periodic functions $f$ and to analyze the influence of each variable.
A frequently used concept is the following, \cite{CaMoOw97, Holtz11, LiOw06}.
\begin{definition}\label{def:anova-terms}
Let $f$ be in $L_2(\T^d)$. For a subset $\vec u\subseteq[d]$ we define the \textit{ANOVA (Analysis of variance) terms} by 
\begin{equation}\label{eq:anova-terms}
f_{\vec u}(\vec x_{\vec u})=\int_{\T^{d-|\vec u|}}f(\vec x)\d \vec x_{\vec u^c}-\sum_{\vec v\subset \vec u}f_{\vec v}(\vec x_{\vec v}).
\end{equation}
The\textit{ ANOVA decomposition} of a function ${f\colon \T^d\to \C}$ is then given by
\begin{equation}\label{eq:anova-decomp}
f(\vec x)=f_{\varnothing}+\sum_{i=1}^d f_{\{i\}}(x_i)+\sum_{i\neq j=1}^d f_{\{i,j\}}(x_i,x_j)+\cdots+f_{[d]}(\vec x)=\sum_{\vec u\subseteq[d]}f_{\vec u}(\vec x_{\vec u}).
\end{equation}
\end{definition} 
The terms~\eqref{eq:anova-terms} are the unique decomposition~\eqref{eq:anova-decomp}, 
such that they have mean zero and are pairwise orthogonal.
Additionally, the decomposition~\eqref{eq:decomp_f} of a function $f\in L_2(\T^d)$ in terms of wavelets,
can be written in ANOVA-terms, i.e.
\begin{equation*}
f_{\vec u}(\vec x_{\vec u}) = \sum_{(\vec j,\vec k)\in I^{\vec u}}a_{\vec j,\vec k}\psi_{\vec j,\vec k}^\per(\vec x), \quad I^{\vec u}:=\{(\vec j,\vec k)\mid \vec j_{\vec u^c} = -\vec 1,\vec k\in \I_{\vec j}\}.
\end{equation*}
The same connection holds true for a truncated wavelet decomposition with $|\vec j|_1 \leq n$
That means, all hyperbolic indices $(\vec j,\vec k)\in I_n$ can be decomposed in a disjoint union of index-sets belonging to one ANOVA term with index $\vec u\subset [d]$, i.e.
\begin{equation}\label{eq:I_n_decomp}
\displaystyle
I_n = \bigcupdot_{\vec u \subseteq[d]} I_n^{\vec u},\quad  I_n^{\vec u} = \{(\vec j,\vec k)\in I_n\mid \vec j_{\vec u^c} =-\vec 1 \}.
\end{equation}
This connection is illustrated in Figure~\ref{fig:indices_3d} for $d=3$ and $n=3$. The crucial property is that for an index $\vec u$ the corresponding functions $\psi^\per_{\vec j,\vec k}$ have to be constant in all directions $i\notin \vec u$, i.e.\,$\vec j_{\vec u^c} = -\vec 1$.
For further details see~\cite[Section 4]{LiPoUl21}.\\

 To get a notion of the importance of single terms compared to the entire function, we define the \textit{variance} of a function by 
\begin{equation*}\sigma^2(f):=\int_{\T^d}\left|f(\vec x)-\int_{\T^d}f(\vec x')\d \vec x'\right|^2\d \vec x=\int_{\T^d}|f(\vec x)|^2\d \vec x-f_\varnothing^2.
\end{equation*}
The idea of the ANOVA decomposition is to analyze which combinations of the input variables $x_j$ play a role for 
the approximation of $f$. The variances 
of the ANOVA terms indicate their importance, hence we do the following.
For subsets $\vec u\subseteq [d]$ with $\vec u\neq \varnothing$ the \textit{global sensitivity indices} (\changed{GSI})
\cite{So01} are defined as
\begin{equation}\label{eq:def_gsi} 
S(\vec u,f):=\frac{\sigma^2(f_{\vec u})}{\sigma^2(f)}\in [0,1],
\end{equation}
where the variance of the ANOVA term $f_{\vec u}$ is 
\begin{equation*}\sigma^2(f_{\vec u})= \int_{\T^{|\vec u|}}|f_{\vec u}(\vec x_{\vec u})|^2\d \vec x_{\vec u},
\end{equation*}
since the mean of the ANOVA terms is zero.
The $L_2(\T^d)$-orthogonality of the ANOVA terms implies that the variance of $f(\vec x)$ for $L_2(\T^d)$-functions $f$ can be decomposed as
$$\sigma^2(f)=\sum_{\stackrel{\vec u\subseteq [d]}{\vec u\neq \varnothing}}\sigma^2(f_{\vec u}).$$
This implies
$$\sum_{\stackrel{\vec u\subseteq [d]}{\vec u\neq \varnothing}}S(\vec u,f)=1.$$
The global sensitivity index $S(\vec u,f)$ represents the proportion of the variance of $f(\vec x)$ explained by the interaction between 
the variables indexed by $\vec u$. These indices can also be computed using only the wavelet coefficients of a function 
with the connection~\eqref{eq:I_n_decomp}.

\section{Transformations of functions to the torus}\label{sec:trafo}
In this chapter we introduce our main approach: the transformation procedure. We transform a function defined on some domain $\Omega$ to the torus, use the well-studied approximation operator for periodic functions and transform the result back to a function defined on $\Omega$.\\ 

\textbf{The univariate setting} \\
The basis for our transformation is the \textit{cumulative distribution function} $F:\Omega\to [0,1]$,
which fulfills 
\begin{equation}\label{eq:diff_eq}
\frac{\d}{\d y}F (y) = \rho(y).
\end{equation}
Subsequently, we identify $[0,1]$ with the torus by the bijective mapping $y\mapsto y-\tfrac 12$. 
Similarly to the cumulative distribution function we define the \textit{transformation} 
\begin{equation}\label{eq:Rho}
\Rho \colon \Omega \rightarrow [-\tfrac 12,\tfrac 12],\quad 
\Rho (y) := \begin{cases}
\displaystyle\int_{-\sfrac 12}^y \rho(t) \d t -\tfrac 12 &\text{if } \Omega =\T,\\
\displaystyle\int_{-\infty}^y \rho(t) \d t -\tfrac 12 &\text{if } \Omega =\R,\\
 \displaystyle\eta + \left(1-\eta\right)\int_{0}^y \rho(t) \d t -\tfrac 12,\quad 0\leq \eta\ll 1 &\text{if } \Omega =[0,1].
\end{cases}
\end{equation}
This transformation $\Rho$ gives us a possibility to transfer the function $f$ to a function $f\circ \Rho^{-1}$, 
which has its domain on the torus. Since we require that the density $\rho$ is positive, the cumulative distribution 
function is strictly monotone increasing and has a well-defined inverse function $\Rho^{-1}$.\\
In the case of a non-periodic function, 
we have to use an extension with parameter $\eta$, since otherwise the transformed function is not even continuous. Our transformation $\Rho$ transforms the function $f$ from $[0,1]$ to $[-\sfrac 12 +\eta,\sfrac 12]$, 
which means we extend the function on the boundary $[-\sfrac 12,-\sfrac 12 +\eta]$ to receive a periodic function $\tilde f\colon \T\rightarrow \C $ with 
$$\tilde f \big|_{[-\sfrac 12 +\eta,\sfrac 12]} = f\circ\Rho^{-1}.$$ 
We give more details about this extension in Section~\ref{sec:extension_cube}. The connection between the functions is illustrated here:
\begin{center}
\begin{tikzcd}[column sep=scriptsize]
\Omega \arrow[dr,"f"'] \arrow[rr, "\Rho",shift left]{}
& & \T \arrow[dl,"f\circ \Rho^{-1} \text{ or } \tilde f"]  \arrow[ll, "\Rho^{-1}",shift left]{}\\  & \C 
\end{tikzcd}
\end{center}
The transformation $\Rho$ has the property that 
\begin{equation*}
\frac{\d}{\d y}\Rho (y) 
=\begin{cases} \rho(y) &\text{ if } \Omega \in\{\T,\R\},\\
(1-\eta)\,\rho(y) &\text{ if } \Omega =[0,1].
\end{cases}
\end{equation*}
A variable substitution shows for $\Omega\in\{ \T,\R\}$ the important relation
\begin{align}\label{eq:norms_are_equal}
\norm{f}_{L_2(\Omega,\rho)}^2 
&= \int_{\Omega }|f(y)|^2\rho(y)\d y 
= \int_{\T}|f(\Rho^{-1}(x))|^2\rho(\Rho^{-1}(x))(\Rho^{-1})'(x)\d x\notag\\
&= \int_{\T}|f(\Rho^{-1}(x))|^2\rho(\Rho^{-1}(x))\tfrac{1}{\Rho'(\Rho^{-1}(x))}\d x
=\int_{\T}|f(\Rho^{-1}(x))|^2\d x
= \norm{f\circ \Rho^{-1}}_{L_2(\T)}^2.
\end{align}
In the case where $\Omega =[0,1]$ we have 
\begin{equation*}
\norm{f }^2_{L_2([0,1],\rho)}=\int_{0}^1 |f(y)|^2\rho( y)\d  y = \frac{1}{(1-\eta)}\int_{-\sfrac 12 +\eta}^{\sfrac 12} |f(\Rho^{-1}( x))|^2\d \vec x= \frac{1}{(1-\eta)}\norm{ f\circ \Rho^{-1}}^2_{L_2([-\sfrac 12 +\eta, \sfrac 12 ])}.
\end{equation*}
This relation between the $L_2$-norms motivates to transform the samples $\Y\subset \Omega$ to the transformed samples 
$\X = \Rho(\Y)$ and then use an approximation operator on $\T$. In this paper we use the operator $S_n^{\X}$, defined in~\eqref{eq:def_S_n^X}. But in general the transformation can be applied to any approximation operator on $\T$. At the end we receive the approximation  
\begin{equation}\label{eq:S_nXR}
(S_n^{\X} (f\circ \Rho^{-1}))\circ \Rho, 
\end{equation}
which takes the given sample points $\Y$ and gives back a function defined on $\Omega$. \\

In fact, we change the function and approximate the transformed function $f\circ \Rho^{-1}$, which is a function on $\T$. For the approximation operator on $\T$ it is known that, the smoother the function, the better the approximation. Indeed, it is not clear whether the transformed function $f\circ \Rho^{-1}$ inherits the smoothness of the function $f$ itself, if the density is smooth enough.    
So far, we do not know which regularity the transformed function $f\circ \Rho^{-1}$ has. 
In the following we will show that, if we request more regularity from the density $\rho$, the regularity in the Sobolev and 
Besov norm is preserved under the transformation, i.e.\,our aim is to introduce suitable 
weighted norms on $\Omega$, namely $\norm{f}_{H^s(\Omega,\rho)}$ and $\norm{f}_{\bB^s_{2,\infty}(\Omega,\rho)}$, such that
\begin{align}\norm{f\circ\Rho^{-1}}_{H^s(\T)} &\leq \norm{f}_{H^s(\Omega,\rho)}\label{eq:aim1}\\
\norm{f\circ\Rho^{-1}}_{\bB^s_{2,\infty}(\T)} &\leq \norm{f}_{\bB^s_{2,\infty}(\Omega,\rho)}.\label{eq:aim2}
\end{align}

\textbf{The multivariate setting}\\
In the multivariate setting we consider the domain $\Omega =  \times_{i=1}^d\Omega_i $ with $\Omega_i\in \{\T,\R, [0,1]\}$ for $i\in [d]$. We require that the input variables $y_i$ are independent, which means that the density $\rho(\vec y)$ is a product measure, 
\begin{equation}\label{eq:rho_d}
\rho(\vec y)=\prod_{i=1}^d \rho_i(y_i).
\end{equation}
We build up a $d$-dimensional transformation $\Rho\colon \Omega \to\T^d$ from one-dimensional transformations~\eqref{eq:Rho} by
\begin{equation}\label{eq:Rho_diff}
\Rho(\vec y) := \left(\Rho_1(y_1),\ldots,\Rho_d(y_d)\right) \quad \text{ with }\quad  
\frac{\d}{\d y_i}\Rho_i (y_i) = 
\begin{cases}
\rho_i(y_i)&\text{ if } \Omega_i \in \{\T,\R\},\\
(1-\eta)\rho_i(y_i)&\text{ if } \Omega_i = [0,1].
\end{cases}
\end{equation} 
From time to time we use the notation 
$\Rho_{\vec u} (\vec y_{\vec u})= \left(\Rho_i(y_i)\right)_{i\in \vec u} $, which is similar to the notation for vectors with index $\vec u$.
The inverse transformation is given by
\begin{equation}\label{eq:diff_Rho_-1}
\Rho^{-1}(\vec x) = \left(\Rho_1^{-1}(x_1),\ldots,\Rho^{-1}_d(x_d)\right) \quad \text{ with }\quad \frac{\d}{\d x_j}\Rho^{-1}_i (x_i) 
= \delta_{i,j}\begin{cases}
\tfrac{1}{\rho_i(\Rho_i^{-1}(x_i))} &\text{ if } \Omega_i \in \{\T,\R\},\\
\tfrac{1}{(1-\eta)\,\rho_i(\Rho_i^{-1}(x_i))} &\text{ if } \Omega_i = [0,1].
\end{cases}
\end{equation} 
The relation that we will use through this paper can be seen in this illustration:
\begin{center}
\begin{tikzcd}[column sep=scriptsize]
\text{samples: } &\vec y\in \Y   \arrow[rr, "\Rho",shift left]{} & & \vec x\in \X \arrow[ll, "\Rho^{-1}",shift left]{}\\
\text{domain: }&\Omega \arrow[dr,"f"'] \arrow[rr, "\Rho",shift left]{}
& & \T^d \arrow[dl,"f\circ \Rho^{-1}\text{ or } \tilde f"]  \arrow[ll, "\Rho^{-1}",shift left]{}\\  && \C
\end{tikzcd}
\end{center}
By the observation that the Jacobi matrix
$$(\Dx (\Rho^{-1})(\vec x))_{i,j}:=\frac{\partial }{\partial x_j}\Rho^{-1}_i(\vec x) = \delta_{i,j}\begin{cases}
\tfrac{1}{\rho_i(\Rho_i^{-1}(x_i))} &\text{ if } \Omega_i \in \{\T,\R\},\\
\tfrac{1}{(1-\eta)\,\rho_i(\Rho_i^{-1}(x_i))} &\text{ if } \Omega_i = [0,1].
\end{cases}$$
is a diagonal matrix because of the product structure~\eqref{eq:Rho_diff} of our transformation, it follows that, similar as in the univariate case~\eqref{eq:norms_are_equal},
\begin{align}\label{eq:L2_equality}
\norm{f}_{L_2(\Omega,\rho)}^2 
&= \int_{\Omega }|f(\vec y)|^2\rho(\vec y)\d \vec y 
=  \int_{\T^d}|f(\Rho^{-1}( \vec x))|^2\rho(\Rho^{-1}(\vec x)) |\det\left(\Dx (\Rho^{-1})(\vec x)\right)|\d \vec x\notag\\
&= \int_{\T^d}|f(\Rho^{-1}(\vec x))|^2\d \vec x
= \frac{1}{(1-\eta)^{\tilde{d}}}\norm{f\circ \Rho^{-1}}_{L_2(\T^d)}^2,
\end{align}
where $\tilde d=\left|\{i\in [d]\mid \Omega_i=[0,1]\}\right|$ is the dimension of the non-periodic variables of $f$. 
The norm equality~\eqref{eq:L2_equality} ensures that the transformation $\Rho$ preserves the $L_2$-norm of the function $f$, up to some factor for the non-periodic setting. \\

One main advantage of our transformation approach is that we have in addition an semi-orthogonal system on $\Omega$ with respect to $\rho$ i.e.
$$\langle \psi^\per_{\vec j,\vec k}\circ \Rho^{-1}, \psi^\per_{\vec j',\vec k'}\circ \Rho^{-1}\rangle_{\rho} = \delta_{\vec j,\vec j'}\, \langle \psi^\per_{\vec j,\vec k}\psi^\per_{\vec j,\vec k'}\rangle. $$
The next chapter is dedicated to the introduction of weighted function spaces on $\Omega$, such that also the smoothness of the function is inherited, i.e.\,we want to generalize the equations~\eqref{eq:aim1} and \eqref{eq:aim2} to the function spaces of dominating mixed derivatives. 

\subsection{The transformation \boldmath{$\Rho$} meets the ANOVA decomposition}\label{sec:ANOVA_trafo}
For an $L_2$-function on the domain $\Omega =  \times_{i=1}^d\Omega_i $ with $\Omega_i\in \{\T,\R,[0,1]\}$ for $i=1,\ldots, d$ with 
respect to the density $\rho$ it is possible to define a generalized ANOVA decomposition. We assume in this paper that the input variables $y_i$ 
are independent, which means that $\rho(\vec y)$ has a product structure~\eqref{eq:rho_d}. Hence, for an ANOVA index 
$\varnothing\neq \vec u\subset\{1,\ldots,d\}$ we define the marginal distributions 
\begin{equation*}
\rho_{\vec u} \colon \Omega_{\vec u} \rightarrow \R,\quad 
\rho_{\vec u}(\vec y_{\vec u}):= \prod_{i\in \vec u}\rho_i(y_i).
\end{equation*}
Then the ANOVA decomposition with respect to the measure $\rho$ is defined by
\begin{equation}\label{eq:ANOVA_terms_transform}
f(\vec y) = \sum_{\vec u\subseteq[d]} f_{\vec u}(\vec y_{\vec u}),
\end{equation}
where the ANOVA terms are expressed, analogously to~\eqref{eq:anova-terms} by a recursive formula 
\begin{align}\label{eq:ANOVA_rho}
f_\varnothing &= \int_{\Omega}f(\vec y)\rho(\vec y)\d \vec y\notag\\
f_{\vec u}(\vec y_{\vec u})&=\int_{\Omega_{\vec u^c}}f(\vec y)\rho_{\vec u^c}(\vec y_{\vec u^c})\d \vec y_{\vec u^c}-\sum_{\vec v\subset \vec u}f_{\vec v}(\vec y_{\vec v}),
\end{align}
see also~\cite{GiKuSl22,KuSlWaWo09,Ra14} for the case in $\R^d$.
Our main idea is to transform a function $f$ from $\Omega$ to the torus $\T^d$. Using Definition~\ref{def:anova-terms}, 
we have a decomposition of periodic functions on the torus, i.e.\,for the function $f\circ \Rho^{-1}$. If we 
transform this decomposition back to $\Omega$, we receive the decomposition~\eqref{eq:ANOVA_terms_transform}. This can be seen by the following.
\begin{lemma}
Let $\Omega_i\in \{\T,\R\}$ for all $i\in [d]$. The ANOVA terms defined in~\eqref{eq:ANOVA_rho} are the same as the transformed terms of the periodic function $f\circ\Rho^{-1}$ with the transformation defined in~\eqref{eq:Rho}, i.e.
$$ (f\circ \Rho^{-1})_{\vec u}(\Rho_{\vec u}(\vec y_{\vec u }))=f_{\vec u}(\vec y_{\vec u}).$$
\end{lemma}
\begin{proof}
We defined the ANOVA terms on $\Omega$ as well as the terms on $\T^d$ recursively, see~\eqref{eq:anova-terms} and~\eqref{eq:ANOVA_terms_transform}. Hence, we show by induction over the order $|\vec u|$ with the help of the substitution $\Rho(\vec y)=\vec x$:
\begin{align*} 
(f\circ \Rho^{-1})_{\varnothing}  &= \int_{\T^d} f(\Rho^{-1}(\vec x)) \dx \vec x =\int_{\Omega}f(\vec y)\rho(\vec y)\d \vec y = f_\varnothing ,\\
(f\circ \Rho^{-1})_{\vec u}(\Rho_{\vec u}(\vec y_{\vec u })) &= (f\circ \Rho^{-1})_{\vec u}(\vec x_{\vec u })
= \int_{\T^{|\vec u^c|}}f(\Rho^{-1}(\vec x))\d \vec x_{\vec u^c}-\sum_{\vec v\subset \vec u}(f\circ\Rho^{-1})_{\vec v}(\vec x_{\vec v})\\
& = \int_{\Omega_{\vec u^c}}f(\vec y)\rho_{\vec u^c}(\vec y_{\vec u^c})\d \vec y_{\vec u^c}-\sum_{\vec v\subset \vec u}(f\circ\Rho^{-1})_{\vec v}(\Rho_{\vec u}(\vec y_{\vec v}))\\
 &=\int_{\Omega_{\vec u^c}}f(\vec y)\rho_{\vec u^c}(\vec y_{\vec u^c})\d \vec y_{\vec u^c}-\sum_{\vec v\subset \vec u}f_{\vec v}(\vec y_{\vec v}) ={f_\vec u}(\vec y_{\vec u}).
\end{align*}
This gives the assertion.
\end{proof}
The decomposition \eqref{eq:ANOVA_terms_transform} of $f\in L_2(\Omega,\rho)$ preserves the orthogonality of the ANOVA-terms, since a simple substitution similar to~\eqref{eq:L2_equality} shows that 
\begin{align*}
\langle f_{\vec u},f_{\vec v}\rangle_{\rho} 
&=\int_{\Omega} f_{\vec u}(y_{\vec u})f_{\vec v}(\vec y_{\vec v})\rho(\vec y)\dx \vec y 
= \int_{\T^d} (f\circ \Rho^{-1})_{\vec u}(\vec x_{\vec u})(f\circ \Rho^{-1})_{\vec v}(\vec x_{\vec u}) \dx \vec x \\
&= \begin{cases}
0 &\text{ if } \vec v\neq \vec u,\\
\norm{(f\circ \Rho^{-1})_{\vec u}}_{L_2(\T^d)}^2&\text{ if } \vec v = \vec u.
\end{cases}
\end{align*}
Hence, the variance of the ANOVA-term of the transformed function $(f\circ\Rho^{-1})_{\vec u}$ is equal to the variance of $f_{\vec u}$ with respect to the density $\rho$,
\begin{equation}\label{eq:sigma_rho}
\sigma_{\rho}^2(f_{\vec u}) := \int_{\Omega_{\vec u}} |f_{\vec u}(\vec y_{\vec u})|^2\prod_{i\in \vec u}\rho_i(y_i)\d \vec y_{\vec u}.
\end{equation}
Analogously to the unweighted case~\eqref{eq:def_gsi}, we define the \textit{global sensitivity indices} for functions defined on $\Omega$ by 
\begin{equation}\label{eq:gsi_rho}
S(\vec u,f):=\frac{\sigma_{\rho}^2(f_{\vec u})}{\sigma_{\rho}^2(f)}\in [0,1].
\end{equation}

\begin{Remark}[\changed{Dependent input variables}]
\changed{
Note that there is no natural way to decompose $f$ into ANOVA terms for dependent input variables. Consider the extremal case where $y_i = y_j$ for $i\neq j$: It is not possible to say which proportion of the variance belong to $f_{i}$ or $f_{j}$. 
Problems arise when the input variables are correlated. If we integrate over some distribution, when in reality features are dependent, we 
create a new data set that deviates from the joint distribution and extrapolates to unlikely combinations of features, which can indicate unwanted variances for feature decompositions. \\
Thus, there has to be a precomputation step to avoid such dependencies. It would be possible to preprocess the given data by a PCA and a linear data transformation. Furthermore, there are approaches to generalize the ANOVA
 decomposition to dependent variables, see for example~\cite{Ho07, Ra142}. The generalized ANOVA decomposition is very difficult to estimate, and the generalization of our approach to this setting is behind the scope of our paper 
and provides an opportunity for further research. }
\end{Remark}

\section{Weighted Function spaces}\label{sec:weighted_function_spaces}
In this chapter we introduce weighted function spaces on $\Omega$, which generalize smoothness from periodic 
functions to functions defined on $\Omega$ using the transformation $\Rho$. The general idea is to study the smoothness of the concatenated function $f\circ \Rho^{-1}$, \changed{since in the periodic setting we know from~\cite{LiPoUl21} the higher the smoothness, the better approximation results using hyperbolic wavelet regression. On the torus, there are results for the Sobolev and Besov regularity, which is based on the wavelet coefficient decay in these spaces. For that reason, we study for these two cases, which functions on $\Omega$ are transformed to a smooth function on $\T^d$.}
In Section~\ref{sec:sobolev} we study the Sobolev norm~\eqref{eq:Hmnorm} or~\eqref{Fourier_char}, which requests the norm inequality~\eqref{eq:aim1} to preserve smoothness. For that reason we introduce in Definition~\ref{def:Hm_Omega} a weighted Sobolev norm on $\Omega$, which can be 
generalized to fractional smoothness by~\eqref{eq:Hsmix_rho}. In Section~\ref{sec:other} we present the relation to function spaces already known from the literature.
To preserve Besov regularity with the transformation $\Rho$, we use in Section~\ref{sec:besov} the 
characterization~\eqref{dyadic_besov} and we introduce the weighted Besov spaces in Definition~\ref{def:besov_on_Omega}, which fulfill~\eqref{eq:aim2}.\\
First, we will study in the following the case where $\Omega_i\in\{\T,\R\}$ for all $i\in [d]$. In Section~\ref{sec:cube} we then show that the non-periodic setting is similar, up to some slight modifications.

\subsection{Weighted Sobolev Spaces}
\label{sec:sobolev}
For measuring the smoothness of the transformed function $f\circ\Rho^{-1}$, we 
have to calculate the derivatives of the concatenation $f\circ \Rho^{-1}$, see Definition~\eqref{eq:Hmnorm}. 
We use the transformation~\eqref{eq:Rho} and consider in this subsection only the case $\Omega_i\in \{\T,\R\}$. The slight modification for non-periodic functions, i.e.\,$\eta >0$ is described in Section~\ref{sec:cube}.\\

\textbf{The univariate setting}\\
We use Faá di Bruno's formula, which generalizes the chain rule for $\alpha\geq 1$ to
\begin{equation*}\frac{\dx^\alpha}{\dx x^{\alpha}} f(\Rho^{-1}(x)) = \sum_{i= 1}^\alpha f^{(i)}(y) \,B_{\alpha,i}((\Rho^{-1})^{(1)}(x),(\Rho^{-1})^{(2)}(x),\ldots,(\Rho^{-1})^{(\alpha-i+1)}(x)),
\end{equation*}
where $B_{\alpha,i}$ are the Bell polynomials
\begin{equation*}B_{\alpha,i}(x_1,\ldots,x_{\alpha-i+1}) =\sum \frac{\alpha!}{j_1!j_2\cdots j_{\alpha-i+1}!}\prod_{k=1}^{\alpha-i+1}\left(\frac{x_k}{k!}\right)^{j_{k}},
\end{equation*}
where the sum is taken over all sequences $j_1,j_2,\ldots , j_{\alpha-i+1}$ of non-negative integers, such that these two conditions are satisfied:
\begin{equation*}
\sum_{k=1}^{\alpha-i+1} j_k = i, \quad
\sum_{k= 1}^{\alpha-i+1} k\cdot j_k = \alpha. 
\end{equation*}
We use the substitution $\Rho(x)=y$. Note that for the differentials of the inverse transformation $\Rho^{-1}$ holds
\begin{align*}(\Rho^{-1})^{(1)} (x)&= \frac{1}{\rho(y)},\quad   (\Rho^{-1})^{(2)} (x)= \frac{-\rho^{(1)}(y)}{\rho^{3}(y)},\quad (\Rho^{-1})^{(3)} (x)= \frac{\rho^{(2)}(y)}{\rho^{4}(y)}+ \frac{3\,(\rho^{(1)}(y))^2}{\rho^5(y)}, \\
\frac{\dx}{\dx x}\left(\frac{1}{\rho(y)}\right)^k &= k \left(\frac{1}{\rho(y)}\right)^{k-1}\,\left(-\frac{1}{\rho^2(y)}\right)\,\frac{\rho^{(1)}(y)}{\rho(y)} = -k \left(\frac{1}{\rho(y)}\right)^{k+2}\rho^{(1)}(y).\notag
\end{align*} 
Thus, every derivative of $(\Rho^{-1})^{(k)}(x)$ can be expressed by a term containing only derivatives 
of $\rho(y)$ up to order $k-1$ as well as power of $\frac{1}{\rho(y)}$ up to order $2k-1$. This allows us to shorten the notation by
\begin{equation}\label{eq:Bell_short}
B_{\alpha,i}(y) := B_{\alpha,i}((\Rho^{-1})^{(1)}(x),(\Rho^{-1})^{(2)}(x),\ldots,(\Rho^{-1})^{(\alpha-i+1)}(x)). 
\end{equation}
With this notation Faá di Bruno's formula reads as
\begin{equation*}\frac{\dx^\alpha}{\dx x^{\alpha}} f(\Rho^{-1}(x)) = \sum_{i= 1}^\alpha B_{\alpha,i}(y)\,f^{(i)}(y).
\end{equation*}
The Bell polynomial $B_{\alpha,i}$ can be expressed in terms of derivatives of $\rho(y)$ up to order $\alpha-1$ and powers of $\rho$ up to order $4\alpha-2i-1$.
We have for small indices,
\begin{align}\label{eq:small_examples}
B_{1,1}(y) &= \tfrac{1}{\rho(y)}, \,&&&&\notag\\
B_{2,1}(y) &= -\tfrac{\rho'(y)}{\rho^3(y)}, \, &B_{2,2}(y) &= \tfrac{1}{\rho^2(y)},\, &\\
B_{3,1}(y) &= -\tfrac{\rho^{(2)}(y)}{\rho^4(y)}+\tfrac{3(\rho'(y))^2}{\rho^5(y)}, \, &B_{3,2}(y) &= \tfrac{-3\rho'(y)}{\rho^4(y)},\, &B_{3,3}(y) &= \tfrac{1}{\rho^3(y)}\notag.
\end{align}
The $L_2(\T)$-norm of the derivatives of $f\circ \Rho^{-1}$ can thus be expressed as 
\begin{align*}
\left\|\frac{\dx^\alpha}{\dx x^\alpha}f\circ \Rho^{-1}(x)\right\|_{L_2(\T)}^2
&= \int_{\T}\left|\sum_{k= 1}^\alpha B_{\alpha,k}(y) \,\Dx ^k f(y)  \right|^2\d x
\leq\sum_{k= 1}^\alpha \int_{\Omega}\left|B_{\alpha,k}(y)\,\Dx^k f(y)  \right|^2\rho(y)\d y \\
&= \sum_{k= 1}^\alpha \left\|\Dx^k f(y) \right\|_{L_2(\Omega,|B_{\alpha,k}(y)|^2\rho(y))}^2. 
\end{align*}
For the Sobolev norm we have to sum over $\alpha$ and interchange the sums, which yields
\begin{align*}
\norm{f\circ\Rho^{-1}}_{H^m(\T)}^2 
&= \norm{f}_{L_2(\Omega,\rho)}^2+ \sum_{\alpha =1}^m\left\|\frac{\dx^\alpha}{\dx x^\alpha}f\circ \Rho^{-1}(x)\right\|_{L_2(\T)}^2\\
 &\leq\norm{f}_{L_2(\Omega,\rho)}^2+ \sum_{\alpha =1}^m \sum_{k= 1}^\alpha \left\|\Dx^k f(y) \right\|_{L_2(\Omega,|B_{\alpha,k}(y)|^2\rho(y))}^2\\
&=\norm{f}_{L_2(\Omega,\rho)}^2 + \sum_{k= 1}^m \sum_{\alpha =k}^{m}\left\|\Dx^k f(y) \right\|_{L_2(\Omega,|B_{\alpha,k}(y)|^2\rho(y))}^2.
\end{align*}
This motivates to generalize the Sobolev norm to functions defined on $\Omega$ by the following definition.
\begin{definition}\label{def:HmNorm_1d}
For $m\in \N$ we define the function space 
\begin{equation*}
H^m(\Omega,\rho):=\left\{f:\Omega\to \C\mid \norm{f}_{H^m(\Omega,\rho)}<\infty\right\},
\end{equation*}
where the norm is defined by
\begin{equation}\label{eq:Hmnorm1}
\norm{f}_{H^m(\Omega,\rho)}^2:=\sum_{k=0}^m \left\|\Dx ^k f  \right\|^2_{L_2(\Omega,\Upsilon_{m,k})}
\end{equation}
and the density $\Upsilon_{m,k}$ is defined by 
\begin{equation}\label{eq:Upsi}
\Upsilon_{m,k}(y) := \begin{cases}\sum_{\alpha=k}^{m} |B_{\alpha,k}(y)|^2 \rho(y) &\text{if }1\leq k\leq m,\\
																	\rho(y) &\text{if }k=0.
																	\end{cases}
\end{equation}
\end{definition}Note that we have for $m\geq2$ the useful recursion formula,
$$\Upsilon_{m,k}(y)= \Upsilon_{m-1,k}(y)+|B_{m,k}(y)|^2\rho(y).$$
We state the previous definition for the cases $1\leq m\leq 3$ explicitly:
\begin{align*}
\norm{f}^2_{H^1(\Omega,\rho)} &= \norm{f}_{L_2(\Omega,\rho)}^2 + \norm{f'}^2_{L_2(\Omega,\sfrac{1}{\rho})}\\
\norm{f}^2_{H^2(\Omega,\rho)} & =  \norm{f}_{L_2(\Omega,\rho)}^2+ \norm{f'}^2_{L_2\left(\Omega,\sfrac{1}{\rho}+\sfrac{(\rho')^2}{\rho^{5}} \right)} +  \norm{f''}^2_{L_2(\Omega,\sfrac{1}{\rho^{3}})}\\
\norm{f}^2_{H^3(\Omega,\rho)} &= \norm{f}_{L_2(\Omega,\rho)}^2+\norm{f'}^2_{L_2\left(\Omega,\sfrac{1}{\rho}+\sfrac{(\rho')^2}{\rho^{5}}+ \sfrac{(\rho'')^2}{\rho^{7}}-\sfrac{6\rho''(\rho')^2}{\rho^{8}}+\sfrac{9(\rho')^4}{\rho^{9}} \right)} \\
&\quad \quad+\norm{f''}^2_{L_2(\Omega,\sfrac{1}{\rho^{3}}+\sfrac{9(\rho')^2}{\rho^{7}})}+ \norm{f'''}^2_{L_2(\Omega,\sfrac{1}{\rho^{5}})}.
\end{align*}
This can also be interpreted as: the function $f$ cannot have a large $L_2$-norm of its derivatives up to order $m$ in areas 
where $\rho$ is small. One can not expect to capture such functions, since where the density is low, we can not approximate 
derivatives of the function $f$ well. 
Note that only in special cases, where the derivatives of the density $\rho$ and the density itself are bounded 
from below and above, this defined norm is equivalent to a norm defined using the derivatives of a function, weighted with the density $\rho$.
\begin{lemma}
Let $m\in \N$ be positive and the density $\rho$ fulfill $0<c_1\leq\norm{\rho^{(i)}}_{L_\infty(\Omega)}\leq c_2<\infty$ for $i=0,\ldots, m-1$. Then the norm in~\eqref{eq:Hmnorm1} is equivalent to the norm 
$$\norm{f}_{H^m(\Omega,\rho)}=\sum_{k=0}^m\left\|\Dx^k f\right\|_{L_2(\Omega,\rho)}.
$$
\end{lemma}
\begin{proof}
Since every derivative $\Dx^k \Rho^{-1}$ can be expressed in terms of derivatives of $\rho$ up to order $k-1$ 
and we also assume that $\rho$ itself is bounded from above and below, the terms $\Upsilon_{m,k}(y)$ can be bounded by
$$0<C\rho(y)\leq\Upsilon_{m,k}(y)\leq C'\rho(y)<\infty,$$
with some constants $0<C,C'<\infty$. This yields the assertion.
\end{proof}

\begin{Example}
\label{ex:rhos}
	We give three examples for distributions on $\Omega = \R$. We plotted the density function, the transformation $\Rho$ as well as the densities $\Upsilon_{m,k}(y)$ from~\eqref{eq:Upsi} in Figures~\ref{fig:example1}, \ref{fig:example_cauchy} and \ref{fig:example_laplace}.
	\begin{itemize}
	\item[i)]
	\textbf{Standard normal distribution on $\R$}\\
	The density
	\begin{equation}\label{eq:rho_N}
	\rho_N(y) = \tfrac{1}{\sqrt{2\pi}} \e^{-y^2/2}
	\end{equation}
	is the \textit{standard normal distribution}. 
	Because of this very smooth density we expect that this transformation passes the smoothness of $f$ to $f\circ \Rho^{-1}$. 
	The corresponding transformation is 
	$$\Rho(y) = \frac 12 \,\erf\left(\frac{y}{\sqrt 2}\right), \quad \text{where} \quad \erf(y) = \frac{2}{\sqrt \pi}\int_{y}^\infty \e^{-t^2}\d t. $$
\item[ii)]
	\textbf{Cauchy distribution on $\R$}\\
	The density
	\begin{equation}\label{eq:rho_C}
	\rho_C(y) = \frac{1}{\pi\,(1+y^2)}
	\end{equation}
	is a \textit{Cauchy distribution}.
	The corresponding transformation is 
	$$\Rho(y) =  \frac{1}{\pi }\arctan{y}. $$	
	\item[iii)]
	\textbf{Laplace distribution $\R$}\\
	The density
	\begin{equation}\label{eq:rho_L}
	\rho_L(y) = \frac{1}{8} \e^{-\frac{|y-2|}{4}}
	\end{equation}
	is a \textit{Laplace distribution}, the corresponding transformation is
	$$\Rho(y) =  \frac 12 \sgn(y-2)\left(1-\e^{-\frac{|y-2|}{4}}\right). $$

		\end{itemize}

\end{Example}
	
\newcommand\gauss[2]{1/(#2*sqrt(2*pi))*exp(-((x-#1)^2)/(2*#2^2))} \newcommand\gaussdiff{(-2*x)/(sqrt(2*pi))*exp(-((x)^2)/(2))} \newcommand\gaussdifftwo{(-2*x)^2/(sqrt(2*pi))*exp(-((x)^2)/(2))+ -2/(sqrt(2*pi))*exp(-((x)^2)/(2))}

\begin{figure}[ht]
\centering
\begin{minipage}[t]{.29\linewidth}
\begin{tikzpicture}
\begin{axis}[scale only axis,width = 0.7\textwidth,
style={
  mark=none,domain=-5:5,samples=100,smooth},
	legend =false,
] 
\addplot[mark=none,black] {\gauss{0}{1}};
\end{axis}
\end{tikzpicture}
\subcaption{Density function $\rho_N(y)$.}
\end{minipage}
\begin{minipage}[t]{.29\linewidth}
\begin{tikzpicture}
\begin{axis}[scale only axis, width = 0.7\textwidth, style={
  mark=none,domain=-5:5,samples=100,smooth},
	legend =false 
] 
\addplot[domain = -5:5,black,smooth]coordinates {
(-5.0,-0.49999971334842813)
(-4.979959919839679,-0.49999968201278555)
(-4.959919839679359,-0.4999996473886275)
(-4.939879759519038,-0.4999996091462021)
(-4.919839679358717,-0.4999995669243552)
(-4.8997995991983965,-0.49999952032771366)
(-4.8797595190380765,-0.499999468923633)
(-4.859719438877756,-0.4999994122388918)
(-4.839679358717435,-0.49999934975611376)
(-4.819639278557114,-0.4999992809098963)
(-4.799599198396794,-0.4999992050826248)
(-4.779559118236473,-0.4999991215999482)
(-4.759519038076152,-0.4999990297258918)
(-4.739478957915832,-0.4999989286575804)
(-4.719438877755511,-0.4999988175195441)
(-4.69939879759519,-0.4999986953575765)
(-4.679358717434869,-0.4999985611321136)
(-4.659318637274549,-0.4999984137111001)
(-4.6392785571142285,-0.49999825186230673)
(-4.619238476953908,-0.49999807424506126)
(-4.599198396793587,-0.49999787940135315)
(-4.579158316633267,-0.4999976657462685)
(-4.559118236472946,-0.4999974315577116)
(-4.539078156312625,-0.49999717496536455)
(-4.519038076152305,-0.4999968939388357)
(-4.498997995991984,-0.49999658627494326)
(-4.478957915831663,-0.4999962495840792)
(-4.458917835671342,-0.49999588127559313)
(-4.438877755511022,-0.49999547854213655)
(-4.4188376753507015,-0.49999503834290016)
(-4.398797595190381,-0.4999945573856772)
(-4.37875751503006,-0.499994032107681)
(-4.35871743486974,-0.4999934586550399)
(-4.338677354709419,-0.49999283286089297)
(-4.318637274549098,-0.4999921502220012)
(-4.298597194388778,-0.4999914058737895)
(-4.278557114228457,-0.49999059456372746)
(-4.258517034068136,-0.4999897106229558)
(-4.238476953907815,-0.4999887479360578)
(-4.218436873747495,-0.49998769990887465)
(-4.198396793587174,-0.499986559434256)
(-4.1783567134268536,-0.4999853188556348)
(-4.158316633266533,-0.4999839699283102)
(-4.138276553106213,-0.49998250377831766)
(-4.118236472945892,-0.49998091085876184)
(-4.098196392785571,-0.4999791809034818)
(-4.078156312625251,-0.4999773028779152)
(-4.05811623246493,-0.49997526492702177)
(-4.038076152304609,-0.4999730543201238)
(-4.018036072144288,-0.4999706573925153)
(-3.997995991983968,-0.499968059483687)
(-3.9779559118236474,-0.49996524487201166)
(-3.9579158316633265,-0.4999621967057278)
(-3.937875751503006,-0.4999588969300565)
(-3.917835671342685,-0.49995532621028294)
(-3.8977955911823647,-0.49995146385062783)
(-3.8777555110220443,-0.49994728770873353)
(-3.8577154308617234,-0.49994277410558263)
(-3.837675350701403,-0.4999378977306662)
(-3.817635270541082,-0.4999326315422147)
(-3.7975951903807617,-0.49992694666230103)
(-3.7775551102204408,-0.4999208122666253)
(-3.7575150300601203,-0.49991419546878574)
(-3.7374749498997994,-0.49990706119884054)
(-3.717434869739479,-0.4998993720759646)
(-3.697394789579158,-0.4998910882750021)
(-3.6773547094188377,-0.499882167386719)
(-3.6573146292585172,-0.49987256427155563)
(-3.6372745490981964,-0.49986223090668597)
(-3.617234468937876,-0.49985111622618633)
(-3.597194388777555,-0.4998391659541223)
(-3.5771543086172346,-0.4998263224303649)
(-3.5571142284569137,-0.49981252442894897)
(-3.5370741482965933,-0.499797706968795)
(-3.5170340681362724,-0.4997818011166173)
(-3.496993987975952,-0.4997647337818512)
(-3.476953907815631,-0.49974642750343695)
(-3.4569138276553106,-0.49972680022830845)
(-3.43687374749499,-0.4997057650814424)
(-3.4168336673346693,-0.4996832301273366)
(-3.396793587174349,-0.49965909812279563)
(-3.376753507014028,-0.4996332662609162)
(-3.3567134268537075,-0.49960562590617946)
(-3.3366733466933867,-0.4995760623205712)
(-3.3166332665330662,-0.4995444543806696)
(-3.2965931863727453,-0.4995106742856575)
(-3.276553106212425,-0.4994745872562356)
(-3.256513026052104,-0.49943605122443524)
(-3.2364729458917836,-0.4993949165143495)
(-3.216432865731463,-0.49935102551382915)
(-3.1963927855711423,-0.4993042123372111)
(-3.176352705410822,-0.49925430247917807)
(-3.156312625250501,-0.49920111245987303)
(-3.1362725450901805,-0.4991444494614256)
(-3.1162324649298596,-0.4990841109560768)
(-3.096192384769539,-0.4990198843261221)
(-3.0761523046092183,-0.49895154647592893)
(-3.056112224448898,-0.4988788634363186)
(-3.036072144288577,-0.4988015899616431)
(-3.0160320641282565,-0.49871946911992393)
(-2.995991983967936,-0.49863223187646255)
(-2.975951903807615,-0.498539596671373)
(-2.9559118236472948,-0.49844126899153124)
(-2.935871743486974,-0.4983369409374804)
(-2.9158316633266534,-0.49822629078587594)
(-2.8957915831663326,-0.498108982548104)
(-2.875751503006012,-0.49798466552575105)
(-2.8557114228456912,-0.49785297386365535)
(-2.835671342685371,-0.49771352610131814)
(-2.81563126252505,-0.4975659247235039)
(-2.7955911823647295,-0.49740975571091045)
(-2.775551102204409,-0.49724458809184097)
(-2.755511022044088,-0.49706997349586207)
(-2.7354709418837677,-0.49688544571048476)
(-2.715430861723447,-0.49669052024195626)
(-2.6953907815631264,-0.49648469388130373)
(-2.6753507014028055,-0.4962674442768225)
(-2.655310621242485,-0.49603822951425236)
(-2.635270541082164,-0.49579648770593576)
(-2.6152304609218437,-0.4955416365903021)
(-2.595190380761523,-0.49527307314307006)
(-2.5751503006012024,-0.49499017320160754)
(-2.555110220440882,-0.4946922911039334)
(-2.535070140280561,-0.4943787593438905)
(-2.5150300601202407,-0.4940488882440581)
(-2.49498997995992,-0.49370196564801616)
(-2.4749498997995993,-0.49333725663360367)
(-2.4549098196392785,-0.4929540032488541)
(-2.434869739478958,-0.49255142427231735)
(-2.414829659318637,-0.4921287149995061)
(-2.3947895791583167,-0.49168504705722965)
(-2.374749498997996,-0.4912195682475959)
(-2.3547094188376754,-0.49073140242348073)
(-2.3346693386773545,-0.4902196493972736)
(-2.314629258517034,-0.48968338488471475)
(-2.2945891783567136,-0.4891216604856437)
(-2.2745490981963927,-0.4885335037034729)
(-2.2545090180360723,-0.48791791800519285)
(-2.2344689378757514,-0.4872738829237014)
(-2.214428857715431,-0.4866003542042293)
(-2.19438877755511,-0.4858962639966081)
(-2.1743486973947896,-0.48516052109509356)
(-2.1543086172344688,-0.48439201122742126)
(-2.1342685370741483,-0.483589597394723)
(-2.1142284569138274,-0.4827521202638828)
(-2.094188376753507,-0.4818783986138515)
(-2.0741482965931866,-0.4809672298373726)
(-2.0541082164328657,-0.48001739049950165)
(-2.0340681362725452,-0.47902763695421724)
(-2.0140280561122244,-0.47799670602033933)
(-1.993987975951904,-0.47692331571787233)
(-1.9739478957915833,-0.47580616606579085)
(-1.9539078156312626,-0.4746439399421771)
(-1.933867735470942,-0.47343530400750244)
(-1.9138276553106213,-0.4721789096917236)
(-1.8937875751503006,-0.4708733942457329)
(-1.87374749498998,-0.4695173818575675)
(-1.8537074148296593,-0.4681094848336377)
(-1.8336673346693386,-0.4666483048450855)
(-1.813627254509018,-0.46513243423922956)
(-1.7935871743486973,-0.4635604574158901)
(-1.7735470941883769,-0.4619309522682219)
(-1.7535070140280562,-0.46024249168750997)
(-1.7334669338677355,-0.4584936451312056)
(-1.7134268537074149,-0.4566829802533)
(-1.6933867735470942,-0.4548090645959461)
(-1.6733466933867736,-0.4528704673410504)
(-1.653306613226453,-0.4508657611203649)
(-1.6332665330661322,-0.44879352388241356)
(-1.6132264529058116,-0.4466523408143931)
(-1.593186372745491,-0.44444080631698635)
(-1.5731462925851702,-0.4421575260298311)
(-1.5531062124248498,-0.43980111890518647)
(-1.5330661322645291,-0.4373702193271389)
(-1.5130260521042085,-0.43486347927349533)
(-1.4929859719438878,-0.4322795705173126)
(-1.4729458917835672,-0.4296171868648196)
(-1.4529058116232465,-0.42687504642629914)
(-1.4328657314629258,-0.42405189391630965)
(-1.4128256513026052,-0.4211465029794442)
(-1.3927855711422845,-0.418157678537651)
(-1.3727454909819639,-0.41508425915496555)
(-1.3527054108216432,-0.41192511941534427)
(-1.3326653306613228,-0.40867917230913203)
(-1.312625250501002,-0.4053453716235493)
(-1.2925851703406814,-0.40192271433244553)
(-1.2725450901803608,-0.3984102429804374)
(-1.25250501002004,-0.3948070480564311)
(-1.2324649298597194,-0.3911122703514208)
(-1.2124248496993988,-0.38732510329536185)
(-1.1923847695390781,-0.3834447952678296)
(-1.1723446893787575,-0.3794706518771087)
(-1.1523046092184368,-0.37540203820229867)
(-1.1322645290581161,-0.37123838099297785)
(-1.1122244488977955,-0.3669791708209421)
(-1.092184368737475,-0.36262396417851794)
(-1.0721442885771544,-0.358172385517953)
(-1.0521042084168337,-0.3536241292264053)
(-1.032064128256513,-0.3489789615310839)
(-1.0120240480961924,-0.34423672232914354)
(-0.9919839679358717,-0.3393973269370043)
(-0.9719438877755511,-0.3344607677538464)
(-0.9519038076152304,-0.32942711583413187)
(-0.9318637274549099,-0.32429652236411993)
(-0.9118236472945892,-0.3190692200374749)
(-0.8917835671342685,-0.313745524325216)
(-0.8717434869739479,-0.3083258346354221)
(-0.8517034068136272,-0.302810635358287)
(-0.8316633266533067,-0.29720049679231836)
(-0.811623246492986,-0.29149607594768623)
(-0.7915831663326653,-0.2856981172229544)
(-0.7715430861723447,-0.2798074529516731)
(-0.751503006012024,-0.2738250038155641)
(-0.7314629258517034,-0.26775177912130504)
(-0.7114228456913828,-0.26158887693819766)
(-0.6913827655310621,-0.2553374840943046)
(-0.6713426853707415,-0.24899887602894358)
(-0.6513026052104208,-0.24257441649974473)
(-0.6312625250501002,-0.23606555714280686)
(-0.6112224448897795,-0.22947383688482065)
(-0.591182364729459,-0.22280088120637312)
(-0.5711422845691383,-0.21604840125599648)
(-0.5511022044088176,-0.2092181928148822)
(-0.531062124248497,-0.20231213511253993)
(-0.5110220440881763,-0.19533218949404796)
(-0.4909819639278557,-0.1882803979399065)
(-0.4709418837675351,-0.18115888143987477)
(-0.45090180360721444,-0.17396983822254067)
(-0.4308617234468938,-0.1667155418427389)
(-0.41082164328657317,-0.15939833912929902)
(-0.3907815631262525,-0.15202064799596512)
(-0.37074148296593185,-0.14458495511868788)
(-0.35070140280561124,-0.13709381348283775)
(-0.3306613226452906,-0.12954983980423512)
(-0.3106212424849699,-0.12195571182822638)
(-0.2905811623246493,-0.11431416551136275)
(-0.27054108216432865,-0.10662799209055396)
(-0.250501002004008,-0.09890003504487353)
(-0.23046092184368738,-0.09113318695548371)
(-0.21042084168336672,-0.08333038626942618)
(-0.1903807615230461,-0.07549461397328867)
(-0.17034068136272545,-0.06762889018300432)
(-0.15030060120240482,-0.059736270656273203)
(-0.13026052104208416,-0.05181984323430896)
(-0.11022044088176353,-0.04388272421980971)
(-0.09018036072144289,-0.03592805469822995)
(-0.07014028056112225,-0.02795899680958836)
(-0.050100200400801605,-0.019978729978184454)
(-0.03006012024048096,-0.011990447107714831)
(-0.01002004008016032,-0.003997350749376838)
(0.01002004008016032,0.003997350749376838)
(0.03006012024048096,0.011990447107714831)
(0.050100200400801605,0.019978729978184454)
(0.07014028056112225,0.02795899680958836)
(0.09018036072144289,0.03592805469822995)
(0.11022044088176353,0.04388272421980971)
(0.13026052104208416,0.05181984323430896)
(0.15030060120240482,0.059736270656273203)
(0.17034068136272545,0.06762889018300432)
(0.1903807615230461,0.07549461397328867)
(0.21042084168336672,0.08333038626942618)
(0.23046092184368738,0.09113318695548371)
(0.250501002004008,0.09890003504487353)
(0.27054108216432865,0.10662799209055396)
(0.2905811623246493,0.11431416551136275)
(0.3106212424849699,0.12195571182822638)
(0.3306613226452906,0.12954983980423512)
(0.35070140280561124,0.13709381348283775)
(0.37074148296593185,0.14458495511868788)
(0.3907815631262525,0.15202064799596512)
(0.41082164328657317,0.15939833912929902)
(0.4308617234468938,0.1667155418427389)
(0.45090180360721444,0.17396983822254067)
(0.4709418837675351,0.18115888143987477)
(0.4909819639278557,0.1882803979399065)
(0.5110220440881763,0.19533218949404796)
(0.531062124248497,0.20231213511253993)
(0.5511022044088176,0.2092181928148822)
(0.5711422845691383,0.21604840125599648)
(0.591182364729459,0.22280088120637312)
(0.6112224448897795,0.22947383688482065)
(0.6312625250501002,0.23606555714280686)
(0.6513026052104208,0.24257441649974473)
(0.6713426853707415,0.24899887602894358)
(0.6913827655310621,0.2553374840943046)
(0.7114228456913828,0.26158887693819766)
(0.7314629258517034,0.26775177912130504)
(0.751503006012024,0.2738250038155641)
(0.7715430861723447,0.2798074529516731)
(0.7915831663326653,0.2856981172229544)
(0.811623246492986,0.29149607594768623)
(0.8316633266533067,0.29720049679231836)
(0.8517034068136272,0.302810635358287)
(0.8717434869739479,0.3083258346354221)
(0.8917835671342685,0.313745524325216)
(0.9118236472945892,0.3190692200374749)
(0.9318637274549099,0.32429652236411993)
(0.9519038076152304,0.32942711583413187)
(0.9719438877755511,0.3344607677538464)
(0.9919839679358717,0.3393973269370043)
(1.0120240480961924,0.34423672232914354)
(1.032064128256513,0.3489789615310839)
(1.0521042084168337,0.3536241292264053)
(1.0721442885771544,0.358172385517953)
(1.092184368737475,0.36262396417851794)
(1.1122244488977955,0.3669791708209421)
(1.1322645290581161,0.37123838099297785)
(1.1523046092184368,0.37540203820229867)
(1.1723446893787575,0.3794706518771087)
(1.1923847695390781,0.3834447952678296)
(1.2124248496993988,0.38732510329536185)
(1.2324649298597194,0.3911122703514208)
(1.25250501002004,0.3948070480564311)
(1.2725450901803608,0.3984102429804374)
(1.2925851703406814,0.40192271433244553)
(1.312625250501002,0.4053453716235493)
(1.3326653306613228,0.40867917230913203)
(1.3527054108216432,0.41192511941534427)
(1.3727454909819639,0.41508425915496555)
(1.3927855711422845,0.418157678537651)
(1.4128256513026052,0.4211465029794442)
(1.4328657314629258,0.42405189391630965)
(1.4529058116232465,0.42687504642629914)
(1.4729458917835672,0.4296171868648196)
(1.4929859719438878,0.4322795705173126)
(1.5130260521042085,0.43486347927349533)
(1.5330661322645291,0.4373702193271389)
(1.5531062124248498,0.43980111890518647)
(1.5731462925851702,0.4421575260298311)
(1.593186372745491,0.44444080631698635)
(1.6132264529058116,0.4466523408143931)
(1.6332665330661322,0.44879352388241356)
(1.653306613226453,0.4508657611203649)
(1.6733466933867736,0.4528704673410504)
(1.6933867735470942,0.4548090645959461)
(1.7134268537074149,0.4566829802533)
(1.7334669338677355,0.4584936451312056)
(1.7535070140280562,0.46024249168750997)
(1.7735470941883769,0.4619309522682219)
(1.7935871743486973,0.4635604574158901)
(1.813627254509018,0.46513243423922956)
(1.8336673346693386,0.4666483048450855)
(1.8537074148296593,0.4681094848336377)
(1.87374749498998,0.4695173818575675)
(1.8937875751503006,0.4708733942457329)
(1.9138276553106213,0.4721789096917236)
(1.933867735470942,0.47343530400750244)
(1.9539078156312626,0.4746439399421771)
(1.9739478957915833,0.47580616606579085)
(1.993987975951904,0.47692331571787233)
(2.0140280561122244,0.47799670602033933)
(2.0340681362725452,0.47902763695421724)
(2.0541082164328657,0.48001739049950165)
(2.0741482965931866,0.4809672298373726)
(2.094188376753507,0.4818783986138515)
(2.1142284569138274,0.4827521202638828)
(2.1342685370741483,0.483589597394723)
(2.1543086172344688,0.48439201122742126)
(2.1743486973947896,0.48516052109509356)
(2.19438877755511,0.4858962639966081)
(2.214428857715431,0.4866003542042293)
(2.2344689378757514,0.4872738829237014)
(2.2545090180360723,0.48791791800519285)
(2.2745490981963927,0.4885335037034729)
(2.2945891783567136,0.4891216604856437)
(2.314629258517034,0.48968338488471475)
(2.3346693386773545,0.4902196493972736)
(2.3547094188376754,0.49073140242348073)
(2.374749498997996,0.4912195682475959)
(2.3947895791583167,0.49168504705722965)
(2.414829659318637,0.4921287149995061)
(2.434869739478958,0.49255142427231735)
(2.4549098196392785,0.4929540032488541)
(2.4749498997995993,0.49333725663360367)
(2.49498997995992,0.49370196564801616)
(2.5150300601202407,0.4940488882440581)
(2.535070140280561,0.4943787593438905)
(2.555110220440882,0.4946922911039334)
(2.5751503006012024,0.49499017320160754)
(2.595190380761523,0.49527307314307006)
(2.6152304609218437,0.4955416365903021)
(2.635270541082164,0.49579648770593576)
(2.655310621242485,0.49603822951425236)
(2.6753507014028055,0.4962674442768225)
(2.6953907815631264,0.49648469388130373)
(2.715430861723447,0.49669052024195626)
(2.7354709418837677,0.49688544571048476)
(2.755511022044088,0.49706997349586207)
(2.775551102204409,0.49724458809184097)
(2.7955911823647295,0.49740975571091045)
(2.81563126252505,0.4975659247235039)
(2.835671342685371,0.49771352610131814)
(2.8557114228456912,0.49785297386365535)
(2.875751503006012,0.49798466552575105)
(2.8957915831663326,0.498108982548104)
(2.9158316633266534,0.49822629078587594)
(2.935871743486974,0.4983369409374804)
(2.9559118236472948,0.49844126899153124)
(2.975951903807615,0.498539596671373)
(2.995991983967936,0.49863223187646255)
(3.0160320641282565,0.49871946911992393)
(3.036072144288577,0.4988015899616431)
(3.056112224448898,0.4988788634363186)
(3.0761523046092183,0.49895154647592893)
(3.096192384769539,0.4990198843261221)
(3.1162324649298596,0.4990841109560768)
(3.1362725450901805,0.4991444494614256)
(3.156312625250501,0.49920111245987303)
(3.176352705410822,0.49925430247917807)
(3.1963927855711423,0.4993042123372111)
(3.216432865731463,0.49935102551382915)
(3.2364729458917836,0.4993949165143495)
(3.256513026052104,0.49943605122443524)
(3.276553106212425,0.4994745872562356)
(3.2965931863727453,0.4995106742856575)
(3.3166332665330662,0.4995444543806696)
(3.3366733466933867,0.4995760623205712)
(3.3567134268537075,0.49960562590617946)
(3.376753507014028,0.4996332662609162)
(3.396793587174349,0.49965909812279563)
(3.4168336673346693,0.4996832301273366)
(3.43687374749499,0.4997057650814424)
(3.4569138276553106,0.49972680022830845)
(3.476953907815631,0.49974642750343695)
(3.496993987975952,0.4997647337818512)
(3.5170340681362724,0.4997818011166173)
(3.5370741482965933,0.499797706968795)
(3.5571142284569137,0.49981252442894897)
(3.5771543086172346,0.4998263224303649)
(3.597194388777555,0.4998391659541223)
(3.617234468937876,0.49985111622618633)
(3.6372745490981964,0.49986223090668597)
(3.6573146292585172,0.49987256427155563)
(3.6773547094188377,0.499882167386719)
(3.697394789579158,0.4998910882750021)
(3.717434869739479,0.4998993720759646)
(3.7374749498997994,0.49990706119884054)
(3.7575150300601203,0.49991419546878574)
(3.7775551102204408,0.4999208122666253)
(3.7975951903807617,0.49992694666230103)
(3.817635270541082,0.4999326315422147)
(3.837675350701403,0.4999378977306662)
(3.8577154308617234,0.49994277410558263)
(3.8777555110220443,0.49994728770873353)
(3.8977955911823647,0.49995146385062783)
(3.917835671342685,0.49995532621028294)
(3.937875751503006,0.4999588969300565)
(3.9579158316633265,0.4999621967057278)
(3.9779559118236474,0.49996524487201166)
(3.997995991983968,0.499968059483687)
(4.018036072144288,0.4999706573925153)
(4.038076152304609,0.4999730543201238)
(4.05811623246493,0.49997526492702177)
(4.078156312625251,0.4999773028779152)
(4.098196392785571,0.4999791809034818)
(4.118236472945892,0.49998091085876184)
(4.138276553106213,0.49998250377831766)
(4.158316633266533,0.4999839699283102)
(4.1783567134268536,0.4999853188556348)
(4.198396793587174,0.499986559434256)
(4.218436873747495,0.49998769990887465)
(4.238476953907815,0.4999887479360578)
(4.258517034068136,0.4999897106229558)
(4.278557114228457,0.49999059456372746)
(4.298597194388778,0.4999914058737895)
(4.318637274549098,0.4999921502220012)
(4.338677354709419,0.49999283286089297)
(4.35871743486974,0.4999934586550399)
(4.37875751503006,0.499994032107681)
(4.398797595190381,0.4999945573856772)
(4.4188376753507015,0.49999503834290016)
(4.438877755511022,0.49999547854213655)
(4.458917835671342,0.49999588127559313)
(4.478957915831663,0.4999962495840792)
(4.498997995991984,0.49999658627494326)
(4.519038076152305,0.4999968939388357)
(4.539078156312625,0.49999717496536455)
(4.559118236472946,0.4999974315577116)
(4.579158316633267,0.4999976657462685)
(4.599198396793587,0.49999787940135315)
(4.619238476953908,0.49999807424506126)
(4.6392785571142285,0.49999825186230673)
(4.659318637274549,0.4999984137111001)
(4.679358717434869,0.4999985611321136)
(4.69939879759519,0.4999986953575765)
(4.719438877755511,0.4999988175195441)
(4.739478957915832,0.4999989286575804)
(4.759519038076152,0.4999990297258918)
(4.779559118236473,0.4999991215999482)
(4.799599198396794,0.4999992050826248)
(4.819639278557114,0.4999992809098963)
(4.839679358717435,0.49999934975611376)
(4.859719438877756,0.4999994122388918)
(4.8797595190380765,0.499999468923633)
(4.8997995991983965,0.49999952032771366)
(4.919839679358717,0.4999995669243552)
(4.939879759519038,0.4999996091462021)
(4.959919839679359,0.4999996473886275)
(4.979959919839679,0.49999968201278555)
(5.0,0.49999971334842813)
 };
\end{axis}
\end{tikzpicture}
\subcaption{The transformation $\Rho(y)$.}
\end{minipage}
\begin{minipage}[t]{.29\linewidth}
\begin{tikzpicture}
\begin{semilogyaxis}[scale only axis,width = 0.7\textwidth,
style={
  mark=none,domain=-10:10,samples=100,smooth},
legend entries={$\Upsilon_{1,1}$ ,$\Upsilon_{2,1}$, $\Upsilon_{2,2}$,$\Upsilon_{3,1}$,$\Upsilon_{3,2}$,$\Upsilon_{3,3}$ },
	legend to name={legend:Y_gauss} 
] 
\addplot[mark=none,blue,domain=-1:1] {(\gauss{0}{1})^(-1)}; 
\addplot[mark=none,green,domain=-1:1] {(\gauss{0}{1})^(-1)+ (\gauss{0}{1})^(-5)*\gaussdiff*\gaussdiff};  
\addplot[mark=none,green,dashed, domain=-1:1] {(\gauss{0}{1})^(-3)};  
\addplot[mark=none,orange, domain=-1:1] {(\gauss{0}{1})^(-1) + (\gauss{0}{1})^(-5)*\gaussdiff*\gaussdiff + (\gauss{0}{1})^(-7)*\gaussdifftwo*\gaussdifftwo -6\gaussdiff*\gaussdifftwo*(\gauss{0}{1})^(-8)+9\gaussdiff^4(\gauss{0}{1})^(-9) }; 
\addplot[mark=none,orange,dashed, domain=-1:1] {(\gauss{0}{1})^(-3)+ \gaussdiff*\gaussdiff*(\gauss{0}{1})^(-7)}; 
\addplot[mark=none,orange,dotted, domain=-1:1] {(\gauss{0}{1})^(-5)}; 
\end{semilogyaxis}
\end{tikzpicture}
\subcaption{The functions $\Upsilon_{m,k}(y)$ defined in~\eqref{eq:Upsi}. }
\end{minipage}
\begin{minipage}[t]{.09\linewidth}
\vspace{-3.2cm}
\ref{legend:Y_gauss}
\end{minipage}

\caption{The standard normal distribution $\rho_N$ on $\R$.}
	\label{fig:example1}
	\end{figure}
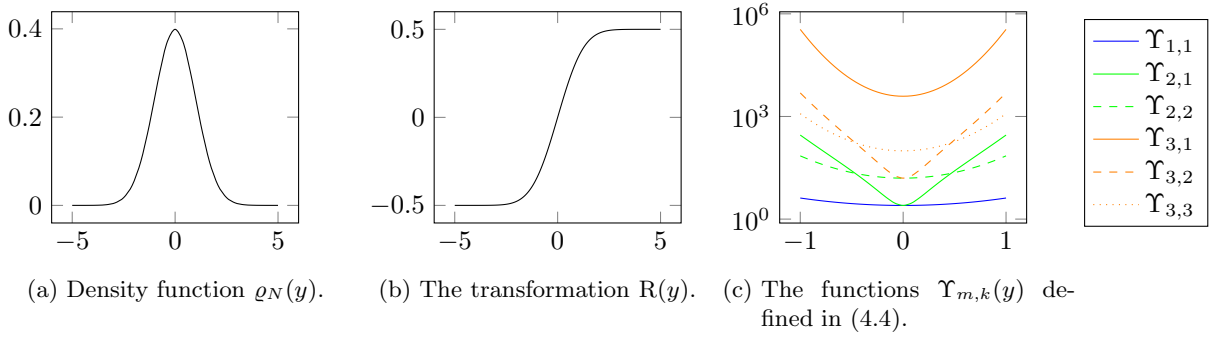 	\newcommand\cauchy{1/(pi*(1+(x)^2))} \newcommand\cauchydiff{(-2*x)/(pi*(1+(x)^2)^2)} \newcommand\cauchydifftwo{(6 (x)^2 - 2)/(pi* ((x)^2 + 1)^3)}
\newcommand\cauchyRho{1/ (pi) * rad(atan(x))} 

\begin{figure}[ht]
\centering
\begin{minipage}[t]{.29\linewidth}
\begin{tikzpicture}
\begin{axis}[scale only axis,width = 0.7\textwidth,
style={
  mark=none,domain=-5:5,samples=100,smooth},
	legend =false,
] 
\addplot[mark=none,black] {\cauchy};
\end{axis}
\end{tikzpicture}
\subcaption{Density function $\rho_C(y)$.}
\end{minipage}
\begin{minipage}[t]{.29\linewidth}
\begin{tikzpicture}
\begin{axis}[scale only axis, width = 0.7\textwidth, style={
  mark=none,domain=-5:5,samples=100,smooth},
	legend =false 
] 
\addplot[mark=none,black,domain = -5:5] {\cauchyRho};
\end{axis}
\end{tikzpicture}
\subcaption{The transformation $\Rho(y)$.}
\end{minipage}
\begin{minipage}[t]{.29\linewidth}
\begin{tikzpicture}
\begin{semilogyaxis}[scale only axis,width = 0.7\textwidth,
style={
  mark=none,domain=-10:10,samples=100,smooth},
legend entries={$\Upsilon_{1,1}$ ,$\Upsilon_{2,1}$, $\Upsilon_{2,2}$,$\Upsilon_{3,1}$,$\Upsilon_{3,2}$,$\Upsilon_{3,3}$ },
	legend to name={legend:Y_cauchy} 
] 
\addplot[mark=none,blue,domain=-1:1] {(\cauchy)^(-1)}; 
\addplot[mark=none,green,domain=-1:1] {(\cauchy)^(-1)+ (\cauchy)^(-5)*\cauchydiff*\cauchydiff};  
\addplot[mark=none,green,dashed, domain=-1:1] {(\cauchy)^(-3)};  
\addplot[mark=none,orange, domain=-1:1] {(\cauchy)^(-1) + (\cauchy)^(-5)*\cauchydiff*\cauchydiff + (\cauchy)^(-7)*\cauchydifftwo*\cauchydifftwo -6\cauchydiff*\cauchydifftwo*(\cauchy)^(-8)+9\cauchydiff^4(\cauchy)^(-9) }; 
\addplot[mark=none,orange,dashed, domain=-1:1] {(\cauchy)^(-3)+ \cauchydiff*\cauchydiff*(\cauchy)^(-7)}; 
\addplot[mark=none,orange,dotted, domain=-1:1] {(\cauchy)^(-5)}; 

\end{semilogyaxis}
\end{tikzpicture}
\subcaption{The functions $\Upsilon_{m,k}(y)$ defined in~\eqref{eq:Upsi}. }
\end{minipage}
\begin{minipage}[t]{.09\linewidth}
\vspace{-3.2cm}
\ref{legend:Y_cauchy}
\end{minipage}

\caption{The Cauchy distribution $\rho_C$ on $\R$.}
	\label{fig:example_cauchy}
	\end{figure} 	\newcommand\laplace{1/8 * exp(-abs(x-2)/4)} \newcommand\laplacediff{ - ((x - 2)* exp(- abs(x - 2)/4))/(32 *abs(x - 2))} \newcommand\laplacedifftwo{ 1/8 *(1/16 *exp(- abs(x - 2)/4)}
\newcommand\laplaceRho{(x < 2) * (1/ 2 * exp((x-2)/4)-1/2)  +(x > 2) * (1/2-1/ 2 * exp(-(x-2)/4))  } 

\begin{figure}[ht]
\centering
\begin{minipage}[t]{.29\linewidth}
\begin{tikzpicture}
\begin{axis}[scale only axis,width = 0.7\textwidth,
style={
  mark=none,domain=-2:6,samples=100,smooth},
	legend =false,
	yticklabel style={
        /pgf/number format/fixed,
        /pgf/number format/precision=2
},
 ytick distance = 0.1,
] 

\addplot[mark=none,black,domain = -2:6] {\laplace};
\end{axis}
\end{tikzpicture}
\subcaption{Density function $\rho_L(y)$.}
\end{minipage}
\begin{minipage}[t]{.29\linewidth}
\begin{tikzpicture}
\begin{axis}[scale only axis, width = 0.7\textwidth, style={
  mark=none,domain=-5:9,samples=100,smooth},
	legend =false 
] 
\addplot[mark=none,black,domain = -5:9] {\laplaceRho};

\end{axis}
\end{tikzpicture}
\subcaption{The transformation $\Rho(y)$.}
\end{minipage}
\begin{minipage}[t]{.29\linewidth}
\begin{tikzpicture}
\begin{semilogyaxis}[scale only axis,width = 0.7\textwidth,
style={
  mark=none,domain=-10:10,samples=100,smooth},
legend entries={$\Upsilon_{1,1}$ ,$\Upsilon_{2,1}$, $\Upsilon_{2,2}$,$\Upsilon_{3,1}$,$\Upsilon_{3,2}$,$\Upsilon_{3,3}$ },
	legend to name={legend:Y_laplace} 
] 
\addplot[mark=none,blue,domain=0:4] {(\laplace)^(-1)}; 
\addplot[mark=none,green,domain=0:4] {(\laplace)^(-1)+ (\laplace)^(-5)*\laplacediff*\laplacediff};  
\addplot[mark=none,green,dashed, domain=0:4] {(\laplace)^(-3)};  
\addplot[mark=none,orange, domain=0:4] {(\laplace)^(-1) + (\laplace)^(-5)*\laplacediff*\laplacediff + (\laplace)^(-7)*\laplacedifftwo*\laplacedifftwo -6\laplacediff*\laplacedifftwo*(\laplace)^(-8)+9\laplacediff^4(\laplace)^(-9) }; 
\addplot[mark=none,orange,dashed, domain=0:4] {(\laplace)^(-3)+ \laplacediff*\laplacediff*(\laplace)^(-7)}; 
\addplot[mark=none,orange,dotted, domain=0:4] {(\laplace)^(-5)}; 

\end{semilogyaxis}
\end{tikzpicture}
\subcaption{The functions $\Upsilon_{m,k}(y)$ defined in~\eqref{eq:Upsi}. }
\end{minipage}
\begin{minipage}[t]{.09\linewidth}
\vspace{-3.2cm}
\ref{legend:Y_laplace}
\end{minipage}

\caption{The Laplace distribution $\rho_L$ on $\R$.}
	\label{fig:example_laplace}
	\end{figure} 
	
So far, we characterized for natural numbers $m$ function spaces where the transformation $\Rho$ preserves the smoothness. The definition of the Sobolev norm using 
the decay of the Fourier coefficients in~\eqref{Fourier_char} allows us to study functions of fractional smoothness. 
Hence, we define fractional smoothness for functions defined on $\Omega$.
\begin{definition}
Let $s>0$. Then we define
\begin{equation*}
H^s(\Omega, \rho):=\left\{f:\Omega\to \C\mid \norm{f}_{H^s(\Omega,\rho)}<\infty\right\},
\end{equation*}
where the norm is defined by
\begin{equation*}
\norm{f}_{H^s(\Omega,\rho)}^2:=\sum_{ k\in\Z}|c_k^\rho(f)|^2(1+|k|^2)^{s}
\end{equation*}
with the Fourier coefficients $c_k^\rho(f)$ of the transformed function $c_k^\rho(f):=c_k(f\circ \Rho^{-1})$ and the Fourier coefficients for periodic functions are defined in~\eqref{eq:Fourier}. 
\end{definition}
\begin{Remark}
The norm in the previous definition is for $m=s\in \N_0$ equivalent to the norm from Definition~\ref{def:HmNorm_1d}, 
since the terms $\Upsilon_{m,k}(y)$ are chosen such that 
$$\sum_{k=1}^{m}\norm{\Dx^k (f\circ \Rho^{-1})}^2_{L_2(\T)} =\sum_{k=1}^m \left\|\Dx ^k f  \right\|^2_{L_2(\Omega,\Upsilon_{m,k})},$$ 
because of the norm equality of the norms~\eqref{eq:Hmnorm} and \eqref{Fourier_char}, see~\cite{KuSiUl14}.
\end{Remark}

\textbf{The multivariate setting}\\
The theory from the one-dimensional case can be transferred to the $d$-dimensional setting. Again, we have to use~\eqref{eq:Hmnorm} and have to estimate norms of derivatives of the transformed function $f\circ \Rho^{-1}$. 
Using the equations~\eqref{eq:diff_Rho_-1}, we have that 
\begin{equation}\label{eq:chain_d}
\Dx^{\vec \alpha}(f\circ \Rho^{-1})(\vec x) =\sum_{\vec k=\vec 1}^{\vec \alpha}B_{\vec \alpha,\vec k}(\vec y)\, \Dx^{\vec k}f (\vec y), 
\end{equation}
where we define the multivariate analogon to~\eqref{eq:Bell_short} by 
\begin{equation}\label{eq:Bell_short_multi}
B_{\vec \alpha,\vec k}(\vec y) := \prod_{i=1}^d B_{\alpha_i,k_i}(y_{k_i}).
\end{equation}
This motivates to generalize Definition~\ref{def:HmNorm_1d} to multi-variate Sobolev spaces of mixed dominating smoothness by
\begin{definition}\label{def:Hm_Omega}
For $m\in \N$ and $\Omega_i\in \{\T,\R\}$ we define the function space 
\begin{equation*}
H^m_{\mix}(\Omega,\rho):=\left\{f:\Omega\to \C\mid \norm{f}_{H^m_{\mix}(\Omega,\rho)}<\infty\right\},
\end{equation*}
where the norm is defined by
\begin{equation*}
\norm{f}_{H^m_{\mix}(\Omega,\rho)}^2=\sum_{0\leq\norm{\vec k}_{\infty}\leq m} \left\|\Dx^{\vec k} f(\vec y) \right\|^2_{L_2(\Omega,\Upsilon_{m,\vec k})}
\end{equation*}
and the density $\Upsilon_{m,\vec k}(y)$ is defined by 
$$\Upsilon_{m,\vec k}(\vec y) := \prod_{i=1}^d \Upsilon_{m,k_i}(y_{k_i}),$$
where the one-dimensional functions $\Upsilon_{m,k_i}$ are defined in~\eqref{eq:Upsi}.
\end{definition}This $H^m_{\mix}(\Omega,\rho)$-norm is equivalent to a norm definition using the decay of the Fourier coefficients of the transformed function $c^\rho_{\vec k}(f\circ\Rho^{-1})$ like in~\eqref{Fourier_char}, which can be generalized to thr case of fractional smoothness by 
\begin{equation}\label{eq:Hsmix_rho}
H^s_{\mix}(\Omega, \rho):=\left\{f:\Omega\to \C\mid \norm{f}_{H^s_{\mix}{(\Omega,\rho)}}<\infty\right\},
\end{equation}
where the norm is defined by
\begin{equation}\label{eq:Ups_multi}
\norm{f}_{H^s_{\mix}(\Omega,\rho)}^2=\sum_{\vec k\in\Z^d}|c^\rho_{\vec k}(f)|^2\prod_{i=1}^d(1+| k_i|^2)^{s}.
\end{equation}
The function spaces $H^s_{\mix}(\Omega,\rho)$ are defined such that the transformed function $f \circ \Rho^{-1}$ inherits 
the smoothness of the function $f$, i.e.
\begin{equation}\label{eq:sobolev_ineq}
\norm{f}_{H^s_{\mix}(\Omega,\rho)} \lesssim \norm{f\circ \Rho^{-1}}_{H^s_{\mix}(\T^d)}. 
\end{equation}

\subsection{Weighted Function spaces in the literature}\label{sec:other}
There is a huge literature about weighted function spaces. We restrict ourselves to a few references closely related to our approach.
Sobolev norms as defined in Definition~\ref{def:HmNorm_1d}, where the norms of the derivatives are 
measured with respect to a different density are also 
considered in~\cite{GiKuSl22} in the case $m=1$. For the one-dimensional case on $\Omega=\R$ the authors showed that the norms 
$\norm{f}_{\mathcal{H}}^2=\norm{f}_{L_2(\R,\rho)}^2 + \tfrac 1\gamma \norm{\Dx^1 f}_{L_2(\R,\psi)}^2$ and 
$\norm{f}_{\mathcal{W}}^2 = |\int_{-\infty}^\infty f(y) \rho(y)\d y|^2+\tfrac 1\gamma \norm{\Dx^1 f}_{L_2(\R,\psi)}^2 $ 
are equivalent under certain conditions on the density $\psi$. We are in the special case where $\psi(y) =\sfrac{1}{\rho(y)}$ 
and meet the conditions in our examples. The authors also showed a norm equivalence for multivariate functions.\\
Another example of weighted norms is~\cite{NuSu21}, where the authors weight the derivatives of functions with some exponential 
term in order to integrate functions on $\R$ numerically.\\

In~\cite{TriebelIII} several weighted function spaces were introduced. On $\R^d$ the authors defined the weighted counterpart to the classical Sobolev spaces by introducing a weight function $w$ for weighting all derivatives of the function $f$ similarly, i.e.\,in our notation
$$\norm{f}_{W^{m}_2(\R^d,w)} := \left(\sum_{|\vec \alpha|_1\leq m }\norm{\Dx^{\vec \alpha}f}_{L_2(\R^d,w^2)}^2\right)^{\sfrac 12}.$$
An admissible weight function $w$ has the following properties, see \cite[Definition 6.1]{TriebelIII}.
\begin{definition}
The class $W^d$ of admissible weight functions $w$ is the collection of all positive functions $w\in C^{\infty}(\R^d)$ with the properties
\begin{itemize}
	\item[i)] For all $\vec \gamma\in \N^d_0$ there is a positive constant $c_{\vec \gamma}$ with 
	\begin{equation}	\label{eq:weight1}
	|\Dx^{\vec \gamma}w(\vec y)|\leq c_{\vec\gamma}w(\vec y) \text{ for all }\vec y\in \R^d.
		\end{equation}
		\item[i)] There are two constants $c>0$ and $\alpha\geq 0$ such that 
	\begin{equation*}	0<w(\vec y_1)\leq cw(\vec y_2)\left(1+|\vec y_1-\vec y_2|^2\right)^{\sfrac{ \alpha}{2}} \text{ for all }\vec y_1\in \R^d, \vec y_2\in \R^d.
		\end{equation*}
\end{itemize}
\end{definition}
We have the following connection between these weighted spaces $W^m_2(\R^d,w)$ and the spaces $H^m(\R^d,\rho)$ defined in this paper.
\begin{theorem}\label{thm:HW}
Let $m\in \N$, $\rho\in W^d$ be an admissible weight function with $\rho(\vec y)<\infty $ for all $\vec y\in\R^d$. Then
$$W^m_2(\R^d,\rho^{-m+\sfrac12})\subset H^m(\R^d,\rho) .$$
\end{theorem}
\begin{proof}
We begin with the one-dimensional case. First we show that for the Bell polynomials defined in~\eqref{eq:Bell_short} there holds
\begin{equation*}|B_{\alpha,k}(y)|\lesssim |\rho(y)|^{-\alpha}
\end{equation*}
by induction over $\alpha$.
The condition~\eqref{eq:weight1} of an admissible weight function and the examples~\eqref{eq:small_examples} show that this is true for $\alpha\leq 3$. 
The Bell polynomials can be described recursively by, \cite{Be34},
\begin{equation}\label{eq:proof3}
B_{\alpha,k}(y) = \sum_{j=1}^{\alpha-k+1}\binom{\alpha-1}{j-1}B_{\alpha-1,k-1}((\Rho^{-1})^{(1)}(x),(\Rho^{-1})^{(2)}(x),\ldots,(\Rho^{-1})^{(\alpha-j-k+2)}(x))(\Rho^{-1})^{(j)}(\Rho(y)).
\end{equation}
The derivatives of $\Rho^{-1}$ can be bounded by
\begin{equation}\label{eq:proof2}
|(\Rho^{-1})^{(j)}(\Rho(y))|\lesssim \frac{1}{\rho(y)}, \text{ for } j\geq 1
\end{equation}
since $(\Rho^{-1})^{(1)}(\Rho(y))=\frac{1}{\rho(y)}$ and inductively every further derivative (of the sum of several 
terms which are fractions of polynomials of derivatives of $\rho$ and a power of $\rho$ in the denominator) either 
increases the power of $\rho$ in the denominator while adding a $\rho'$ in the nominator or just increases the derivatives 
of $\rho$ (but not the degree of the polynomial) in the nominator. The condition~\eqref{eq:weight1} shows than~\eqref{eq:proof2}, 
which gives by induction and using~\eqref{eq:proof3} the result~\eqref{eq:proof2}. This again gives us that 
$$|\Upsilon_{m,k}(y)|^2\lesssim\frac{1}{(\rho(y))^{2m-1}},$$
which is also true in case where $m=0$ because of the condition that $\rho(y)<\infty$.
The choice ${w=\rho^{-m+\sfrac12} \in W^d}$, gives then 
$$\norm{\Dx^{k}f}_{L_2(\R^d,w^2)}\leq\norm{\Dx^{k}f}_{L_2(\R^d,\Upsilon_{m,k})}$$
for $0\leq k\leq m$.
The multivariate case follows from the fact that $\rho$ is assumed to be a product density, such that the weight $w$ 
is also a product weight, since the densities $\Upsilon_{m,\vec k}(\vec y)$, defined in~\eqref{eq:Ups_multi} are also products of the one-dimensional functions. 
\end{proof}
Note that for instance the constant function is in all $H^m(\R^d,\rho)$, but not in $W^m(\R^d,\rho^{-m+\sfrac12})$, 
since $\lim_{x\rightarrow \pm \infty}\rho(x)=0$.
The Cauchy distribution~\eqref{eq:rho_C} belongs to the set of admissible weight functions. It is also possible to 
extend this theory of weighted function spaces to other weight functions, for instance exponential weights. Then one has 
to change the definition of admissible weights, see~\cite[Remark 6.4]{TriebelIII}, but the connection to our function spaces via Theorem~\ref{thm:HW} is nevertheless possible.

\subsection{Weighted Besov spaces}\label{sec:besov}
So far we studied the smoothness of the transformed function $f\circ\Rho$ in Sobolev spaces with dominating mixed derivatives. The advantage of Besov-Nikolskij spaces 
compared to the Sobolev spaces is that they are a much more general tool in describing the smoothness properties of functions.\\ 
Be aware that we will use in this chapter the indices $\vec j$ and $\vec k$ not as wavelet indices, but $\vec k$ as index of the Fourier coefficients and $\vec j$ as index of dyadic blocks in which we decompose the indices $\vec k$.\\
We use the Fourier analytic characterization of 
the spaces~\eqref{dyadic_besov}. Therefore we introduce the dyadic blocks 
\begin{equation*}\label{eq:def_J}
J_j=\begin{cases}
\{k\in \Z\mid 2^{j-1}\leq |k|<2^{j}\}& \text{ if } j\geq 1,\\
\{0\}& \text{ if } j=0.\end{cases}
\end{equation*}
For $\vec j \in \N_0^d$ we define 
\begin{equation*}J_{\vec j}:= J_{j_1} \times ... \times J_{j_d},
\end{equation*}
if all components belong to $\N_0$. Using these dyadic blocks, we decompose the Fourier series of the function $f$ into 
\begin{equation}\label{eq:decomp_dyadic}
f=\sum\limits_{\vec j \in \N_0^d} f_{\vec j}({\vec x}) \text{ with } f_{\vec j}({\vec x}):=\sum\limits_{\vec k \in J_{\vec j}} c_{\vec k}(f)\e^{2\pi \im\langle\vec k,{\vec x}\rangle}.
\end{equation}
Furthermore we introduce the Fourier coefficients for functions defined on $\Omega$, using the Fourier coefficients of the 
transformed function $f\circ\Rho^{-1}$ from~\eqref{eq:Fourier} by
\begin{equation*}
c^\rho_{\vec k}(f)=\int_{\T^d}(f\circ \Rho^{-1})( \vec x)\,\e^{-2\pi \im \langle \vec k,\vec x\rangle}\d  \vec x = \int_{\Omega}f(\vec y) \e^{-2\pi \im \langle\vec  k,\Rho(\vec y)\rangle}\rho(\vec y)\d \vec y.
\end{equation*}
Therefore, 
\begin{equation*}
f(\vec y) = \sum_{\vec k\in\Z^d}c_{\vec k}^\rho(f)\,\e^{2\pi \im\langle\Rho(\vec y) ,\vec k\rangle}. 
\end{equation*}
This yields immediately by~\eqref{eq:L2_equality}.
\begin{equation*}\sum_{\vec k\in \Z^d} |c_{\vec k}^\rho(f)|^2 =\norm{f}^2_{L_2(\Omega,\rho)}.
\end{equation*}
For the Definition in~\eqref{dyadic_besov} we have to split the periodic function $f\circ \Rho^{-1}$ in dyadic blocks belonging to the indices $\vec j$. For a fixed index $\vec j \in \N_0^d $ we define $\vec u = \{i\in [d]\mid j_i>0\}$. Similar to the decomposition with wavelet functions one can describe a connection of Fourier coefficients and ANOVA terms, see~\cite{PoSc19a}, for a function $g\in L_2(\T^d)$ 
\begin{equation*}c_{\vec k}(g_{\vec u}) \neq 0 \Leftrightarrow \supp \vec k := \{i\in [d]\mid k_i\neq 0\} = \vec u.
\end{equation*}
and 
\begin{equation}\label{eq:fourier_ANOVA2}
c_{\vec k}(g) = c_{\vec k_{\vec u}}(g_{\supp k}).
\end{equation}
The splitting of the function $f\circ \Rho^{-1}$ in dyadic blocks $\vec j\neq \vec 0$ gives for 
$$\vec \alpha \in \{\alpha \in \N_0^d \mid \supp \vec \alpha = \vec u ,\sum_{i\in \vec u}\alpha_i = m\},$$
using $|c_k(g)| = \tfrac{1}{k} c_k(g')$ and~\eqref{eq:chain_d} that
\begin{align*}
\allowdisplaybreaks
\norm{(f\circ \Rho^{-1})_{\vec j}}_{L_2(\T^d)}^2 
&= \sum_{\vec k\in J_{\vec j}}|c_{\vec k}(f\circ\Rho^{-1})|^2
\stackrel{\eqref{eq:fourier_ANOVA2}}{=} \sum_{\vec k\in J_{\vec j}}|c_{\vec k_{\vec u}}((f\circ\Rho^{-1})_{\vec u})|^2\\
&= \sum_{\vec k\in J_{\vec j}} \frac{1}{\prod_{i\in \vec u} |k_i|^{2\alpha_i}}\,\left|c_{\vec k_{\vec u}}\left(\frac{\d^{\vec \alpha_{\vec u}}}{\d x^{\vec \alpha_{\vec u}}}(f\circ\Rho^{-1})_{\vec u}\right)\right|^2\\
&\stackrel{\eqref{eq:chain_d}}{=} \sum_{\vec k\in J_{\vec j}} \frac{1}{\prod_{i\in \vec u} |k_i|^{2\alpha_i}}\,\left|c_{\vec k_{\vec u}}\left(\sum_{\vec \beta=\vec 1}^{\vec \alpha_{\vec u}} B_{\vec \alpha_{\vec u},\vec \beta}(\vec y_{\vec u})\,\Dx^{\vec \beta} f_{\vec u}(\vec y_{\vec u})\right)\right|^2\\
&\leq\sum_{\vec \beta=\vec 1}^{\vec \alpha_{\vec u}}\sum_{\vec k\in J_{\vec j}} \frac{1}{\prod_{i\in \vec u} |k_i|^{2\alpha_i}}\,\left|\int_{\T}B_{\vec \alpha_{\vec u},\vec \beta}(\Rho_{\vec u}^{-1}(\vec x_{\vec u}))\,\Dx^{\vec \beta} f(\Rho_{\vec u}^{-1}(\vec x_{\vec u}))\e^{-2\pi \im \langle \vec x_{\vec u} ,\vec k_{\vec u}\rangle}\d \vec x_{\vec u}\right|^2\\
&=\sum_{\vec \beta=\vec 1}^{\vec \alpha_{\vec u}} \sum_{\vec k\in J_{\vec j}} \frac{1}{\prod_{i\in \vec u} |k_i|^{2\alpha_i}}\,\left|\int_{\Omega}B_{\vec \alpha_{\vec u},\vec \beta}(\vec y_{\vec u})\,\Dx^{\vec \beta} f_{\vec u}(\vec y_{\vec u})\e^{-2\pi \im \langle\Rho_{\vec u}(\vec y_{\vec u}) ,\vec k\rangle}\rho_{\vec u}(\vec y_{\vec u})\d \vec y_{\vec u}\right|^2\\
&=\sum_{\vec \beta=\vec 1}^{\vec \alpha_{\vec u}} \sum_{\vec k\in J_{\vec j}} \frac{1}{\prod_{i\in\vec u} |k_i|^{2\alpha_i}}\,\left|c_{\vec k_{\vec u}}^{\rho_{\vec u}}\left(B_{\vec \alpha_{\vec u},\vec \beta}(\vec y_{\vec u})\,\Dx^{\vec \beta} f_{\vec u}(\vec y_{\vec u})\right)\right|^2,
\end{align*}
where we used the notation from~\eqref{eq:Bell_short_multi} and $\beta\in \N^{|\vec u|}$.
In the special case $\vec j=\vec 0$ we have $\vec k=\vec 0$ and 
$$\norm{(f\circ \Rho^{-1})_{\vec j}}_{L_2(\T^d)}^2 = \int_{\T^d}f\circ \Rho^{-1}(\vec x)\d \vec x = f_{\varnothing}$$
Therefore, we define a Besov norm for functions defined on $\Omega$ by
\begin{definition}\label{def:besov_on_Omega}
Let $s>\sfrac{1}{2}$ and $m=\left\lfloor s\right\rfloor$. Then we define 
\begin{equation*}\norm{f}_{\bB_{2,\infty}^s(\Omega,\rho)}^2 := \max\left\{f_{\varnothing},\sup_{\stackrel{\vec j\in \N^d}{\vec j\neq \vec 0}}2^{2|\vec j|_1s} \sup_{\stackrel{\vec \alpha_{\vec u}\in \N^{|\vec u|}}{|\vec \alpha_{\vec u}| = m}} \sup_{\vec 1\leq \vec \beta\leq \vec \alpha_{\vec u}} \sum_{\vec k\in J_{\vec j}}\frac{1}{\prod_{i\in \vec u} |k_i|^{2\alpha_i}}\left|c_{\vec k_{\vec u}}^{\rho_{\vec u}}\left(B_{\vec \alpha_{\vec u},\vec \beta}(\vec y_{\vec u})\,\Dx^{\vec \beta} f_{\vec u}(\vec y_{\vec u})\right)\right|^2\right\}.
\end{equation*}
\end{definition}
This definition yields the estimate
\begin{equation}\label{eq:besov_ineq}
\norm{f\circ\Rho^{-1}}_{\bB_{2,\infty}^s(\T^d)} \lesssim \norm{f}_{\bB_{2,\infty}^s(\Omega,\rho)}.
\end{equation}
In contrast to the norm $\norm{f\circ\Rho^{-1}}_{\bB_{2,\infty}^s(\T^d)}$, which deals with the transformed function on $\T^d$, the 
norm defined in Definition~\ref{def:besov_on_Omega} considers products of the function $f$ with terms consisting 
of powers and derivatives of the density $\rho$.

\subsection{A note on non-periodic functions}\label{sec:cube}
So far we studied function spaces for the case $\Omega_i\in \{\T,\R\}$ for all $i\in [d]$. A key difficulty in the approximation on $\Omega = [0,1]$ is the non-periodicity of the function, i.e.\,the behaviour at the boundary. 
No matter how the density $\rho$ looks like, a transformation which equals the cumulative distribution function can not ensure 
that the function $f\circ \Rho^{-1}$ is periodic. Therefore we introduced in the transformation~\eqref{eq:Rho} the extension parameter $\eta$. 
We denote the extension of the function $f$ by $\tilde{f}$ with $\tilde f\big|_{[-\sfrac 12+\eta,\sfrac 12]} = f\circ \Rho^{-1}$. 
The $L_2$-norm of the function $f$ itself behaves like the $L_2$-norm of the transformation up to some factor, see~\eqref{eq:L2_equality}. The same is true for the derivatives. For $\alpha\in \N$ with $\alpha \leq m$ we have
\begin{align}\label{eq:norm_eq_cube2}
\norm{\tilde{f}^{(\alpha) }}^2_{L_2([0,1])}&=\norm{\Dx^{(\alpha)}\tilde f}^2_{L_2([-\sfrac 12,-\sfrac 12+\eta])}+  
\frac{1}{(1-\eta)} \sum_{k=1}^{\alpha}\norm{\Dx^{k}f}^2_{L_2([0,1],|B_{\alpha,k}(y)|^2\rho(y))}\notag\\
&= \norm{\Dx^{(\alpha)}\tilde f}^2_{L_2([-\sfrac 12,-\sfrac 12+\eta])}+  
\sum_{k=1}^{\alpha}\norm{\Dx^{k}f}^2_{L_2([0,1],|B_{\alpha,k}(y)|^2\tfrac{\rho(y)}{(1-\eta)})}.
\end{align}
Hence, for every non-periodic direction we get the additional factor $(1-\eta)$. This gives us the reasonable assumption that extension $\tilde f$ at the boundary has to have a small Sobolev norm. For more details see Section~\ref{sec:extension_cube}. Since the factor $(1-\eta)$ is a constant for fixed $\eta$, we define for the sake of simplicity for function on $[0,1]^d$ the weighted function spaces in the same manner as for periodic functions with the following definition.  
\begin{definition}\label{def:Hm_cube}
For $m\in \N$ let the extension of $f\circ\Rho^{-1}$ to the boundary be in the Sobolev space $H^m_{\mix}(\T^d\backslash [-\sfrac 12+\eta,\sfrac 12]^d)$. Then we define the function space 
\begin{equation*}
H^m_{\mix}([0,1]^d,\rho):=\left\{f:[0,1]^d\to \C\mid \norm{f}_{H^m_{\mix}([0,1]^d,\rho)}<\infty\right\},
\end{equation*}
where the norm is defined by
\begin{equation*}
\norm{f}_{H^m_{\mix}([0,1]^d,\rho)}^2=\sum_{0\leq\norm{\vec k}_{\infty}\leq m} \left\|\Dx^{\vec k} f(\vec y) \right\|^2_{L_2([0,1]^d,\Upsilon_{m,\vec k})}
\end{equation*}
and the density $\Upsilon_{m,\vec k}(y)$ is defined by 
$$\Upsilon_{m,\vec k}(\vec y) := \prod_{i=1}^d \Upsilon_{m,k_i}(y_{k_i}),$$
where the one-dimensional functions $\Upsilon_{m,k_i}$ are defined in~\eqref{eq:Upsi}.
\end{definition}
Of course, this definition can be mixed with Definition~\ref{def:Hm_Omega} for a mixed function, which has different domains $\Omega_i$ for different $i\in[d]$. Analogously, the function spaces for fractional smoothness $s$ and for mixed Besov spaces are defined like in the periodic case.\\
If the function $f$, has certain smoothness $m$ on the interval $[0,1]$, the transformed function inherits the smoothness of this function, as in the cases for periodic functions. This means, if the Sobolev norm of the extentended function $\tilde f$ at boundary is finite, we have 
\begin{equation*}
f\in H^s_\mix([0,1]^d,\rho)\Leftrightarrow f\circ \Rho^{-1} \in H^s_{\mix}([-\tfrac 12+\eta, \tfrac 12]^d)\Leftrightarrow \tilde f \in H^s_{\mix}(\T^d).
\end{equation*}
Let us finish this excursion to the non-periodic setting with an example for the distribution $\rho$. 
\begin{Example}\label{ex:beta}
	\textbf{Beta distribution on the interval $[0,1]$}\\
	Let $\Omega =[0,1]$ be the unit interval. For $\alpha>0$ we define by
	\begin{equation}\label{eq:rho_B}
	\rho_{B,\alpha}(y) = \frac{\Gamma(2\alpha)}{\Gamma(\alpha)^2}y^{\alpha-1}(1-y)^{\alpha-1}
	\end{equation}
	the \textit{beta distribution} with the shape parameter $\alpha$, where $\Gamma$ be the Gamma function.  
	For $\alpha=1$ this is the uniform distribution. For $\alpha>1$ 
	the cumulative distribution function is the regularized incomplete beta function, so the transformation in the case $\eta = 0$ reads
	$$\Rho(y) = \frac{\Gamma(2\alpha)}{\Gamma(\alpha)^2} \int_{0}^y(t^2-t)^{\alpha-1}\d t -\tfrac 12,$$
	which can be computed analytically for fixed $\alpha$.
	These functions are plotted in Figure~\ref{fig:example2} for different parameters $\alpha$.
		\end{Example}
		\newcommand\rhocube{(1/3.14) *(x-x^2)^(-1/2)} 
\newcommand\rhocubes{6*(x-x^2)} 
\newcommand\rhocubet{30*(x-x^2)^2} 
\newcommand\rhocubediff{(1/3.14) *(-1 + 2*x)/(2 (-(-1 + x)* x)^(3/2))} 
\newcommand\rhocubedifftwo{(1/3.14)/4 *(8*x^2-8*x+3)/(x*(1-x))^(5/2)} 
\newcommand\rhocubesdiff{6*(1-2*x)} 
\newcommand\rhocubesdifftwo{-12}
\newcommand\rhocubetdiff{30 *2 *x *(1 - 3 *x + 2 *x^2)} 
\newcommand\rhocubetdifftwo{30 * 2 *(1 - 6 *x + 6 *x^2)}

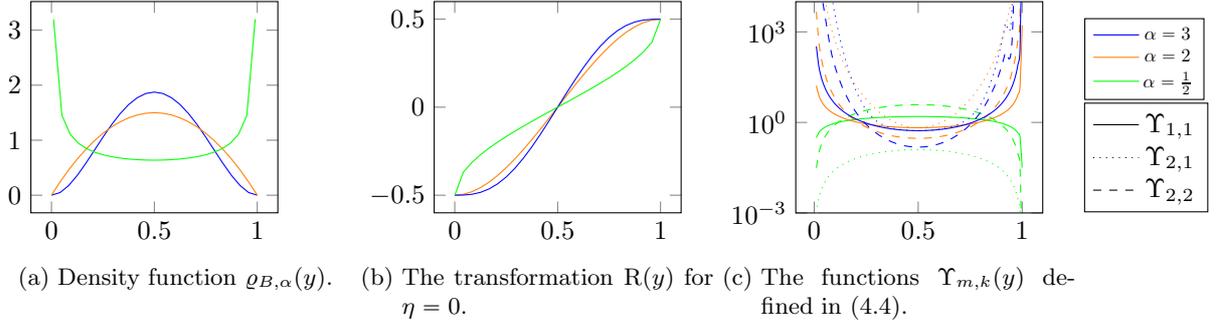
\begin{figure}[ht]
\centering
\begin{minipage}[t]{0.29\textwidth}
\begin{tikzpicture}
\begin{axis}[scale only axis,width = 0.7\textwidth,
legend style={nodes={scale=0.7, transform shape}},
style={
  mark=none,	},
	legend entries={$\alpha=3$,$\alpha=2$,$\alpha=\tfrac 12$},
	legend pos=north west,
	legend to name={legend:Y_cube1} 
] 
\addplot[mark=none,blue,domain =0:1] {30*(x-x^2)^2};
\addplot[mark=none,orange,domain =0:1] {6*(x-x^2)};
\addplot[mark=none,green,domain =0.01:0.99] {(1/3.14) *(x-x^2)^(-1/2)};\end{axis}
\end{tikzpicture}
\subcaption{Density function $\rho_{B,\alpha}(y)$.}
\end{minipage}
\begin{minipage}[t]{0.29\textwidth}
\begin{tikzpicture}
\begin{axis}[scale only axis, width = 0.7\textwidth, 
legend style={nodes={scale=0.7, transform shape}},
style={
  mark=none,},
legend=false,
	yticklabel style={
        /pgf/number format/fixed,
        /pgf/number format/precision=2
},
] 
\addplot[mark=none,orange,domain =0:1] {6*(1/2*x^2-1/3*x^3)-1/2};
\addplot[mark=none,blue,domain =0:1] {30*(x^5/5-x^4/2+x^3/3)-1/2};
\addplot[mark=none,green,domain =0:1] {-2*acos(x^(1/2))/180 +1/2};

\end{axis}
\end{tikzpicture}
\subcaption{The transformation $\Rho(y)$ for $\eta=0$.}
\end{minipage}
\begin{minipage}[t]{.29\linewidth}
\begin{tikzpicture}
\begin{semilogyaxis}[scale only axis,width = 0.7\textwidth,
style={
  mark=none,domain=0:1,samples=100,smooth},
ymax = 10^4,
	ymin = 10^(-3),
	legend entries={$\Upsilon_{1,1}$ ,$\Upsilon_{2,1}$, $\Upsilon_{2,2}$ },
	 legend image post style={black},
	legend to name={legend:Y_cube},
] 

\addplot[mark=none,green,domain=0:1,samples=100] {1/(\rhocube)};
\addplot[mark=none,green,dotted,domain=0:1,samples=100] {(\rhocube)^(-1)+ (\rhocube)^(-5) * (\rhocubediff*\rhocubediff) }; 
\addplot[mark=none,green,dashed, domain=0:1,samples=100] {(\rhocube)^(-3)};

\addplot[mark=none,orange,domain=0:1,samples=100] {(\rhocubes)^(-1)}; 
\addplot[mark=none,orange,dotted,domain=0:1,samples=100] {(\rhocubes)^(-1)+ (\rhocubes)^(-5) * (\rhocubesdiff*\rhocubesdiff) };  
\addplot[mark=none,orange,dashed, domain=0:1,samples=100] {(\rhocubes)^(-3)}; 

\addplot[mark=none,blue,domain=0:1,samples=100] {(\rhocubet)^(-1)}; 
\addplot[mark=none,blue,dotted,domain=0:1] {(\rhocubet)^(-1)+ (\rhocubet)^(-5) *(\rhocubetdiff^2) };  
\addplot[mark=none,blue,dashed, domain=0:1] {(\rhocubet)^(-3)}; 
\end{semilogyaxis}

\end{tikzpicture}

\subcaption{The functions $\Upsilon_{m,k}(y)$ defined in~\eqref{eq:Upsi}. }
\end{minipage}
\begin{minipage}[t]{.09\linewidth}
\vspace{-3.05cm}
\centering
\ref{legend:Y_cube1} 
\ref{legend:Y_cube}
\end{minipage}

\caption{The beta distribution $\rho_{B,\alpha}$ on the interval $[0,1]$ for $\alpha\in \{\sfrac 12,2,3\}$.}
	\label{fig:example2}
	\end{figure} 
\begin{Remark}
The beta distribution $\rho_{B,\sfrac 12}$ coincides with the Chebychev distribution, which is defined on $[-1,1]$. In this case \cite[Section 10.3]{KaUlVo19},\cite{CoDo21} propose the Chebyshev polynomials as basis in $L_2([-1,1]^d,\rho_{B,\sfrac 12})$. This coincides with using our transformation $\Rho$ and the cosine basis on $\T^d$.
\end{Remark}

\section{First Setting: Known density \boldmath{$\rho$}}\label{sec:known_rho}
In this chapter we study the case where we assume that the underlying density $\rho$ of the samples is known and it is a tensor product 
density~\eqref{eq:rho_d}. With this information we use the transformation~\eqref{eq:Rho}, 
transform the given samples $\Y\subset \Omega$ to the transformed samples $ \X=\Rho(\Y)$ on the torus and apply the approximation operator 
$S_n^{\X}$, given in~\eqref{eq:def_S_n^X}. With the introduction of weighted function spaces in the previous chapter, we 
estimate the error of this approximation if the function $f$ itself is in a weighted function space. In fact, we formulate the following Theorem.
\begin{theorem}\label{thm:decay_trafo}
Let $\Omega_i\in \{\T,\R\}$, the density $\rho_i$ be in $C^{m-1}(\Omega_i)$ for $i\in [d]$ and let $m\in \N$ be the order of vanishing moments of the wavelet $\psi$. Let the function fulfill for all $i\in [d]$ where $\Omega_i = \R$ that 
\begin{equation*}
\lim_{y_i\to \infty}f(\vec y) = \lim_{y_i\to -\infty}f(\vec y).
\end{equation*}
Let furthermore $M$ be the number of samples satisfying 
$M\gtrsim r N \log N$, where $N=|I_n|$ is the number of wavelet functions and $r>1$. Let $\Y = (\vec y_j)_{j=1}^M\subset \Omega$ be drawn i.i.d,.\,at random according to $\rho$, $f\in C(\Omega)$ a continuous function and the samples $\Y$ transformed to $\X = \Rho(\Y)$ using~\eqref{eq:rho_d}. In the case where $1/2<s<m$ we have 
\begin{equation}\label{eq:thm_s<m}
\P\left(\norm{f-(S_n^{\X} (f\circ \Rho^{-1}))\circ \Rho }_{L_2(\Omega,\rho)}\lesssim 2^{-ns}n^{(d-1)/2} \norm{f}_{\bB_{2,\infty}^s(\Omega, \rho)}\right)\geq 1-2\,M^{-r}.
\end{equation}
and in the case $s=m$ we have
\begin{equation}\label{eq:thm_s=m}
\P\left(\norm{f-(S_n^{\X} (f\circ \Rho^{-1}))\circ \Rho }_{L_2(\Omega,\rho)}\lesssim 2^{-ns}n^{(d-1)/2} \norm{f}_{H^s_{\mix}(\Omega, \rho)}\right)\geq 1-2\,M^{-r}.
\end{equation}
\end{theorem}
\begin{proof}
The theory from in~\cite{LiPoUl21} studies the behaviour of periodic functions on $\T^d$. Because of the assertions, the function $f\circ \Rho^{-1}$ is a function on $\T^d$, and the samples $\X=\Rho(\Y)$ are uniform i.i.d. Hence, \cite[Corollary 3.22]{LiPoUl21} is applicable to the function $f\circ \Rho^{-1}$. Together with the definitions of the Sobolev and Besov norms for functions of mixed smoothness, for which we have~\eqref{eq:sobolev_ineq} as well as~\eqref{eq:besov_ineq}, this yields the assertion. 
\end{proof}
This theorem is about the case where we do not have non-periodic variables $y_i$ involved. We introduced in~\eqref{eq:Rho} the extension parameter $\eta$ to also deal with non-periodic functions. In the next chapter we will give a similar result for the non-periodic case. Of course, these results can be mixed if the domain has different parts $\Omega_i$.

\subsection{Extensions of non-periodic functions}\label{sec:extension_cube}
Here we study the case where $\Omega = [0,1]^d$.
Of course, one can interpret a function on $[0,1]^d$ as (possibly non-continuous) function on the torus by gluing the endpoints together. 
This coincides with the transformation~\eqref{eq:Rho} with no extension $\eta=0$. Since the function $f\circ \Rho^{-1}$ is then non-continuous, 
Theorem~\ref{thm:decay_trafo} does not give us a reasonable error decay. For that reason the aim of this section is to show that the transformation idea also works for non-periodic functions with a reasonable choice of the extension parameter $\eta$ in the transformation~\eqref{eq:Rho}.

For Fourier approximation there is the approach of \textit{Fourier 
extension}~\cite{AdHu20,Boyd10,Hu10}, where the function is continued outside of the interval $[0,1]$ to 
a smooth function. See also~\cite{AdHu19} for a nice overview and the connection to the frame approach. We use a similar approach by introducing the extension parameter $\eta$, which allows us to extend $\tilde f$ on the boundary $[-\sfrac 12,-\sfrac 12+\eta]$ in an appropriate way. On one hand, this gives better approximation rates, but on the hand the stability gets worse. The occurring problem is that we have to bound the condition of the approximation matrix $\vec A$, see~\eqref{eq:A}. To circumvent this, 
in the mentioned literature the authors use for instance the truncated singular value decomposition. 
We do not want to set up the whole matrix $\vec A$, but rather use a least squares algorithm which gets only the result of a 
matrix vector multiplication with $\vec A$ and $\vec A^\top$. We will see in this chapter that an appropriate choice of the extension parameter ensures stability. \\

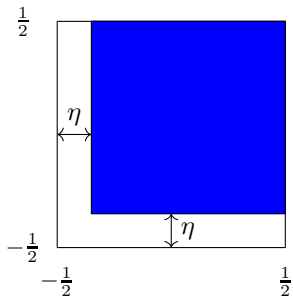
\begin{figure}[htp]
\centering
\begin{tikzpicture}[scale=3]
\draw[black] (0,0) -- (1,0);
\draw[black] (0,0) -- (0,1);
\draw[black] (0,1) -- (1,1);
\draw[black] (1,0) -- (1,1);
\draw[blue] (0.15,0.15) -- (1,0.15);
\draw[blue] (0.15,0.15) -- (0.15,1);
\draw[blue] (0.15,1) -- (1,1);
\draw[blue] (1,0.15) -- (1,1);
\node[scale=0.9] at (0,-0.15) {$-\tfrac 12$};
\node[scale=0.9] at (1,-0.15) {$\tfrac 12$};
\node[scale=0.9] at (-0.15,0) {$-\tfrac 12$};
\node[scale=0.9] at (-0.15,1) {$\tfrac 12$};
\node[scale=0.9,blue] at (0.575,0.575) {$\tilde\Omega$};
\filldraw[fill=blue,opacity=0.2] (0.15,0.15) rectangle (1,1);
\draw[<->] (0.15,0.5) --node[above] {$\eta$} (0,0.5);
\draw[<->] (0.5,0.15) --node[right] {$\eta$} (0.5,0);
\end{tikzpicture}
\caption{Illustration of the two-dimensional extension}
\label{fig:domain_2d}
\end{figure}

Remember, rather than taking the cumulative distribution function of the density $\rho$, we use the modification 
\begin{equation}\label{eq:Rho_1}
\Rho(y) = \eta + (1-\eta) \int_{0}^y \rho(t)\d t,
\end{equation}
where $0<\eta\ll 1$ is some extension parameter. In fact, we get the transformed samples $\X=\Rho(\Y)$, which are uniformly 
distributed on the cube
 $\tilde{\Omega}:=[-\sfrac 12+\eta,\sfrac 12]^d $. This procedure transforms and compresses the original 
function $f$ into the box $\tilde{\Omega}$ and allows to extend this function to a function $\tilde{f}$ defined on $\T^d$. 
Figure~\ref{fig:domain_2d} shows an illustration of the two-dimensional domains. On the boundary $[-\sfrac 12,-\sfrac 12 +\eta]$ we extend the function $f$. As mentioned in the discussion before Definition~\ref{def:Hm_cube}, it is reasonable to choose an extension which has Sobolev smoothness.
 For the following result we introduce the notation of restricting function spaces $V\subset L_2(\T^d)$ to some smaller domain by
\begin{equation*}
V\big|_{[a,b]^d} := \big\{g\big|_{[a,b]^d}\colon g\in V\big\}
\end{equation*}
for some cube $[a,b]^d\subset \T^d$.
\begin{lemma}Let $m$ be the order of the Chui-Wang wavelets and $n$ the maximal wavelet level. Denote the hyperbolic function spaces by 
$$V_n^{(d)}:=\lin\{\psi_{\vec j,\vec k}^\per\mid (\vec j,\vec k)\in I_n\}\subset L_2(\T^d).$$ 
Let the function $g\colon \tilde \Omega \rightarrow \C$ fulfill the boundary conditions
\begin{align*}g\big|_{x_i=\tfrac 12} &\in V_n^{(d-1)}\big|_{[-\tfrac 12+\eta,\tfrac 12]}, \quad g\big|_{x_i=-\tfrac 12+\eta} \in V_n^{(d-1)}\big|_{[-\tfrac 12+\eta,\tfrac 12]} \text{ for all } i\in [d],\\
 |g^{(\alpha \vec e_i)}(\vec x)|&<\infty, \text{ for all } \vec x \text{ with } x_i\in \{-\tfrac 12+\eta,\tfrac 12\}\text{ for all } \alpha=0,\ldots,m-2, \end{align*}
where $\vec e_i$ denotes the unit vector $(\vec e_i)_j = \delta_{i,j}$.
Let furthermore the extension parameter fulfill 
\begin{equation}\label{eq:condition_eta}
\eta \geq \frac{m-1}{2^{\left\lceil \tfrac nd\right\rceil +1}},
\end{equation}
where $\left\lceil a\right\rceil$ denotes the smallest integer, which is bigger than $a$. 
Then there exists an extension 
$$\tilde{g}\colon [-\tfrac 12,-\tfrac 12+\eta]^d\rightarrow \C, \quad \tilde g \in V_n^{(d)}\big|_{[-\tfrac 12,-\tfrac 12+\eta]^d}. $$

\end{lemma}
\begin{proof}
Let us begin with the one-dimensional case. The function space $V_n^{(1)}$ is the space of all spline functions, which 
are piece-wise polynomials of degree $m$ with discontinuities of the $(m-1)$-th derivative only at the grid points 
$\{\tfrac{k}{2^{n+1}}\mid  -2^{n}\leq k\leq 2^n\}$. These grid points divide the domain $\T$ naturally in pieces of length $2^{-(n+1)}$. The function $g$ is defined on $2^{n+1}-(m-1)$ of them. The function $\tilde g$ has to be a piece-wise polynomial of degree $m-1$ with $m-1$ pieces. To construct the coefficients 
of the function $\tilde g$, one has to solve a system of linear equations, which are independent. In fact, we have $m\,(m-1)$ 
coefficients and $(2 + (m-2))\,(m-1) = m\, (m-1)$ constraints (at the boundary and the conditions for piece-wise polynomials). 
This system has always a solution. Figure~\ref{fig:proof_1} is an illustration of the one-dimensional case. For $m=2$ a simple linear interpolation between $g(-\tfrac 12)$ and $g(-\tfrac 12+\eta)$ does the job.\\
\begin{figure}[ht]
\centering
\begin{tikzpicture}[scale=1]
\begin{axis}[scale only axis, width = 0.25\textwidth,
xmin=-1/2-0.1,
xmax=1/2+0.1,
ymin=0.5,
ymax=4,
ytick=false,
xtick={-1/2,1/2 },
xticklabels={$-\tfrac 12$, $\frac{1}{2}$},
]

\addplot[mark=none,black,domain=-1/4:1/2,smooth] { 2+3*x +0.5*cos(10*180*x)}node [black, midway,yshift=18pt] {$ g$};
\draw[blue] (axis cs:-0.5,3) -- (axis cs:-0.25,1.25) node [blue, midway,yshift=18pt] {$\tilde g$};
\draw[<->] (axis cs:-0.25,0.9) --node[midway,yshift=-10pt] {$\eta$} (axis cs:-0.5,0.9);
\filldraw (axis cs:-0.25,1.25) circle (3pt);
\filldraw (axis cs:0.5,3) circle (3pt);
\filldraw (axis cs:-0.5,3) circle (3pt);
\draw[dashed] (axis cs:-0.25,0.5) -- (axis cs:-0.25,4);
\draw[dashed] (axis cs:-0.5,0.5) -- (axis cs:-0.5,4);
\draw[dashed] (axis cs:0.5,0.5) -- (axis cs:0.5,4);
\end{axis}
\end{tikzpicture} 
\caption{The extension  of the function $g$ to $\tilde g$ for the dimension $d=1$ and the order of wavelets $m=2$.}
\label{fig:proof_1}
\end{figure}
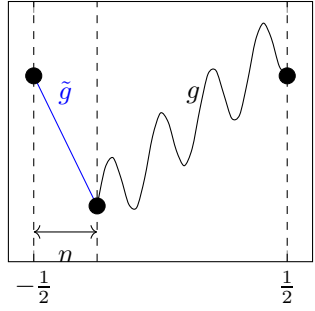  
For the multivariate case, we have to observe, that we use the index-set $I_n$ of hyperbolic structure. 
That means we need an index $\vec j$ with $|\vec j|_1\leq n$, such that $\eta \geq \frac{m-1}{2^{(j_i+1)}}$ for all $i = 1,\ldots, d$. Since all indices $j_i$ are natural numbers, multiplication and taking the $d$-th root of these inequalities leads to the condition~\eqref{eq:condition_eta}.

\end{proof} 

\begin{Remark}
If we choose not the hyperbolic index-set $I_n$, but the tensor index-set 
$$\{(\vec j,\vec k)\mid -\vec 1\leq \vec j \leq n\,\vec 1,\vec k\in \I_{\vec j}\},$$
the corresponding tensor product spaces are
$$V^{(d),\square}_n\colon = \bigotimes_{i=1}^d V_n^{(1)} = \lin\{\psi_{\vec j,\vec k}^\per\mid \norm{\vec j}_{\infty}\leq n,\vec k\in \I_{\vec j}\}.$$
For the case $d=1$ this coincides with the previous lemma. But for $d\geq 2$ we need in the previous proof only an index $\vec j$ with $|\vec j|_{\infty}\leq n$, such that $\eta \geq \frac{m-1}{2^{(j_i+1)}}$ for all $i = 1,\ldots, d$. This leads to the condition $\eta\geq \frac{m-1}{2^{n+1}}$.
\end{Remark}
\begin{Remark}
Consider the ANOVA decomposition~\eqref{eq:ANOVA_terms_transform} of a function $f\in L_2([0,1]^d,\rho)$. One ANOVA term $f_{\vec u}(\vec y_{\vec u})$ is a 
function, which depends on only $|\vec u|$ variables. Transforming $f_{\vec u}$ to $f_{\vec u}\circ \Rho_{\vec u}^{-1}$ needs only a $|\vec u|$- dimensional extension. For that reason, it is enough to choose $\eta \geq \tfrac{m-1}{2^{\left\lceil \tfrac{n}{|\vec u|}\right\rceil +1}}$ for the transformation $\Rho_{\vec u}$. We will go more into details in Chapter~\ref{sec:high-dim}. 
\end{Remark}
Motivated by the previous lemma, we make the reasonable choice
\begin{equation}\label{eq:eta}
\eta = \frac{m-1}{2^{\left\lceil \tfrac nd\right\rceil +1}}.
\end{equation}
In the remaining part of this section we will focus on the one-dimensional case, because this can be also applied to high-dimensional functions with only low-dimensional interactions. The following ideas can be generalized to $d>1$, but in this case we have to omit boundary wavelets to ensure stability.  \\

\textbf{The one-dimensional case}\\
The following lemma shows that the projection operator~\eqref{eq:defPn} applied to the extended function $\tilde f$ indeed inherits the approximation 
rate set by the order of the wavelets if the non-periodic function $f$ is smooth enough on $[0,1]$. Compare the following result with the periodic setting~\eqref{eq:P_n_2}, the only difference is the term $\tfrac{1}{(1-\eta)^{\sfrac 32}}$. 
\begin{theorem}\label{thm:nice_nonperiodic}
Let $d=1$ and the maximal wavelet level $n\in \N$. We choose the extension parameter $\eta$ as in~\eqref{eq:eta} and the transformation~\eqref{eq:Rho_1}, which gives $\tilde{\Omega} = [-\sfrac 12+\eta,\sfrac 12]$.
Then we have for the approximation error of the projection operator $P_n$ defined in~\eqref{eq:defPn} for functions $f\in H^m([0,1],\rho)$ that
\begin{equation*}
\norm{f-(P_n{\tilde f})\circ\tilde{\Rho}}_{L_2([0,1],\rho)}\lesssim  \frac{2^{-nm}}{(1-\eta)^{\sfrac 32}}  \norm{f}_{H^m([0,1],\rho)},
\end{equation*} 
where $\tilde f = f\circ \tilde \Rho^{-1}$.
\end{theorem}
\begin{proof}
For $j\in \N_0$ let us split the indices $k$ into the sets, depending on the support of the wavelet functions
\begin{align*}
I_{\text{in}} &= \{k\mid \supp \psi_{j,k}^\per\subset \tilde \Omega\},\\
I_{\text{bo}} &= \{k\mid \supp \psi_{j,k}^\per\subset [-\tfrac 12,-\tfrac 12+\eta]\},\\
I_{\text{r}} &=  \I_j\backslash (I_{\text{in}}\cup I_{\text{bo}}).
\end{align*}
First, we have a look at the wavelet coefficients of the extended function $\tilde{f}$ with $j>n$. 
From~\cite[Lemma 3.4]{LiPoUl21} we have that 
\begin{equation}\label{eq:sca_f_psi}
\langle \tilde f ,\psi_{j,k}^\per\rangle = 2^{j/2}2^{-jm}\int_{I_{j,k}}\overline{\tilde f^{(m)}}(x)\Psi_m(2^jx-k)\d x,
\end{equation}
where $I_{j,k} =\supp \psi_{j,k}^\per \subset \T$ and the function $\Psi_m$ is defined in~\eqref{eq:Psi_m}.
Since the extension of $f$ to $\tilde f$, namely $\tilde f\mid_{\T\backslash \tilde \Omega}$ is contained in the space of wavelet functions, the $m$-th derivative of $\tilde f$ is zero at the boundary $(-\sfrac 12,-\sfrac 12+\eta)$ and we have that 
\begin{equation*}
\tilde{f}^{(m)}(x) =  \delta(x+\tfrac 12) F_1(f)+ \delta(x+\tfrac 12-\eta) F_0(f),\quad  x\in [-\tfrac 12 ,-\tfrac 12+\eta],
\end{equation*}
where $\delta(\cdot)$ is the delta distribution and the numbers $F_0(f)$ and $F_1(f)$ depend on the boundary behavior of the function $f$.\\
For the indices $k\in \I_{\text{r}}$ we split the integral~\eqref{eq:sca_f_psi},
\begin{equation}\label{eq:sca_f_psi2}
\langle \tilde f ,\psi_{j,k}^\per\rangle = 2^{j/2}2^{-jm} \left(\int_{I_{j,k}\cap \tilde{\Omega} }\overline{\tilde f^{(m)}}(x)\Psi_m(2^jx-k)\d x + F_1(f)\Psi_{m}(2^j(-\tfrac 12+\eta)-k)\right),
\end{equation} 
if $-\sfrac 12+\eta\in \supp \psi_{j,k}^\per$. The other case where $-\sfrac 12\in \supp \psi_{j,k}^\per$ is analogue. 
The function $\Psi_m$ has the property that $\Psi_m(x)=0$ for $x\in \N$, see Lemma~\ref{lem:Psi_m=0}. Because of the 
choice $\eta = \frac{m-1}{2^{n+1}}$ and the assumption $j>n$, the numbers $2^j(-\tfrac 12+\eta)-k$ are in $\{1,\ldots ,m\}$. 
This allows us to omit the second term in~\eqref{eq:sca_f_psi2} and we get that 
\begin{equation*}
\langle \tilde f ,\psi_{j,k}^\per\rangle = \begin{cases} 
2^{j/2}2^{-jm}\,\int_{I_{j,k}\cap \tilde{\Omega} }\overline{\tilde f^{(m)}}(x)\Psi_m(2^jx-k)\d x, &\text{ if } k\in I_{\text{r}} \cup I_{\text{in}}, \\
0 & \text{ if } k\in I_{\text{bo}}.
\end{cases}
\end{equation*}
Using Cauchy-Schwarz inequality we receive for $k\in I_{\text{r}} \cup I_{\text{in}}$ that 
\begin{align}\label{eq:bound_coeff}
|\langle \tilde f ,\psi_{j,k}^\per\rangle|
&\leq  2^{j/2}2^{-jm} \left(\int_{I_{j,k}\cap \tilde{\Omega} }|\tilde f^{(m)}(x)|^2 \dx x \right)^{1/2}\left(\int_{I_{j,k}\cap \tilde{\Omega} }|\Psi_m(2^jx-k)|^2 \dx x \right)^{1/2}\notag\\
&\lesssim 2^{-jm}  \norm{\tilde{f}^{(m)}}_{L_2(I_{j,k}\cap \tilde{\Omega} )}\max_{x\in [0,2m-1]} \Psi_m(x)\notag\\
&\lesssim 2^{-jm} \frac{1}{(1-\eta)^{\sfrac 12}}\norm{f}_{H^m([0,1],\rho)\big|_{\Rho^{-1}(I_{j,k}\cap \tilde{\Omega} )}},
\end{align}
where the last inequality follows from~\eqref{eq:norm_eq_cube2} and \eqref{eq:Psi_bounded}.
The Riesz basis property~\eqref{eq:Riesz}, which also applies to the dual wavelets yields
\begin{equation}\label{eq:Riesz2}
\sum_{ k\in \I_ j}|\langle f,\psi_{ j, k}^{\per *}\rangle|^2\lesssim \frac{1}{\gamma_m}\left\|\sum_{ k\in \I_j}\langle f,\psi_{ j, k}^{\per *}\rangle\psi_{ j, \vec k}^{\per}\right\|^2_{L_2(\T)}
=\frac{1}{\gamma_m}\left\|\sum_{ k\in \I_j}\langle f,\psi_{  j,  k}^{\per}\rangle\psi_{j, k}^{\per *}\right\|^2_{L_2(\T)}
\lesssim \frac{\delta_m}{\gamma_m}\sum_{ k\in \I_j}|\langle f,\psi_{ j, k}^{\per }\rangle|^2.
\end{equation}
Also because of the Riesz basis property of the wavelet functions we have for $j>n$ that
\begin{equation}\label{eq:Riesz3}
\int_{-\sfrac 12+\eta}^{\sfrac 12} \left|\sum_{k=0}^{2^j-1}a_{j,k} \psi_{j,k}^\per(x)\right|^2\d x 
\leq \int_{-\sfrac 12}^{\sfrac 12} \left|\sum_{k\in I_{\text{in}}\cup I_{\text{r}} }a_{j,k} \psi_{j,k}^\per(x)\right|^2\d x \leq \delta_m \sum_{k\in I_{\text{in}}\cup I_{\text{r}} }|a_{j,k}|^2. \end{equation}
To estimate the error of the projection operator $P_n$ we have to estimate the sum of wavelet coefficients, namely we first insert the definition~\eqref{eq:defPn} of $P_n$,
\begin{align}\label{eq:error_nonper}
\allowdisplaybreaks
\norm{f-(P_n{\tilde f})\circ\Rho}^2_{L_2([0,1],\rho)} &=\frac{1}{1-\eta} \norm{\tilde{f}- P_n{\tilde f}}^2_{L_2(\tilde{\Omega})}
=\frac{1}{1-\eta} \norm{\sum_{|j|>n}\sum_{k\in \I_j}\langle\tilde f,\psi_{j,k}^{\per,*}\rangle\psi_{j,k}^{\per}}^2_{L_2(\tilde{\Omega})}\notag\\
&\leq \frac{1}{1-\eta} \left(\sum_{|j|>n}\norm{\sum_{k\in \I_j}\langle\tilde f,\psi_{j,k}^{\per,*}\rangle\psi_{j,k}^{\per}}_{L_2(\tilde{\Omega})}\right)^2\notag \\
&\stackrel{\eqref{eq:Riesz3}}{\leq} \frac{\delta_m^{\sfrac 12}}{1-\eta} \left(\sum_{|j|>n}\left(\sum_{k\in \I_j} |\langle\tilde f,\psi_{j,k}^{\per,*}\rangle|^2\right)^{\sfrac 12}\right)^2\notag\\
&\stackrel{\eqref{eq:Riesz2}}{\leq} \frac{\delta^{\sfrac 32}_m}{\gamma_m\,(1-\eta)} \left(\sum_{|j|>n}\left(\sum_{k\in \I_j} |\langle\tilde f,\psi_{j,k}^{\per}\rangle|^2\right)^{\sfrac 12}\right)^2\notag\\
&\stackrel{\eqref{eq:bound_coeff}}{\lesssim} \frac{\delta_m^{\sfrac 32}}{\gamma_m\,(1-\eta)^3}\left(\sum_{|j|>n} 2^{-jm}\left(\sum_{k\in \I_j}\norm{f}^2_{H^m([0,1],\rho)\big|_{\Rho^{-1}(I_{j,k}\cap \tilde{\Omega} )}}\right)^{\sfrac 12}\right)^2\notag\\
&\lesssim \frac{\delta_m^{\sfrac 32}}{\gamma_m\,(1-\eta)^3}\norm{f}^2_{H^^m([0,1],\rho)}\left(\sum_{|j|>n} 2^{-jm}\right)^2. 
\end{align}
The last sum of $2^{-jm}$ is bounded by geometric series, 
\begin{equation*}
\sum_{|j|>n} 2^{-jm}\leq \frac{1}{1-2^{-m}} - \frac{1-2^{-(n+1)m}}{1-2^{-m}} = \frac{2^{-(n+1)m}}{1-2^{-m}} \lesssim 2^{-nm}.
\end{equation*}  
Taking the square root in~\eqref{eq:error_nonper} gives the result, where the factor $\tfrac{\delta_m^{\sfrac 32}}{\gamma_m}$ is a constant. 
\end{proof}
Note that we receive at least in the one-dimensional case the same approximation rate as in the periodic setting. In the higher-dimensional setting this is not the case, since we loose the orthogonality between different wavelet levels in the $L_2(\tilde{\Omega})$-norm.
It is also possible to estimate the $L_\infty$-error. Here also the only change is the additional term $\tfrac{1}{1-\eta} $ in comparison to the periodic result,~\eqref{eq:P_n_infty}.
\begin{theorem}\label{thm:nice_nonperiodic2}
Let $d=1$ and the maximal wavelet level $n\in \N$. We choose the extension parameter $\eta$ as~\eqref{eq:eta}, the transformation~\eqref{eq:Rho_1} and denote $\tilde f=f\circ \Rho^{-1}$. 
Then we have for the approximation error of the projection operator 
\begin{equation*}
\norm{f-(P_n{\tilde f})\circ\Rho}_{L_\infty([0,1])}\lesssim  \frac{2^{-n(m-1/2)}}{1-\eta} \norm{f}_{H^m([0,1],\rho)}.
\end{equation*} 
\end{theorem}
\begin{proof}
Similar to the proof of the previous theorem we consider
\begin{align*}
\norm{f-(P_n{\tilde f})\circ\Rho}_{L_\infty([0,1])} &=\sup_{y\in[0,1]} |f(y) - (P_n{\tilde f})\circ\Rho(y)| 
= \sup_{x\in[-\tfrac 12+\eta,\tfrac 12]} |\tilde{f}(x) - (P_n{\tilde f})(x)| \\
&= \sup_{x\in[-\tfrac 12+\eta,\tfrac 12]} |\sum_{|j|>n}\sum_{k\in \I_j}\langle\tilde f,\psi_{j,k}^{\per,*}\rangle\psi_{j,k}^{\per}(x)|.
\end{align*}
Using the same lines as in~\cite[Theorem 3.15]{LiPoUl21} we have
\begin{align*}
\norm{f-(P_n{\tilde f})\circ\Rho}_{L_\infty([0,1])} 
&\lesssim 
\left(\sum_{|j|>n} 2^{-| j|_1(m-1/2)}\right) \, \left( \sup_{| j|>n}2^{|j|s}  \left(\sum_{ k\in \I_{ j}}|\langle \tilde{f}, \psi_{ j, k}^{*,\per}\rangle|^2\right)^{1/2} \right) \\
&\lesssim  \frac{2^{-n(m-1/2)}}{1-\eta} \norm{f}_{H^m([0,1],\rho)},
\end{align*}
which gives the assertion.
\end{proof}
In the following we discuss the numerical properties, that arise when using such an extension. We loose some stability in the sense that the wavelet matrix $\vec A$ has a bigger condition number. But nevertheless, we estimate the eigenvalues of the Moore Penrose inverse from below. The eigenvalues of the expectation matrix 
\begin{equation*}
\vec \Lambda = \left(\int_{\T}\psi_{j,k}^\per(x)\,\psi_{j',k'}^\per(x)\d x\right)_{(j,k)\in I_n,(j',k')\in I_n}
\end{equation*}
are bounded by the Riesz constants $\gamma_m$ and $\delta_m$, see~\cite[Lemma 3.18]{LiPoUl21}. The transformation $\Rho$ changes the expectation matrix to 
\begin{equation*}
\tilde{\vec \Lambda} = \left(\int_{\textcolor[rgb]{0,0,1}{\tilde\Omega}}\psi_{j,k}^\per(x)\,\psi_{j',k'}^\per(x)\d x\right)_{(j,k)\in I_n,(j',k')\in I_n}.
\end{equation*}
If we choose the extension parameter $\eta$ like~\eqref{eq:eta}, it turns out, that the eigenvalues of $\tilde{\vec \Lambda}$ do not differ much from the eigenvalues of the initial expectation matrix $\vec\Lambda$. We show this numerically. We lose the orthogonality of the wavelets of different levels. But the entries of the matrix $\tilde{\vec \Lambda}$ differ from the entries of $\vec \Lambda$ only at indices where $\supp \psi_{j,k}^\per \cap [-\tfrac 12,-\tfrac 12+\eta] \neq \varnothing$. For different maximal level $n$, we construct the matrix $\tilde{\vec \Lambda}$ and calculate the extremal eigenvalues $\mu_{\min} (\tilde{\vec \Lambda})$ and $\mu_{\max} (\tilde{\vec \Lambda})$. The results are summarized in Table~\ref{tab:eig_vals}. \begin{table}[hbt]\centering
\begin{tabular}{cc|cccccc|c}
    \hline
    $m$&$n$   & $2$&$3$ &$4$ &$5$ &$6$ &$7$  & $\vec \Lambda$\\
    \hline
		$2$&$\mu_{\min} (\tilde{\vec \Lambda})$\rule{0pt}{11pt} & $0.0896$ &$0.0903$	&$0.0879$ &$ 0.0859$&$0.0848$&  $0.0843$& $0.1481$\\
    &$\mu_{\max} (\tilde{\vec \Lambda})$\rule{0pt}{11pt}& $0.8990$ &$0.9497$	&$0.9748$ &$0.9873$&$0.9937$& $0.9968$&  $1$\\
		\hline
			$3$&$\mu_{\min} (\tilde{\vec \Lambda})$\rule{0pt}{11pt} & $0.0032$ &$0.0025$	&$0.0024$ &$ 0.0024$&$0.0024$&   $0.0025$& $0.0379$\\
    &$\mu_{\max} (\tilde{\vec \Lambda})$\rule{0pt}{11pt} & $0.7735$ &$0.8854$ &$0.9426$&$0.9714$& $0.9857$& $0.9928$& $1$\\
		\hline

  \end{tabular}
\caption{Extremal eigenvalues of the expectation matrix $\tilde{\vec \Lambda}$ for different maximal level $n$ and order $m$ of the wavelets.}
\label{tab:eig_vals}
\end{table}
 For comparison we also give the extremal eigenvalues of the matrix $\vec \Lambda$, which are the lower Riesz-constant $\gamma_m$ and $1$.
We leave the proof that $\mu_{\min} (\tilde{\vec \Lambda})>C>0$ for a constant $C$ as an open problem, but the numeric indicates that this is true.\\
In fact, the choice of $\eta$ must be a balance between an ill-conditioned matrix for big $\eta$ and a large approximation 
error for small $\eta$. The choice~\eqref{eq:eta} does this job. \\

\textbf{Error estimates}\\
Up to now we gave estimates for the $L_2$ and $L_\infty$ error of the projection operator $P_n$,~\eqref{eq:defPn}, instead of the approximation operator $S^{ \X}_n$. To end this subsection, we give also for the non-periodic case an error estimate with high probability, similar to Theorem~\ref{thm:decay_trafo}.
\begin{corollary}\label{cor:nonper_end}
Let $d=1$, $m\in \N$ be the order of vanishing moments of the wavelet $\psi$, $\rho\in C^{m-1}([0,1])$ be a density and and $n\in \N$ be the maximal level of the wavelets. We choose the transformation $\Rho$ as~\eqref{eq:Rho_1} with the extension parameter $\eta$ as in~\eqref{eq:eta}. Let furthermore $\Y = (y_i)_{i=1}^M\subset [0,1]$ with $M\gtrsim r\,2^{n+1} \, (n+1) $ be drawn i.i.d.\,at random according to $\rho$ and $r>1$. Then we denote by $ \X=\Rho(\Y)$ the transformed samples. Let furthermore $\mu_{\min}(\tilde{\vec \Lambda}) \geq C>0$.
Then we have 
\begin{equation*}
\P\left(\norm{f-(S_n^{ \X}{\tilde f})\circ\Rho}_{L_2([0,1],\rho)}\lesssim \frac{2^{-nm}}{(1-\eta)^{\sfrac 52}}\,\norm{f}_{H^m([0,1],\rho)}\right)\geq 1-2M^{-r}.
\end{equation*} 
\end{corollary}
\begin{proof}
First, we denote $e_2:=\norm{\tilde{f}-(P_n\tilde{f}) }_{L_2(\tilde \Omega)}$ and $e_\infty:=\norm{\tilde f-(P_n\tilde{f})}_{L_\infty(\tilde \Omega)}$. By using the extension we loose the orthogonality of the wavelet functions even for different levels, since the $L_2$-norm is then defined on $\tilde \Omega$ and not on the whole torus. Therefore, we have to modify the proof of~\cite[Theorem 3.20]{LiPoUl21} slightly. We have
\begin{align*}\norm{f-(S_n^{ \X} \tilde f) \circ{ \Rho}}_{L_2([0,1],\rho)} 
&= \frac{1}{1-\eta} 
\norm{\tilde f-S_n^{ \X} \tilde f}_{L_2(\tilde \Omega)}
\leq \frac{1}{1-\eta} \left(e_2+\norm{P_n \tilde f-S_n^{ \X} \tilde f}_{L_2(\tilde \Omega)}\right)\notag\\
&=\frac{1}{1-\eta} \left(e_2+\norm{S_n^{ \X} (P_n \tilde f- \tilde f)}_{L_2(\T^d)} \right)
\leq \frac{1}{1-\eta} \left(e_2+\norm{S_n^\X}_2 \norm{P_n f- f}_{\ell_2(\X)}\right).
\end{align*} 
For the last term we use the same lines as in~\cite[Theorem 3.20]{LiPoUl21}, which are based on Bernstein's inequality to get
\begin{equation*}
\P\left(\norm{\tilde f( x_i)-P_n\tilde f( x_i)}_{\ell_2(\X)} \geq  \left(M \sqrt{\tfrac{2\, e_\infty^2e_2^2 r\log M}{M} } + \tfrac{2\,e_\infty^2 r \log M}{3M}+e_2^2\right)^{\sfrac 12}\right)\leq M^{-r}.
\end{equation*}  
Taking the event into account that 
$$\norm{S_n^{ \X}}_2\leq \sqrt{\frac{2}{M\mu_{\min}(\tilde{\vec \Lambda})}} <\sqrt{\frac{2}{C}}  $$
with high probability, we obtain by union bound the overall probability exceeding the sum of the probabilities, i.e.
\begin{equation*}
\P\left(\norm{\tilde f-S_n^{\X} \tilde f}_{L_2(\tilde \Omega)} \leq e_2 + \sqrt{\tfrac{2}{C}}\,\left(\tfrac{e_2^2}{ M}+e_2e_{\infty}\sqrt{\tfrac{r\,\log M}{M}}+ e_{\infty}^2\tfrac{r\,\log M}{M}\right)^{\sfrac 12}\right)
\geq 1- 2\,M^{-r}.
\end{equation*}
Collecting the bounds from the occurring terms from Theorems~\ref{thm:nice_nonperiodic} and~\ref{thm:nice_nonperiodic2} as well as logarithmic oversampling, which means $\tfrac{\log M}{M}\lesssim 2^{-n}$, we end up with 
\begin{align*}
&\norm{f-(S_n^{ \X} \tilde f) \circ{ \Rho}}_{L_2([0,1],\rho)} 
= \frac{1}{1-\eta} 
\norm{\tilde f-S_n^{\X} \tilde f}_{L_2(\tilde \Omega)} \\
&\quad\lesssim  \frac{1}{1-\eta} \left(\frac{2^{-nm}}{(1-\eta)^{\sfrac 32}} + \sqrt{\tfrac{2}{C}}\,\left(\frac{2^{-2nm}}{(1-\eta)^{3}\,2^{n+1}}+\frac{2^{-2nm}}{(1-\eta)^{\sfrac 52}}\,\sqrt{r}+ \frac{2^{-2nm}}{(1-\eta)^2}\,r\right)^{\sfrac 12}\right) \cdot \norm{f}_{H^m([0,1],\rho)}\\
&\quad= \frac{2^{-nm}}{(1-\eta)^{\sfrac 52}} \left(1+\sqrt{\tfrac{2}{C}}\,\left(\frac{1}{(1-\eta)^{\sfrac 32} \, 2^{n+1}}+ \frac{\sqrt{r}}{(1-\eta)}+ \frac{r}{(1-\eta)^{\sfrac 12}}\right)^{\sfrac 12}\right)\cdot \norm{f}_{H^m([0,1],\rho)}
\end{align*}
with high probability.
\end{proof}

\subsection{Numerical Experiments}\label{sec:Numerics1}
In this section we study the approximation behavior numerically for some examples to underpin our findings from Theorem~\ref{thm:decay_trafo} and Corollary~\ref{cor:nonper_end}. We consider examples where $\Omega \in \{\R^d,\T^d,[0,1]^d\}$.  \\
We do the following procedure. For a maximal level $n$, we draw $M\gtrsim N\log N \asymp 2^n\,n^d$ i.i.d.\,samples according to the density $\rho$, 
which coincides with logarithmic oversampling. Also the corresponding function values are given.  
In our approximation, we transfer the samples in the set $\Y$ to the torus, by $\Rho(\Y)=\X$ and apply the approximation operator 
$S_n^{\X}$ given in~\eqref{eq:def_S_n^X}. A good estimator for the $L_2(\Omega,\rho)$-error, is the \textit{root mean square error} (RMSE), which is defined by
\begin{equation}\label{eq:RMSE}
\text{RMSE} = \left(\frac{1}{|\Y_{\text{test}}|}\sum_{\vec y\in \Y_{\text{test}}}|f(\vec y)-(S_n^{\X} (f\circ \Rho^{-1}))\circ \Rho(\vec y)|^2\right)^{1/2},
\end{equation} 
for sample points $\Y_{\text{test}}\subset \Omega$, which are i.i.d.\,according to $\rho$. We use always $|\Y_{\text{test}}|=3\,|\Y| = 3\,M$.
We defined in Examples~\ref{ex:rhos} and~\ref{ex:beta} only one-dimensional densities. In the following we interpret the $d$-dimensional densities as a tensor product in the sense of~\eqref{eq:rho_d}.\\

\textbf{Distributions on $\R^d$}\\
We begin with the normal distribution $\rho_N$ from~\eqref{eq:rho_N} for all one-dimensional densities $\rho_i(y_i)$. As test function we use the Gaussian
\begin{equation}\label{eq:f_gauss}
f: \R^d\to \R, \quad f(\vec y) = \e^{-\norm{\vec y}_2^2},
\end{equation}
which is smooth. But in the sense of~\eqref{eq:Hsmix_rho}, we have to note that in the one-dimensional-case
\begin{align*}
\norm{f}_{H^1(\R,\rho_N)}^2&= \tfrac{1}{\sqrt 5} + \tfrac{8\pi}{3\sqrt 3} \approx 5.28<\infty,\\
\norm{f}_{H^2(\R,\rho_N)}^2&= \tfrac{1}{\sqrt 5} + (\tfrac{8\pi}{3\sqrt 3}+48\pi^2)+144 \pi^2 \approx 1900<\infty,\\
\norm{f}_{H^3(\R,\rho_N)}^2&=\infty.
\end{align*}
Due to the embedding $H^2(\T)\subset \bB_{2,\infty}^2(\T)$ we also have $f\circ\Rho^{-1}\in\bB_{2,\infty}^2(\T)$. To investigate the smoothness further, we have a look at the terms from Definition~\ref{def:besov_on_Omega}, i.e.\,for $s=\sfrac 52$
\begin{equation*}
 \sup_{ \stackrel{j\in \Z}{j\neq 0} }2^{5| j|_1} \sup_{\alpha \leq 2} \sup_{1\leq  \beta\leq \alpha} \sum_{k\in J_{j}}\frac{1}{ |k|^{2\alpha}}\left|c_{k}^\rho\left(\Dx^{ \beta} f(y)B_{ \alpha,\beta}(y)\right)\right|^2.
\end{equation*}
 We have 
\begin{align*}
|c_k^\rho(\Dx^2 f(y) B_{2,2}(y))| &= |c_k^\rho((-2+4y^2)\e^{-y^2} \,2\pi\,\e^{y^2})| = |c_k^\rho((-2+4y^2)\,2\pi)| = |2\pi \,c_k(-2+4(\Rho^{-1}(x))^2)| \asymp \tfrac{ 1}{|k|},
\end{align*}
since $\Rho^{-1}(x)^2$ is a smooth function except at $x=\sfrac 12$. Analogously, we get for the other term 
$$|c_k^\rho(\Dx^1 f(y) B_{2,1}(y))| \asymp \tfrac{ 1}{|k|}.$$
 This yields
\begin{equation*}
\norm{f}_{\bB_{2,\infty}^{\sfrac 52}}^2\lesssim \sup_{j \in \Z} 2^{5|j|} \sum_{k\in J_j} \frac{1}{|k|^4}\frac{1}{|k|^2}\lesssim \sup_{j \in \Z} 2^{5|j|} 2^{|j|-1} 2^{-6|j|}<\infty,
\end{equation*}
since the index sets $J_j$ are defined such that $k>2^{|j|-1}$ and $|J_j|= 2^{|j|}$.
The function $f$ as well as the density $\rho_N$ have a tensor product structure, hence we get 
$f\circ\Rho^{-1}\in \bB^{5/2}_{2,\infty}(\T^d)$ and $f\in \bB^{5/2}_{2,\infty}(\Omega,\rho_N)$. \\

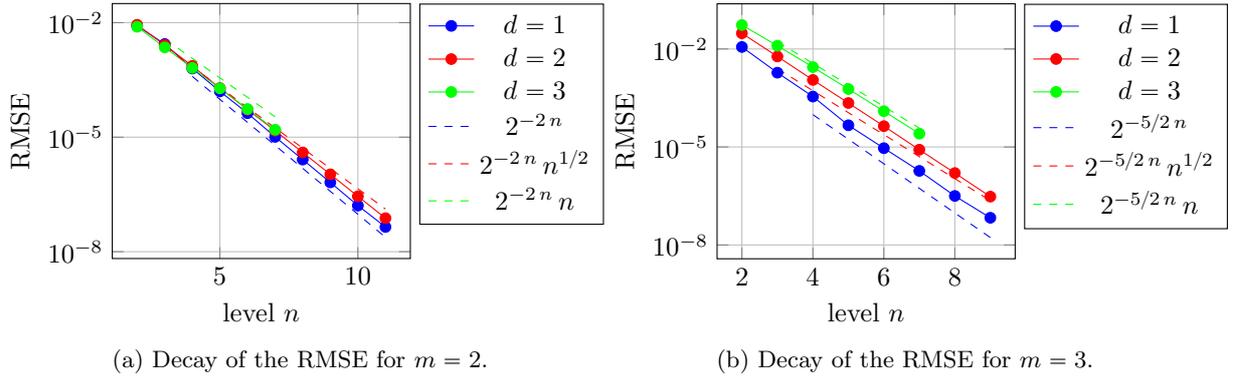
\begin{figure}[ht]
\begin{subfigure}[c]{0.49\textwidth}
\raggedleft
\begin{tikzpicture}[scale=1]
\begin{semilogyaxis}[scale only axis, width = 0.5\textwidth,xlabel = {level $n$},ylabel ={RMSE},
legend entries={$d=1$ ,$d=2$,$d=3$,$2^{- 2\,n}$,$2^{-2\,n  }\,n^{1/2}$,$2^{- 2\,n  }\,n$},
grid=major,
legend pos=outer north east,
]
\addplot[blue,mark=*]  coordinates { 		
(2,0.008258773647310504)
(3,0.002753369251649594)
(4,0.0006285377780795004)
(5,0.00015851710117530227)
(6,4.249536229931639e-5)
(7,1.0192639217091407e-5)
(8,2.5868966354994144e-6)
(9,6.529679140406846e-7)
(10,1.6140706883031577e-7)
(11,4.44205415886394e-8)

};

\addplot[red,mark=*]  coordinates {

(2,0.008774848533975485)
(3,0.0025902834099958694)
(4,0.0007346972726777058)
(5,0.00019378808199678414)
(6,5.270916797130816e-5)
(7,1.4785661139324594e-5)
(8,3.977117573449586e-6)
(9,1.0635947648887423e-6)
(10,2.823607247290891e-7)
(11,7.548457555224825e-8)

};

\addplot[green,mark=*]  coordinates {

(2,0.007801122471648614)
(3,0.0022480672037606224)
(4,0.0006528735380641511)
(5,0.00018992264864414743)
(6,5.446185233621944e-5)
(7,1.564144184443896e-5)

};

\addplot[blue,mark=none,dashed] coordinates {(4,0.1*1/2^8) (11,0.1*1/2^22)};
\addplot[red,mark=none,dashed] coordinates {(3,0.1* 3^1/2 * 2^-2*3 ) (11, 0.1   *11^1/2 * 2^-2*11 )};
\addplot[green,mark=none,dashed] coordinates {(3, 0.08*3 * 2^-2*3 ) (7, 0.08*7  * 2^-2*7 )};

\end{semilogyaxis}
\end{tikzpicture} \subcaption{Decay of the RMSE for $m=2$. }

\end{subfigure}
\begin{subfigure}[c]{0.49\textwidth}
\raggedright
\begin{tikzpicture}[scale=1]
\begin{semilogyaxis}[scale only axis,width = 0.5\textwidth,xlabel = {level $n$},ylabel ={RMSE},
legend entries={$d=1$ ,$d=2$,$d=3$,$2^{- 5/2\,n }$,$2^{-5/2\, n }\,n^{1/2}$,$2^{- 5/2\,n }\,n$},
grid=major,
legend pos=outer north east,
]
\addplot[blue,mark=*]  coordinates { 		
(2,0.011512865308148481)
(3,0.0018810135363476923)
(4,0.00034948931284204375)
(5,4.619655707876724e-5)
(6,9.181459569781514e-6)
(7,1.851280524734041e-6)
(8,3.166773644993778e-7)
(9,6.799754473331349e-8)
};

\addplot[red,mark=*]  coordinates {
(2,0.03020269258510617)
(3,0.00588902611792983)
(4,0.0011284497714056485)
(5,0.00022195921214862452)
(6,4.3380628748332615e-5)
(7,8.254160992616864e-6)
(8,1.5933250166483387e-6)
(9,3.0239660278928993e-7)

};

\addplot[green,mark=*]  coordinates {

(2,0.05398991975248214)
(3,0.012582244742030615)
(4,0.002782951620647981)
(5,0.000594539734852015)
(6,0.00012450839681408157)
(7,2.557344831352366e-5)

};

\addplot[blue,mark=none,dashed] coordinates {(4,0.1*1/2^10) (9,0.1*2^-5/2*9)};
\addplot[red,mark=none,dashed] coordinates {(3,   0.3* 3^1/2 *2^-2.5*3 ) (9,  0.3* 9^1/2* 2^-2.5*9 )};
\addplot[green,mark=none,dashed] coordinates {(3, 3 * 2^-5/2*3 ) (7,  7  * 2^-5/2*7 )};

\end{semilogyaxis}
\end{tikzpicture} \subcaption{Decay of the RMSE for $m=3$. }
\end{subfigure}

\caption{Approximation of the function~\eqref{eq:f_gauss} on $\R^d$ for $d\in \{1,2,3\}$ and the normal distribution $\rho_N$.}
\label{fig:num_gauss}
\end{figure}

We did the approximation for $d\in \{1,2,3\}$ and for the order $m=\{2,3\}$ of vanishing moments of the wavelet. In Figure~\ref{fig:num_gauss} 
we plotted the results. One can see, that we end up with the proposed error decay rates from Theorem~\ref{thm:decay_trafo}. In case where $m=2$, 
we are in the case \eqref{eq:thm_s=m} and get the proposed error decay rate of $2^{-2n}\,n^{(d-1)/2}$. If we increase the order of 
vanishing moments of the wavelets to $m=3$, we are in the case~\eqref{eq:thm_s<m} and receive also numerically the proposed error decay 
of $2^{-5/2n}\,n^{(d-1)/2}$.
The numerical results are even slightly better in some cases. \\

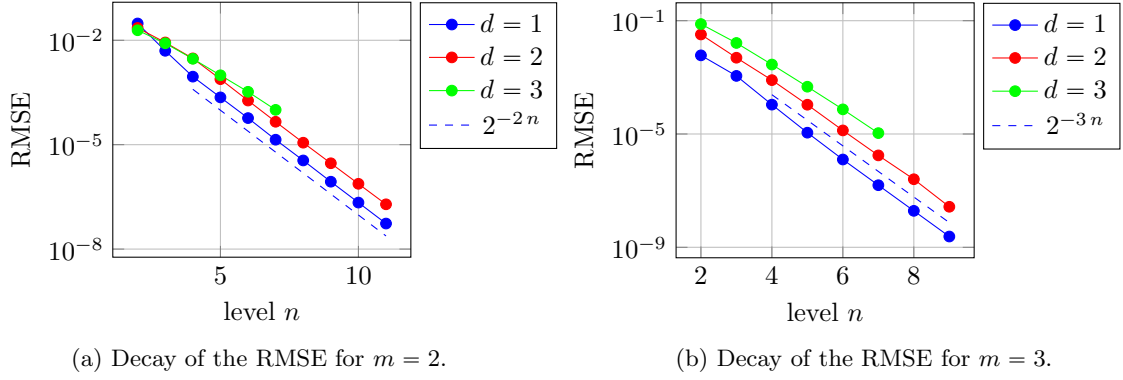
\begin{figure}[ht]
\begin{subfigure}[c]{0.49\textwidth}
\raggedleft
\begin{tikzpicture}[scale=1]
\begin{semilogyaxis}[scale only axis, width = 0.5\textwidth,xlabel = {level $n$},ylabel ={RMSE},
legend entries={$d=1$ ,$d=2$,$d=3$,$2^{- 2\,n}$},grid=major,
legend pos=outer north east,
]
\addplot[blue,mark=*]  coordinates { 		
(2,0.030386846783657842)
(3,0.0049848825626384275)
(4,0.0009013677358339524)
(5,0.0002286451675239334)
(6,5.81156909780464e-5)
(7,1.3898887379239857e-5)
(8,3.5210068137442513e-6)
(9,8.632405199418129e-7)
(10,2.168918375174088e-7)
(11,5.4174112128971905e-8)

};

\addplot[red,mark=*]  coordinates {
(2,0.022601191801632912)
(3,0.008649278731015342)
(4,0.003017125570420653)
(5,0.0007626804129802458)
(6,0.00018550205703782494)
(7,4.552852385112095e-5)
(8,1.1401245200547176e-5)
(9,2.9250421541633564e-6)
(10,7.514448082323972e-7)
(11,1.9350541972610507e-7)

};

\addplot[green,mark=*]  coordinates {

(2,0.019402515254722474)
(3,0.008193812709938844)
(4,0.0029367622503568462)
(5,0.0009865877606407252)
(6,0.000330959366139055)
(7,0.00010014700820896297)

};

\addplot[blue,mark=none,dashed] coordinates {(4,0.1*1/2^8) (11,0.1*1/2^22)};

\end{semilogyaxis}
\end{tikzpicture} \subcaption{Decay of the RMSE for $m=2$. }

\end{subfigure}
\begin{subfigure}[c]{0.49\textwidth}
\raggedright
\begin{tikzpicture}[scale=1]
\begin{semilogyaxis}[scale only axis,width = 0.5\textwidth,xlabel = {level $n$},ylabel ={RMSE},
legend entries={$d=1$ ,$d=2$,$d=3$,$2^{- 3\,n }$},grid=major,
legend pos=outer north east,
]
\addplot[blue,mark=*]  coordinates { 		
(2,0.005980561548994862)
(3,0.0011246882735062309)
(4,0.00010878659091386085)
(5,1.1130405607849689e-5)
(6,1.257742385095359e-6)
(7,1.556168274507245e-7)
(8,1.935857224601924e-8)
(9,2.401538820134692e-9)

};

\addplot[red,mark=*]  coordinates {
(2,0.03263143758417272)
(3,0.004931385580086986)
(4,0.0007982032760036581)
(5,0.00010803236474046022)
(6,1.3368685852431527e-5)
(7,1.756823705262893e-6)
(8,2.528237420880295e-7)
(9,2.6968124886065625e-8)

};

\addplot[green,mark=*]  coordinates {
(2,0.07621934578060434)
(3,0.016371708899960383)
(4,0.0028349520848911843)
(5,0.00046723200456255665)
(6,7.405490754142548e-5)
(7,1.0751997504974189e-5)

};

\addplot[blue,mark=none,dashed] coordinates {(4, 1*2^-3*4) (9,1*2^-3*9)};

\end{semilogyaxis}
\end{tikzpicture} \subcaption{Decay of the RMSE for $m=3$. }
\end{subfigure}

\caption{Approximation of the function~\eqref{eq:f_gauss} on $\R^d$ for $d\in \{1,2,3\}$  with the Cauchy distribution $\rho_C$.}
\label{fig:num_cauchy}
\end{figure}

\begin{figure}[ht]
\begin{subfigure}[c]{0.49\textwidth}
\raggedleft
\begin{tikzpicture}[scale=1]
\begin{semilogyaxis}[scale only axis, width = 0.5\textwidth,xlabel = {level $n$},ylabel ={RMSE},
legend entries={$d=1$ ,$d=2$,$d=3$,$2^{- 2\,n}$},grid=major,
legend pos=outer north east,
]
\addplot[blue,mark=*]  coordinates { 	

(2,0.13173276080304086)
(3,0.028194182675750072)
(4,0.0045767581282938384)
(5,0.001011349520268819)
(6,0.0002377402703301097)
(7,5.911869572715196e-5)
(8,1.4424567679062613e-5)
(9,3.587922166044898e-6)
(10,9.031149452471319e-7)
(11,2.2581181793510832e-7)

};

\addplot[red,mark=*]  coordinates {
(2,0.06428160574607264)
(3,0.05053830772912623)
(4,0.030111047017384644)
(5,0.017092380636817724)
(6,0.004259857096231651)
(7,0.0009509665384978196)
(8,0.00022185447871083934)
(9,5.392076369733504e-5)
(10,1.3798024252863419e-5)
(11,3.5488056257559013e-6)

};

\addplot[green,mark=*]  coordinates {

(2,0.02579091110817439)
(3,0.022978438155647218)
(4,0.01856308919045666)
(5,0.01312035856296072)
(6,0.008050607907023383)
(7,0.004359261762499709)

};

\addplot[blue,mark=none,dashed] coordinates {(4,0.5*1/2^8) (11,0.5*1/2^22)};
\addplot[blue,mark=none,dashed] coordinates {(6, 30 * 2^-2*6 ) (11, 30   * 2^-2*11 )};

\end{semilogyaxis}
\end{tikzpicture} \subcaption{Decay of the RMSE for $m=2$. }

\end{subfigure}
\begin{subfigure}[c]{0.49\textwidth}
\raggedright
\begin{tikzpicture}[scale=1]
\begin{semilogyaxis}[scale only axis,width = 0.5\textwidth,xlabel = {level $n$},ylabel ={RMSE},
ytick = {10^-2,10^-5,10^-8},
legend entries={$d=1$ ,$d=2$,$d=3$,$2^{- 3\,n }$,$2^{-3\, n }\,n^{1/2}$,$2^{- 3\,n }\,n$},
grid=major,
legend pos=outer north east,
]
\addplot[blue,mark=*]  coordinates {
(2,0.04420650476877685)
(3,0.009650904904477866)
(4,0.0008376306574596207)
(5,7.365499747887757e-5)
(6,9.317238399100896e-6)
(7,1.1399487809570894e-6)
(8,1.4062906242031953e-7)
(9,1.7473071750490557e-8)

};

\addplot[red,mark=*]  coordinates {
(2,0.061197022718758175)
(3,0.034823285029110064)
(4,0.01016485939212716)
(5,0.002759088682206801)
(6,0.0005221270445895316)
(7,8.843356184879798e-5)
(8,1.1199206494646483e-5)
(9,1.3064011939393448e-6)

};

\addplot[green,mark=*]  coordinates {
(2,0.025197507685750808)
(3,0.02008871675863307)
(4,0.011268193209981283)
(5,0.006562717168402755)
(6,0.002222664753392578)
(7,0.00065453939823561)

};

\addplot[blue,mark=none,dashed] coordinates {(4, 1*2^-3*4) (9,1*2^-3*9)};
\addplot[red,mark=none,dashed] coordinates {(3, 20* 3^1/2 *2^-3*3 ) (9,  20* 9^1/2* 2^-3*9 )};

\end{semilogyaxis}
\end{tikzpicture} \subcaption{Decay of the RMSE for $m=3$. }
\end{subfigure}

\caption{Approximation of the function~\eqref{eq:f_gauss} on $\R^d$ for $d\in \{1,2,3\}$  with the Laplace distribution $\rho_L$.}
\label{fig:num_laplace}
\end{figure}
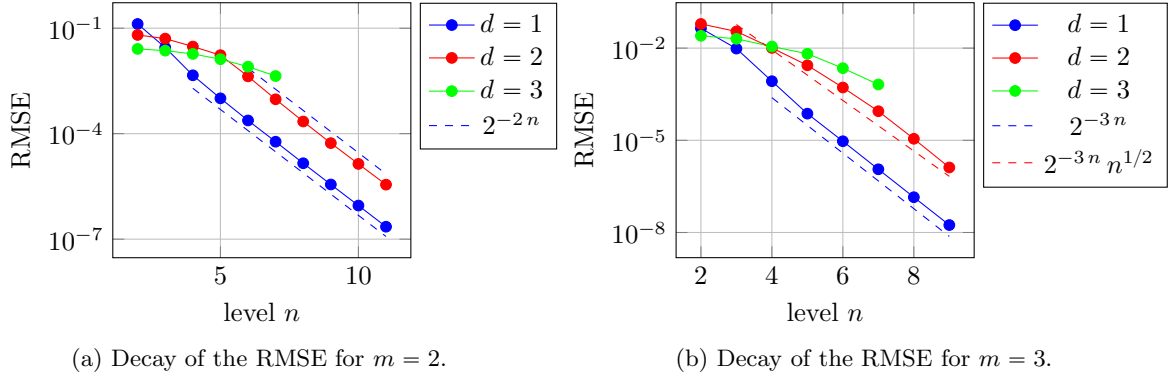

 Note that a different density $\rho(y)$ can lead to a different smoothness of $f\circ\Rho^{-1}$ even 
for the same function $f$. For a tensor product of the one-dimensional Cauchy $\rho_C$, \eqref{eq:rho_C}, or the Laplace distribution $\rho_L$, \eqref{eq:rho_L}, we even have $f\circ{\Rho^{-1}}\in H^{s}_\mix(\Omega,\rho)$ for all $s\in \N$. 
That follows from $\Dx^\alpha f(y)=p(y)\e^{-y^2/2}$ for a polynomial $p(y)$ and the fact that all differentials of the Cauchy as well as the Laplace distribution are polynomials or polynomials with a the factor of the behavior $\e^{k|y-2|}$. Furthermore, all integrals of the form 
$\int_{\R}\e^{-y^2+k|y-2|}p(y)\d y $ are finite. In Figure~\ref{fig:num_cauchy} and \ref{fig:num_laplace} we plotted the approximation results. In this case the order of vanishing moments determines the error decay. In dimension $d=3$ with Laplace distribution we are still in the preasymptotic case.
\\
  

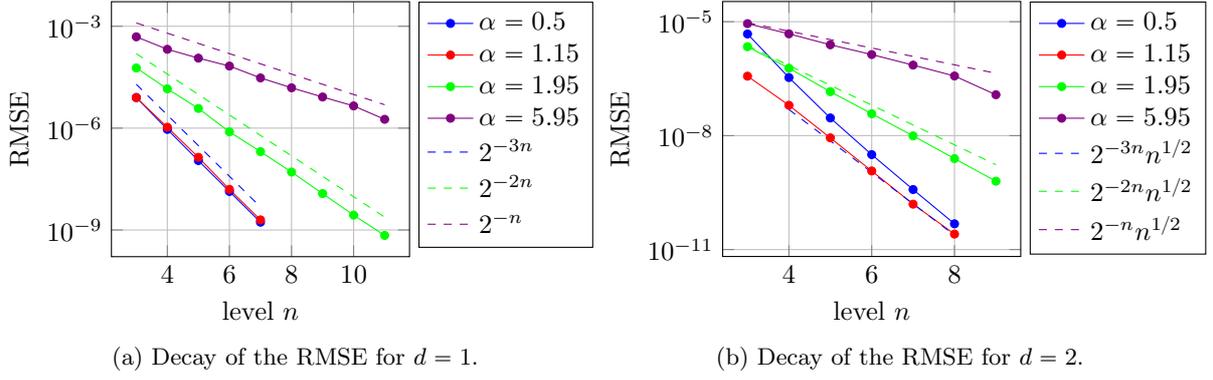
\begin{figure}[ht]
\begin{subfigure}[t]{0.49\textwidth}
\raggedleft
\begin{tikzpicture}[scale=1]
\begin{semilogyaxis}[scale only axis, width = 0.5\textwidth,xlabel = {level $n$},ylabel ={RMSE},
legend entries={$\alpha=0.5$ ,$\alpha=1.15$,$\alpha=1.95$,$\alpha=5.95$,$2^{-3n}$,$2^{-2n}$,$2^{-n}$},
grid=major,
legend pos=outer north east,
legend cell align={left}
]
\addplot[blue,mark=*,mark size=1.5pt]  coordinates { 		

(3,7.973818431916849e-6)
(4,9.163594256662534e-7)
(5,1.107155197829408e-7)
(6,1.3626437749224019e-8)
(7,1.6959705838012245e-9)

};

\addplot[red,mark=*,mark size=1.5pt]  coordinates {

(3,7.905667323890103e-6)
(4,1.0644207945131054e-6)
(5,1.3817514690918609e-7)
(6,1.5722323713252653e-8)
(7,1.982095364146727e-9)

};

\addplot[green,mark=*,mark size=1.5pt]  coordinates {

(3,5.903249933938167e-5)
(4,1.4363292343306572e-5)
(5,3.79432929077689e-6)
(6,7.783621818079864e-7)
(7,2.031086024253903e-7)
(8,5.087253858955471e-8)
(9,1.1831108212933865e-8)
(10,2.737384428298854e-9)
(11,6.913707993942589e-10)

};
\addplot[violet,mark=*,mark size=1.5pt]  coordinates { 		
(3,0.000488345342111391)
(4,0.0002094670808135388)
(5,0.00011510901680976869)
(6,6.727843050238251e-5)
(7,3.0092753642515774e-5)
(8,1.5391706109718624e-5)
(9,8.254709322029301e-6)
(10,4.504416699846754e-6)
(11,1.8154030898758224e-6)

};

\addplot[blue,mark=none,dashed] coordinates {(3,0.01*2^-3*3) (7,0.01*2^-3*7)};
\addplot[green,mark=none,dashed] coordinates {(3, 0.01* 2^-2*3 ) (11, 0.01 * 2^-2*11 )};
\addplot[violet,mark=none,dashed] coordinates {(3, 0.01 * 2^-1*3 ) (11, 0.01 *  2^-1*11 )};

\end{semilogyaxis}
\end{tikzpicture} \subcaption{Decay of the RMSE for $d=1$. }
\end{subfigure}
\begin{subfigure}[t]{0.49\textwidth}
\centering
\begin{tikzpicture}[scale=1]
\begin{semilogyaxis}[scale only axis,width = 0.5\textwidth,xlabel = {level $n$},ylabel ={RMSE},
legend entries={$\alpha=0.5$ ,$\alpha=1.15$,$\alpha=1.95$,$\alpha=5.95$,$2^{-3n}n^{1/2}$,$2^{-2n}n^{1/2}$,$2^{-n}n^{1/2}$},
grid=major,
legend pos=outer north east,
legend cell align={left}
]
\addplot[blue,mark=*,mark size=1.5pt]  coordinates { 		

(3,4.750495168294732e-6)
(4,3.3842174830783615e-7)
(5,2.9194776299243706e-8)
(6,3.164227825505167e-9)
(7,3.803899523618314e-10)
(8,4.7644725435215784e-11)

};

\addplot[red,mark=*,mark size=1.5pt]  coordinates {
(3,3.677266994269429e-7)
(4,6.302048304865215e-8)
(5,8.7365794459546e-9)
(6,1.1726249832195715e-9)
(7,1.5680509442286305e-10)
(8,2.5542507141109844e-11)

};

\addplot[green,mark=*,mark size=1.5pt]  coordinates {

(3,2.1954054493315066e-6)
(4,5.96568507883468e-7)
(5,1.4401615633240354e-7)
(6,3.7697076340962665e-8)
(7,9.843220311048168e-9)
(8,2.4844439954703463e-9)
(9,6.307459038282379e-10)

};

\addplot[violet,mark=*,mark size=1.5pt]  coordinates { 		
(3,8.8939599032935e-6)
(4,4.788929555137152e-6)
(5,2.466402964897063e-6)
(6,1.3642069723428778e-6)
(7,7.19281932179946e-7)
(8,3.7185476198114113e-7)
(9,1.1933595867541537e-7)

};

\addplot[blue,mark=none,dashed] coordinates {(4,0.0001*  4^1/2 * 2^-3*4) (8,0.0001 * 8^1/2 * 2^-3*8)};
\addplot[green,mark=none,dashed] coordinates {(3, 0.0001* 3^1/2 * 2^-2*3 ) (9, 0.0001 * 9^1/2  * 2^-2*9  )};
\addplot[violet,mark=none,dashed] coordinates {(3, 0.00005 *3^1/2 * 2^-1*3 ) (9, 0.00005* 9^1/2 *  2^-1*9 )};

\end{semilogyaxis}

\end{tikzpicture} \subcaption{Decay of the RMSE for $d=2$. }

\end{subfigure}

\caption{Approximation on $\T^d$ for $d\in \{1,2\}$ of the test function~\eqref{eq:f_interval} samples distributed with respect to the beta distribution $\rho_{B,\alpha}$.}
\label{fig:num_interval}
\end{figure}

\textbf{The beta distribution on the torus}\\
Let us consider the beta distribution~\eqref{eq:rho_B}, but shifted by $-\sfrac 12$ to the torus $\T=[-\sfrac 12,\sfrac 12)$, for all one-dimensional densities $\rho_i(x_i)$, 
which also includes the uniformly distribution on $\T$ for $\alpha=1$. 
We choose as test function
\begin{equation}\label{eq:f_interval}
f:\T^d\to \R,\quad f(\vec y ) = \prod_{i=1}^d (y_i-\tfrac 12)^3 (y_i+\tfrac 12)^3,
\end{equation}
which is the tensor product of a polynomial of degree $6$ and has triple zeros at $\sfrac 12$ and is in $H^3_{\mix}(\T^d)$. Depending on the choice 
of the parameter $\alpha$ of the beta distribution, the transformed function has different regularity. We consider first the one-dimensional case. 
The crucial points to decide whether the norm $\norm{f}_{L_2(\T,\rho_{B,\alpha})}$ is finite is at the points with lower regularity $\sfrac 12 \cong - \sfrac 12$. Because of the symmetry of the 
function $f$ as well as the density $\rho$, it is sufficient to have a look at the point $y = \sfrac 12$. There we have the behavior $f^{(i)}\sim (y-\sfrac 12)^{3-i}$ 
for $i=0,1,2,3$. With the same arguments we have that $\rho_{B,\alpha}^{(i)}\sim (y-\sfrac 12)^{\alpha-1-i}$ for $\alpha\neq 1,2$ and $i=0,1,2$.
Furthermore the integral $\int_{0}^1 x^k\d x $ is finite for $k\in \R$ and $k>-1$. 
Hence, Definition~\ref{def:HmNorm_1d} gives that we have to require the following:
\begin{align*}
f\in H^1(\T,\rho_{B,\alpha}) &\Rightarrow f' \in L_2(\T,\Upsilon_{1,1})&&\Rightarrow \alpha<6,\\
f\in H^2(\T,\rho_{B,\alpha}) &\Rightarrow f'' \in L_2(\T,\Upsilon_{2,2}),
											 f' \in L_2([0,1],\Upsilon_{1,2})&&\Rightarrow  \alpha<2,\\
f\in H^3(\T,\rho_{B,\alpha}) &\Rightarrow f''' \in L_2(\T,\Upsilon_{3,3}),
							 f'' \in L_2([0,1],\Upsilon_{2,3}),
							 f' \in L_2([0,1],\Upsilon_{1,3}) &&\Rightarrow \alpha<6/5.
\end{align*}
Since, the function $f$ is a tensor product, we have the same estimates for the multivariate cases. 
Indeed, we used the order of vanishing moments $m=3$ of the wavelets, which limits the maximal error decay rate and in Figure~\ref{fig:num_interval} is the 
resulting numerical approximation decay for different parameters $\alpha$ and $d$. Indeed, if we are below the critical values $\sfrac 65 ,2$ and $6$ for $\alpha$, we get the desired approximation rates of $2^{-3n}$, $2^{-2n}$ and $2^{-n}$ given in Theorem~\ref{thm:decay_trafo}.\\
Even dense samples at the boundary, which coincides with small $\alpha$ can not increase the error decay rate of $3$.
 \\

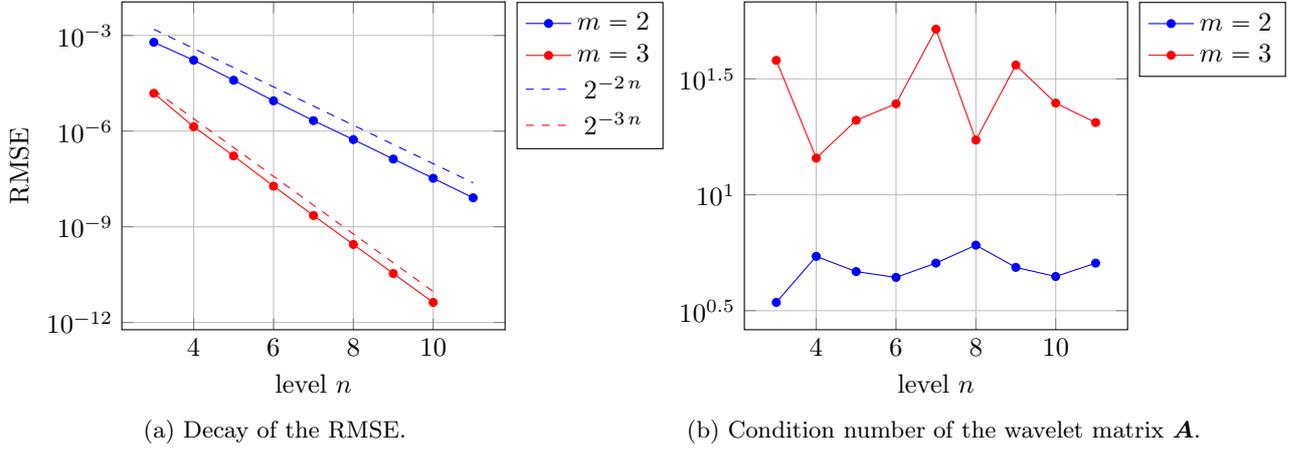
\begin{figure}[ht]
\centering
\begin{minipage}[t]{0.45\linewidth}
\centering
\begin{tikzpicture}[scale=1]
\begin{semilogyaxis}[scale only axis, width = 0.7\textwidth,xlabel = {level $n$},ylabel ={RMSE},
legend entries={$m=2$, $m=3$,$2^{- 2\,n }$,$2^{-3\, n }$},
grid=major,
legend pos=outer north east,
]
\addplot[blue,mark=*,mark size=1.5pt]  coordinates {

(3,0.0006094890591122885)
(4,0.00016683934213692272)
(5,3.928299179973519e-5)
(6,8.889122601308736e-6)
(7,2.1281208870984316e-6)
(8,5.417297372918646e-7)
(9,1.3283456593336577e-7)
(10,3.314440026434848e-8)
(11,8.090276509131133e-9)

};

\addplot[red,mark=*,mark size=1.5pt]  coordinates { 		
(3,1.5323663012361567e-5)
(4,1.3768793482992073e-6)
(5,1.6621205784969875e-7)
(6,1.8779439621595492e-8)
(7,2.258799666669327e-9)
(8,2.7850430322773536e-10)
(9,3.4188696385269416e-11)
(10,4.242800566005292e-12)

};

\addplot[blue,mark=none,dashed] coordinates {(3,0.1*2^-2*3) (11,0.1*2^-2*11)};
\addplot[red,mark=none,dashed] coordinates {(3, 0.01 * 2^-3*3 ) (10, 0.01   * 2^-3*10 )};

\end{semilogyaxis}
\end{tikzpicture} 
\subcaption{Decay of the RMSE. }
\end{minipage}
\hfill
\begin{minipage}[t]{0.45\linewidth}
\begin{tikzpicture}[scale=1]
\begin{semilogyaxis}[scale only axis,width = 0.7\textwidth,xlabel = {level $n$},
legend entries={$m=2$ ,$m=3$},
grid=major,
legend pos=outer north east,
]
\addplot[blue,mark=*,mark size=1.5pt]  coordinates { 		

(3,3.4341262495750016)
(4,5.425466807204421)
(5,4.66264808222643)
(6,4.402375623550835)
(7,5.074021781084119)
(8,6.064636127699306)
(9,4.8625509803600595)
(10,4.442035167716643)
(11,5.0705564700713115)

};

\addplot[red,mark=*,mark size=1.5pt]  coordinates { 		

(3,38.0168591871203)
(4,14.3937456893822)
(5,20.95670820347312)
(6,24.69691340678128)
(7,51.790000878152426)
(8,17.219913986167903)
(9,36.275072549611764)
(10,24.85888191202414)
(11,20.492869662545043)

};

\end{semilogyaxis}
\end{tikzpicture} \subcaption{Condition number of the wavelet matrix $\vec A$. }
\end{minipage}
\caption{Approximation on $[0,1]$ of the test function~\eqref{eq:f_cube} with uniformly distributed points using Chui-Wang wavelets of order $m\in \{2,3\}$.}
\label{fig:num_cube}
\end{figure}

\textbf{Extensions of non-periodic functions on the cube}\\
Here we want to demonstrate the benefits of the extension proposed in Section~\ref{sec:extension_cube}. Let us study the non-periodic 
function 
\begin{equation}\label{eq:f_cube} 
f\colon [0,1] \rightarrow \R, \quad f( y)= y^3.
\end{equation}
and the uniformly distribution 
on the cube, $\rho( y)=1$. Also for this non-periodic function we manage to use the periodic approximation operator and get good approximation results. 
We use a polynomial of degree $3$, since a lower degree together with the order of the wavelets $m$ leads to a function $f\circ\Rho^{-1}$, 
which is in the finite function space, which we use for the approximation and gives us approximation errors near machine precision. The results are plotted in Figure~\ref{fig:num_cube}. We see, that the extension increases the approximation rate to $2^{-mn}$, as proposed in Crollary~\ref{cor:nonper_end}
Additionally, the wavelet matrix is in both cases well-conditioned.    

 \section{Second Setting: Using kernel density estimation for unknown density \boldmath{$\rho$}}\label{sec:unknown_rho}
While working with real world data, the underlying density $\rho(\vec y)$ is possibly a priori not known, we only have the given random sample points $\Y$. 
For that reason we want adapt our strategy to this setting. A transformation of the given data $\Y$ is also in this case a 
useful tool to approximate a function $f$ well. Instead of using the transformation $\Rho$ belonging to the underlying density $\rho$ as in Chapter~\ref{sec:known_rho},
we approximate the underlying density function by a kernel density estimation, \cite{Gra18}. The cumulative distribution function $\r{\Rho}$ of the 
estimated density function gives us a transformation for the underlying data set $\Y$ to the samples $\r{\X} = \r{\Rho}(\Y)$ on the torus. Then we apply our approximation method 
for functions on the torus and at the end we transform the function back to a function defined on $\Omega$. We will restrict 
our study in this chapter to the one-dimensional case, which we apply in Chapter~\ref{sec:high-dim} to high-dimensional functions.
The error between the estimated density and the true density influences the total approximation error.\\

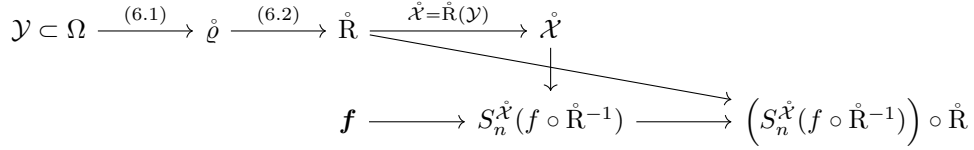
\begin{figure}[ht]
\begin{center}
\begin{tikzcd}[row sep=scriptsize, column sep=scriptsize]
\Y\subset \Omega \arrow[rr,"\eqref{eq:KDE}"]& & \r{\rho} \arrow[rr,"\eqref{eq:Rho_rho}"] & & \r{\Rho} \arrow[rr, "\r{\X} = \r{\Rho}(\Y)"] \arrow[rrrrd]&&\r{\X} \arrow[d]\\
&&&&\vec f \arrow[rr] &&S_n^{\r{\X}}(f\circ\r{\Rho}^{-1})\arrow[rr]&&\left(S_n^{\r{\X}}(f\circ\r{\Rho}^{-1})\right)\circ \r{\Rho}
\end{tikzcd}
\end{center}
\caption{Outline of our approximation procedure.} 
\label{fig:skizze}
\end{figure}

Let us introduce the kernel density estimator
\begin{equation}\label{eq:KDE}
\r{\rho}(y) =\sum_{x\in \X} \frac{1}{\sigma M}\,k\left(\frac{y-x}{\sigma}\right),
\end{equation}
where $k$ is a non-negative kernel function $k:\Omega\to \R$ which is normed by $\int_\Omega k(y)\d y =1$ and $\sigma$ is a smoothing 
parameter. Frequently used kernels 
are the standard normal distribution $\rho_N$, \eqref{eq:rho_N} or B-Splines~\eqref{eq:BSpline}. 
The normalization ensures that also $\r{\rho}$ is normalized. 
We get a transformation $\r{\Rho}:\Omega\to \T$, which fulfills~\eqref{eq:diff_eq}, by integration, 
\begin{align}\label{eq:Rho_rho}
\r{\Rho}(y)
&=\begin{cases}
\int_{0}^y \r{\rho}(t) \d t -\tfrac 12 &\text{if } \Omega =[0,1],\notag\\
\int_{-\infty}^y \r{\rho}(t) \d t -\tfrac 12 &\text{if } \Omega =\R,\\
\int_{-\sfrac 12}^y \r{\rho}(t) \d t -\tfrac 12 &\text{if } \Omega =\T.\notag\\
\end{cases}\\
&=\frac{1}{ M} \sum_{x\in \X}K(y-x) -\tfrac 12,
\end{align}
where $K$ is the antiderivative 
$$
K(y)=\begin{cases}
\int_{0}^y k(t) \d t  &\text{if } \Omega =[0,1],\\
\int_{-\infty}^y k(t) \d t  &\text{if } \Omega =\R,\\
\int_{-\sfrac 12}^y k(t) \d t  &\text{if } \Omega =\T.\\
\end{cases}
$$
Hence, the integral $K$ of the kernel $k$ has to be calculated once in advance. 
Then the calculation of the integral of $\r{\rho}$ is only the evaluation 
of $K$ at different points.
Using this transformation, we get our transformed sample points by 
$\r{\X} =\r{\Rho}(\Y) $. If $\r{\rho}$ is a good approximation to $\rho$, these samples $\r{\X}$ are nearly uniformly 
distributed on $\T$, which allows us to approximate the function $f\circ \r{\Rho}^{-1}$ using an approximation operator on $\T^d$, 
for instance $S_n^{\r{\X}}$ from~\eqref{eq:def_S_n^X}.
Our procedure to get an approximation of $f$ out of the samples 
$\X$ and the corresponding function values $\vec f$ is summarized in Figure~\ref{fig:skizze}. 
For shortening notation, we denote the error function $\r{e}_f\colon \Omega\to\R$ by
$$\r{e}_f=f-\left(S_n^{\r{\X}}(f\circ \r{\Rho}^{-1})\right)\circ\r{\Rho},$$
which we aim to estimate. 
Using the theory of the previous section, we receive a bound for 
$\norm{\r{e}_f}_{L_2(\Omega,\r{\rho})}$ in the norm with density $\r{\rho}$, 
if we assume that we choose the bandwidth $\sigma$ so that $\rho\approx \r{\rho}$, which yields that the samples $\r{\X}=\Rho(\Y)$ 
are distributed uniformly on $\T$. Since the original 
samples $\Y$ are distributed according to $\rho$ and we assume new test points are also samples according to $\rho$, we are interested in the 
$L_2(\Omega,\rho)$-error $\norm{\r{e}_f}_{L_2(\Omega,\rho)}$. Intuitively, if $\rho$ and $\r{\rho}$ are equal 
enough, these two errors have the same behavior. This can be made more precise by
\begin{theorem}\label{thm:bound_KDE}
In the case where $\Omega = \T$ we have that 
\begin{equation*}\norm{\r{e}_f}^2_{L_2(\Omega,\rho)}\leq  \norm{\r{e}_f}^2_{L_2(\T,\r{\rho})}+ \left(\int_{\T}  |\r{e}_f(y)|^4 \d y\right)^{1/2}\,\norm{\r{\rho}(y)-\rho(y)}_{L_2(\T)}.
\end{equation*}
In the case $\Omega=\R$ there is a set $\Eps\subset \R$, such that $\int_\Eps \rho(y)\dx y\leq \epsilon$ for some $\epsilon>0$. Therefore we have
\begin{equation}\label{eq:bound_KDE2}
\norm{\r{e}_f}^2_{L_2(\Omega,\rho)}\leq  \norm{\r{e}_f}^2_{L_2(\R,\r{\rho})}+ \left(\int_{\R\backslash \Eps}  |\r{e}_f(y)|^4 \d y\right)^{1/2}\,\norm{\r{\rho}(y)-\rho(y)}_{L_2(\R )} + \epsilon \sup_{y\in \Eps} |\r{e}_f|^2.
\end{equation}
\end{theorem}
\begin{proof}
The first inequality follows from triangle inequality and Cauchy-Schwarz inequality, 
\begin{align*}
\norm{\r{e}_f}^2_{L_2(\T,\rho)}
&=\int_{-\sfrac 12}^{\sfrac 12}|\r{e}_f(y)|^2 \rho(y)\d y 
=\int_{-\sfrac 12}^{\sfrac 12}|\r{e}_f(y)|^2 \r{\rho}(y)\left(1-\frac{\r{\rho}(y)-\rho(y)}{\r{\rho}(y)}\right)\d y\\
&\leq  \norm{\r{e}_f}^2_{L_2(\T,\r{\rho})}+ \int_{-\sfrac 12}^{\sfrac 12} |\r{e}_f(y)|^2 |\r{\rho}(y)-\rho(y)|\d y \\
&\leq  \norm{\r{e}_f}^2_{L_2(\T,\r{\rho})}+ \left(\int_{-\sfrac 12}^{\sfrac 12}  |\r{e}_f(y)|^4 \d y\right)^{1/2}\,\norm{\r{\rho}(y)-\rho(y)}_{L_2(\T)}.
\end{align*}
In the case where $\Omega = \R$ we can not use these calculation, since $\r{e}_f(y)$ can be non-zero on the whole real axis. 
But the splitting of $\R$ into the set $\Eps$ and the complement $\R\backslash \Eps$ and using the previous estimates, gives 
\begin{align*}
\norm{\r{e}_f}^2_{L_2(\R,\rho)}&= \int_{\R\backslash \Eps} |\r{e}_f(y)|^2 \rho(y)\dx y + \int_{ \Eps} |\r{e}_f(y)|^2 \rho(y)\dx y \\
&\leq  \norm{\r{e}_f}^2_{L_2(\R\backslash \Eps,\r{\rho})}+ \left(\int_{\R\backslash \Eps}  |\r{e}_f(y)|^4 \d y\right)^{1/2}\,\norm{\r{\rho}(y)-\rho(y)}_{L_2(\R\backslash \Eps)} + \epsilon \sup_{y\in \Eps} |\r{e}_f|^2\\
&\leq \norm{\r{e}_f}^2_{L_2(\R,\r{\rho})}+ \left(\int_{\R\backslash \Eps}  |\r{e}_f(y)|^4 \d y\right)^{1/2}\,\norm{\r{\rho}(y)-\rho(y)}_{L_2(\R )} + \epsilon \sup_{y\in \Eps} |\r{e}_f|^2.
\end{align*}
Therefore, it follows the assertion.
\end{proof}
The introduction of the set $\Eps$ in the case where $\Omega=\R$ tackles the behavior, that given data $\Y$ is contained 
in a finite interval $[a,b]$ and $\rho(y)$ is small outside this interval. In fact, we can not expect to approximate a function 
well where we have no information about the function.

The previous theorem shows that a good estimator for the bandwidth $\sigma$ ensures that the mean \textit{integrated squared error} ($\MISE$)
$$\MISE(\r{\rho}) = \E\left(\int_\Omega (\r{\rho}(y)-\rho(y))^2 \d y\right) = \E\left(\norm{\r{\rho}(y)-\rho(y)}_{L_2(\Omega)}^2\right).$$
is small. 
This choice of the smoothing parameter has to be a good trade-off between over- and under-fitting. Consider the two extremal cases, which do not work. This behavior is illustrated in Figure~\ref{fig:pic_fitting}.
\begin{itemize}
	\item If we would choose a small smoothing parameter $\sigma$ ($\sigma = 0.01$ in the Figure~\ref{fig:pic_fitting}) or even as kernel function $k$ the delta-distribution, we would get an equi-spaced sample set $\r{\Rho}(\Y)$. But this posses on the other hand no smooth cumulative distribution function $\r{\Rho}$ and the distribution $\r{\rho}$ is not a good approximation on $\rho$, since $\MISE(\r{\rho})$ would not decay.
	\item If we choose the parameter $\sigma$ in the Kernel function $k$ too big ($\sigma = 1$ in the Figure~\ref{fig:pic_fitting}), the density $\r{\rho}$ would not capture the behavior of the density $\rho$ and the transformed samples would not be uniformly distributed on $\T$.
\end{itemize}
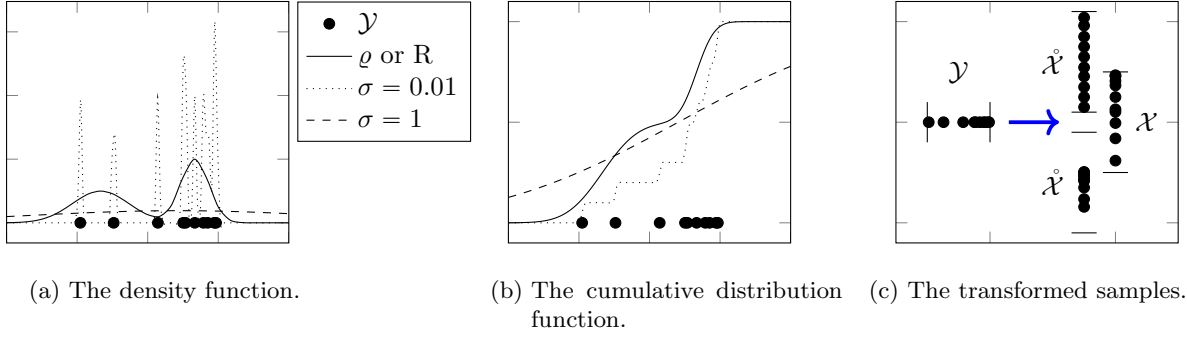
\begin{figure}[ht]
\centering
\begin{minipage}[t]{0.29\textwidth}
\centering
\begin{tikzpicture}[scale=1]
\begin{axis}[scale only axis, width = 0.8\textwidth,legend entries={ $\Y$ ,$\rho \text{ or } \Rho$ ,$\sigma=0.01$,$\sigma=1$},
xticklabel=\empty,yticklabel=\empty,
legend cell align={left},
legend pos=outer north east,
xmin=-1,
xmax=1,
]
\addplot[black,mark=*,only marks]  coordinates{ 		
(0.2683063433795556,0)
(-0.24101800485454633,0)
(0.3949788762882986,0)
(0.3336047093763347,0)
(0.4251451142305751,0)
(0.2516772258276241,0)
(0.4843401169681668,0)
(0.0722047095629883,0)
(0.471656195475935,0)
(-0.4772175315745213,0)

};
\addplot[mark=none,black,domain=-1:1,smooth] { \gauss{-1/3}{0.2}/2+\gauss{1/3}{0.1}/2};

\addplot[black,mark=*,mark=none,smooth,dotted]  coordinates{ 	

(-1.0,0.0)
(-0.9797979797979798,0.0)
(-0.9595959595959596,0.0)
(-0.9393939393939394,0.0)
(-0.9191919191919192,0.0)
(-0.898989898989899,0.0)
(-0.8787878787878788,0.0)
(-0.8585858585858586,5.9986333e-316)
(-0.8383838383838383,2.2472831679113756e-283)
(-0.8181818181818182,1.4217184076008678e-252)
(-0.797979797979798,1.5188668299423024e-223)
(-0.7777777777777778,2.740161328857049e-196)
(-0.7575757575757576,8.348008679878283e-171)
(-0.7373737373737373,4.2947706327558604e-147)
(-0.7171717171717171,3.731190143091719e-125)
(-0.696969696969697,5.4740040190252936e-105)
(-0.6767676767676768,1.3561669809496696e-86)
(-0.6565656565656566,5.673767619317033e-70)
(-0.6363636363636364,4.008484848784774e-55)
(-0.6161616161616161,4.782327041519576e-42)
(-0.5959595959595959,9.634931289132945e-31)
(-0.5757575757575758,3.2779950644863703e-21)
(-0.5555555555555556,1.8832947366130346e-13)
(-0.5353535353535354,1.827168560189075e-7)
(-0.5151515151515151,0.0029935689030507326)
(-0.494949494949495,0.8282286081809764)
(-0.47474747474747475,3.8695594221110627)
(-0.45454545454545453,0.30529742174353147)
(-0.43434343434343436,0.0004067570387683805)
(-0.41414141414141414,9.151606563269766e-9)
(-0.3939393939393939,3.477041939125452e-15)
(-0.37373737373737376,2.230862693308375e-23)
(-0.35353535353535354,1.2874675858605492e-27)
(-0.3333333333333333,1.245552084449469e-18)
(-0.31313131313131315,2.0348791464228387e-11)
(-0.29292929292929293,5.6139129382988975e-6)
(-0.2727272727272727,0.026154296880070455)
(-0.25252525252525254,2.05764641988361)
(-0.23232323232323232,2.7336869055705395)
(-0.21212121212121213,0.06133058679397168)
(-0.1919191919191919,2.323570163584553e-5)
(-0.1717171717171717,1.486568021719065e-10)
(-0.15151515151515152,1.606068651057432e-17)
(-0.13131313131313133,2.930176276787179e-26)
(-0.1111111111111111,9.02762923914535e-37)
(-0.09090909090909091,4.6968209606293354e-49)
(-0.0707070707070707,1.7833719955263923e-44)
(-0.050505050505050504,8.008700785129523e-33)
(-0.030303030303030304,6.073408845329647e-23)
(-0.010101010101010102,7.777737667441451e-15)
(0.010101010101010102,1.6819937898528872e-8)
(0.030303030303030304,0.000614250910522115)
(0.050505050505050504,0.3788066144004145)
(0.0707070707070707,3.9449329149686805)
(0.09090909090909091,0.6937644065681251)
(0.1111111111111111,0.0020603202495239523)
(0.13131313131313133,1.0332556078494887e-7)
(0.15151515151515152,8.750465123651279e-14)
(0.1717171717171717,5.216293230408014e-14)
(0.1919191919191919,7.023251940370631e-8)
(0.21212121212121213,0.0015973718275348076)
(0.23232323232323232,0.6192492531881352)
(0.25252525252525254,5.123575672647782)
(0.2727272727272727,4.053239708442099)
(0.29292929292929293,0.19429641657306426)
(0.31313131313131315,0.4907596580125648)
(0.3333333333333333,3.987954095050176)
(0.35353535353535354,0.5481813296865394)
(0.37373737373737376,0.4192400471794497)
(0.3939393939393939,3.9985710336023463)
(0.41414141414141414,2.813755593083431)
(0.43434343434343436,2.6187856079687406)
(0.45454545454545453,1.0229768476632908)
(0.47474747474747475,6.321529116103472)
(0.494949494949495,2.5371129446741048)
(0.5151515151515151,0.03494001437029257)
(0.5353535353535354,8.917407833479708e-6)
(0.5555555555555556,3.8726635868484656e-11)
(0.5757575757575758,2.8417597649278774e-18)
(0.5959595959595959,3.521559519282367e-27)
(0.6161616161616161,7.369441108547611e-38)
(0.6363636363636364,2.604261038332907e-50)
(0.6565656565656566,1.5541209235828803e-64)
(0.6767676767676768,1.5661574688809924e-80)
(0.696969696969697,2.66524028044432e-98)
(0.7171717171717171,7.659274499113115e-118)
(0.7373737373737373,3.716971442007581e-139)
(0.7575757575757576,3.0460788289088156e-162)
(0.7777777777777778,4.215444499993867e-187)
(0.797979797979798,9.851354270521805e-214)
(0.8181818181818182,3.887754409117089e-242)
(0.8383838383838383,2.590908236959619e-272)
(0.8585858585858586,2.9157854513905325e-304)
(0.8787878787878788,0.0)
(0.898989898989899,0.0)
(0.9191919191919192,0.0)
(0.9393939393939394,0.0)
(0.9595959595959596,0.0)
(0.9797979797979798,0.0)
(1.0,0.0)
};

\addplot[black,mark=*,mark=none,smooth,dashed]  coordinates{
(-1.0,0.1959276757680965)
(-0.9797979797979798,0.20025207046772525)
(-0.9595959595959596,0.20459839516107803)
(-0.9393939393939394,0.2089639295826562)
(-0.9191919191919192,0.21334587085204781)
(-0.898989898989899,0.21774133572591886)
(-0.8787878787878788,0.22214736305025157)
(-0.8585858585858586,0.22656091641209836)
(-0.8383838383838383,0.23097888698958452)
(-0.8181818181818182,0.23539809659835215)
(-0.797979797979798,0.23981530093208564)
(-0.7777777777777778,0.2442271929942042)
(-0.7575757575757576,0.24863040671724432)
(-0.7373737373737373,0.2530215207658943)
(-0.7171717171717171,0.2573970625190769)
(-0.696969696969697,0.2617535122259151)
(-0.6767676767676768,0.2660873073298547)
(-0.6565656565656566,0.2703948469546645)
(-0.6363636363636364,0.274672496545487)
(-0.6161616161616161,0.27891659265757257)
(-0.5959595959595959,0.2831234478848043)
(-0.5757575757575758,0.2872893559196052)
(-0.5555555555555556,0.2914105967353203)
(-0.5353535353535354,0.29548344188168285)
(-0.5151515151515151,0.2995041598835115)
(-0.494949494949495,0.3034690217323413)
(-0.47474747474747475,0.3073743064602724)
(-0.45454545454545453,0.3112163067849231)
(-0.43434343434343436,0.31499133481400793)
(-0.41414141414141414,0.3186957277977169)
(-0.3939393939393939,0.3223258539167632)
(-0.37373737373737376,0.32587811809368605)
(-0.35353535353535354,0.32934896781474693)
(-0.3333333333333333,0.3327348989495467)
(-0.31313131313131315,0.3360324615553107)
(-0.29292929292929293,0.33923826565264886)
(-0.2727272727272727,0.3423489869594949)
(-0.25252525252525254,0.34536137256986177)
(-0.23232323232323232,0.3482722465640253)
(-0.21212121212121213,0.3510785155367619)
(-0.1919191919191919,0.35377717403031955)
(-0.1717171717171717,0.3563653098588981)
(-0.15151515151515152,0.35884010931154653)
(-0.13131313131313133,0.36119886222056674)
(-0.1111111111111111,0.3634389668827271)
(-0.09090909090909091,0.36555793482084775)
(-0.0707070707070707,0.36755339537361925)
(-0.050505050505050504,0.3694231001018524)
(-0.030303030303030304,0.3711649269997344)
(-0.010101010101010102,0.3727768845000807)
(0.010101010101010102,0.37425711526302335)
(0.030303030303030304,0.37560389973806413)
(0.050505050505050504,0.37681565948994195)
(0.0707070707070707,0.37789096027931823)
(0.09090909090909091,0.3788285148898695)
(0.1111111111111111,0.37962718569399134)
(0.13131313131313133,0.38028598694996274)
(0.15151515151515152,0.38080408682408423)
(0.1717171717171717,0.38118080913199814)
(0.1919191919191919,0.3814156347941116)
(0.21212121212121213,0.38150820300077426)
(0.23232323232323232,0.38145831208361297)
(0.25252525252525254,0.3812659200901857)
(0.2727272727272727,0.38093114505989434)
(0.29292929292929293,0.3804542649998782)
(0.31313131313131315,0.3798357175603975)
(0.3333333333333333,0.3790760994100126)
(0.35353535353535354,0.37817616531165676)
(0.37373737373737376,0.3771368269014947)
(0.3939393939393939,0.37595915117324247)
(0.41414141414141414,0.37464435867141144)
(0.43434343434343436,0.37319382139770363)
(0.45454545454545453,0.37160906043554925)
(0.47474747474747475,0.36989174329851504)
(0.494949494949495,0.3680436810090425)
(0.5151515151515151,0.3660668249146764)
(0.5353535353535354,0.3639632632496292)
(0.5555555555555556,0.3617352174501827)
(0.5757575757575758,0.35938503823306084)
(0.5959595959595959,0.35691520144650607)
(0.6161616161616161,0.3543283037043662)
(0.6363636363636364,0.3516270578140297)
(0.6565656565656566,0.3488142880095566)
(0.6767676767676768,0.34589292500180996)
(0.696969696969697,0.34286600085782815)
(0.7171717171717171,0.33973664372206)
(0.7373737373737373,0.3365080723924386)
(0.7575757575757576,0.33318359076457116)
(0.7777777777777778,0.32976658215759075)
(0.797979797979798,0.32626050353543445)
(0.8181818181818182,0.32266887963749213)
(0.8383838383838383,0.31899529703270313)
(0.8585858585858586,0.3152433981112696)
(0.8787878787878788,0.31141687502820203)
(0.898989898989899,0.3075194636129126)
(0.9191919191919192,0.3035549372590362)
(0.9393939393939394,0.29952710080857065)
(0.9595959595959596,0.2954397844443052)
(0.9797979797979798,0.2912968376043386)
(1.0,0.2871021229322797)

};

\end{axis}
\end{tikzpicture} 
\subcaption{The density function.}
\end{minipage}
\hfill
\begin{minipage}[t]{0.29\textwidth}
\begin{tikzpicture}[scale=1]
\begin{axis}[scale only axis,width = 0.8\textwidth,xticklabel=\empty,yticklabel=\empty,
xmin=-1,
xmax=1,
]
\addplot[black,mark=*,only marks]  coordinates{ 		
(0.2683063433795556,0)
(-0.24101800485454633,0)
(0.3949788762882986,0)
(0.3336047093763347,0)
(0.4251451142305751,0)
(0.2516772258276241,0)
(0.4843401169681668,0)
(0.0722047095629883,0)
(0.471656195475935,0)
(-0.4772175315745213,0)
};
\addplot[black,mark=*,mark=none,smooth]  coordinates{

(-1.0,0.00021453016659841868)
(-0.9797979797979798,0.00030697019517622446)
(-0.9595959595959596,0.0004350662449444111)
(-0.9393939393939394,0.0006107712067901248)
(-0.9191919191919192,0.000849335223703457)
(-0.898989898989899,0.0011699610315230433)
(-0.8787878787878788,0.0015965058206767642)
(-0.8585858585858586,0.00215820453818167)
(-0.8383838383838383,0.0028903794644278024)
(-0.8181818181818182,0.003835090310370476)
(-0.797979797979798,0.0050416688066311326)
(-0.7777777777777778,0.0065670728455105585)
(-0.7575757575757576,0.008475988920195748)
(-0.7373737373737373,0.010840609207509828)
(-0.7171717171717171,0.013740012457608489)
(-0.696969696969697,0.0172590869986038)
(-0.6767676767676768,0.021486950391571957)
(-0.6565656565656566,0.026514843806310577)
(-0.6363636363636364,0.032433509573667224)
(-0.6161616161616161,0.039330096365986646)
(-0.5959595959595959,0.0472846760106365)
(-0.5757575757575758,0.05636649625112677)
(-0.5555555555555556,0.06663013145125268)
(-0.5353535353535354,0.07811172458748533)
(-0.5151515151515151,0.09082553522172446)
(-0.494949494949495,0.10476101624909707)
(-0.47474747474747475,0.11988063467006782)
(-0.45454545454545453,0.13611862724075355)
(-0.43434343434343436,0.15338084087841347)
(-0.41414141414141414,0.171545752007488)
(-0.3939393939393939,0.19046669204715003)
(-0.37373737373737376,0.20997523271499308)
(-0.35353535353535354,0.2298856104375768)
(-0.3333333333333333,0.250000000006542)
(-0.31313131313131315,0.270114389589427)
(-0.29292929292929293,0.290024767380029)
(-0.2727272727272727,0.309533308291962)
(-0.25252525252525254,0.3284542491595944)
(-0.23232323232323232,0.34661916298238854)
(-0.21212121212121213,0.363881385033816)
(-0.1919191919191919,0.380119402836904)
(-0.1717171717171717,0.395239093911768)
(-0.15151515151515152,0.40917477579835343)
(-0.13131313131313133,0.4218891195924413)
(-0.1111111111111111,0.43337207153059853)
(-0.09090909090909091,0.4436390317037229)
(-0.0707070707070707,0.4527286637880177)
(-0.050505050505050504,0.4607008653287758)
(-0.030303030303030304,0.46763561890524236)
(-0.010101010101010102,0.47363365396878077)
(0.010101010101010102,0.4788200198036042)
(0.030303030303030304,0.48335168420818636)
(0.050505050505050504,0.4874299485739146)
(0.0707070707070707,0.4913175953306718)
(0.09090909090909091,0.49535910139017475)
(0.1111111111111111,0.5)
(0.13131313131313133,0.5057989404008787)
(0.15151515151515152,0.5134239966882334)
(0.1717171717171717,0.5236244643418828)
(0.1919191919191919,0.537171891827805)
(0.21212121212121213,0.55476999043045)
(0.23232323232323232,0.5769417635559623)
(0.25252525252525254,0.6039116810253937)
(0.2727272727272727,0.6355078560339625)
(0.29292929292929293,0.6711106857625243)
(0.31313131313131315,0.7096682625194316)
(0.3333333333333333,0.7497854698334017)
(0.35353535353535354,0.789876269510301)
(0.37373737373737376,0.8283524405107412)
(0.3939393939393939,0.8638122442803378)
(0.41414141414141414,0.8951924956466164)
(0.43434343434343436,0.9218573137177524)
(0.45454545454545453,0.9436130818228025)
(0.47474747474747475,0.9606565638353592)
(0.494949494949495,0.9734765269166823)
(0.5151515151515151,0.9827353850465466)
(0.5353535353535354,0.9891558839821617)
(0.5555555555555556,0.9934307241726381)
(0.5757575757575758,0.9961635392689663)
(0.5959595959595959,0.9978409512818917)
(0.6161616161616161,0.9988295240382354)
(0.6363636363636364,0.9993889177731319)
(0.6565656565656566,0.9996928437890527)
(0.6767676767676768,0.9998513920640436)
(0.696969696969697,0.9999308069235042)
(0.7171717171717171,0.9999690007982659)
(0.7373737373737373,0.9999866386383025)
(0.7575757575757576,0.9999944597705817)
(0.7777777777777778,0.9999977900998942)
(0.797979797979798,0.9999991519592754)
(0.8181818181818182,0.9999996868466601)
(0.8383838383838383,0.9999998886720718)
(0.8585858585858586,0.9999999618608724)
(0.8787878787878788,0.9999999873864073)
(0.898989898989899,0.9999999959591865)
(0.9191919191919192,0.9999999987383001)
(0.9393939393939394,0.9999999996117179)
(0.9595959595959596,0.9999999998799792)
(0.9797979797979798,0.999999999961665)
(1.0,0.999999999986916)
};

\addplot[black,mark=*,mark=none,smooth,dotted]  coordinates{
(-1.0,0.0)
(-0.9797979797979798,0.0)
(-0.9595959595959596,0.0)
(-0.9393939393939394,0.0)
(-0.9191919191919192,0.0)
(-0.898989898989899,0.0)
(-0.8787878787878788,0.0)
(-0.8585858585858586,1.57187e-319)
(-0.8383838383838383,6.217535354640553e-287)
(-0.8181818181818182,4.1661214735831846e-256)
(-0.797979797979798,4.730590363700891e-227)
(-0.7777777777777778,9.106786734149844e-200)
(-0.7575757575757576,2.973848245244954e-174)
(-0.7373737373737373,1.6484146349772548e-150)
(-0.7171717171717171,1.5522727281220667e-128)
(-0.696969696969697,2.4858632266799565e-108)
(-0.6767676767676768,6.779181243815333e-90)
(-0.6565656565656566,3.153805189154897e-73)
(-0.6363636363636364,2.5089160015984675e-58)
(-0.6161616161616161,3.4243490735217175e-45)
(-0.5959595959595959,8.05780305094111e-34)
(-0.5757575757575758,3.2933102958638077e-24)
(-0.5555555555555556,2.3666624427758405e-16)
(-0.5353535353535354,3.057163988611359e-10)
(-0.5151515151515151,7.429973698161159e-6)
(-0.494949494949495,0.0038098089371470635)
(-0.47474747474747475,0.059754808692097316)
(-0.45454545454545453,0.09883112297503396)
(-0.43434343434343436,0.09999909615657578)
(-0.41414141414141414,0.09999999998583131)
(-0.3939393939393939,0.1)
(-0.37373737373737376,0.1)
(-0.35353535353535354,0.1)
(-0.3333333333333333,0.1)
(-0.31313131313131315,0.1000000000000277)
(-0.29292929292929293,0.10000001045114774)
(-0.2727272727272727,0.10007597671762593)
(-0.25252525252525254,0.11249227420077489)
(-0.23232323232323232,0.18077069277467303)
(-0.21212121212121213,0.19980718247436197)
(-0.1919191919191919,0.1999999544342293)
(-0.1717171717171717,0.1999999999997897)
(-0.15151515151515152,0.2)
(-0.13131313131313133,0.2)
(-0.1111111111111111,0.2)
(-0.09090909090909091,0.2)
(-0.0707070707070707,0.2)
(-0.050505050505050504,0.2)
(-0.030303030303030304,0.2)
(-0.010101010101010102,0.2)
(0.010101010101010102,0.20000000002643004)
(0.030303030303030304,0.20000139374042275)
(0.050505050505050504,0.2015004714437041)
(0.0707070707070707,0.24404754515321092)
(0.09090909090909091,0.296928849974766)
(0.1111111111111111,0.299995000994577)
(0.13131313131313133,0.2999999998298185)
(0.15151515151515152,0.29999999999999993)
(0.1717171717171717,0.30000000000000004)
(0.1919191919191919,0.3000000001144806)
(0.21212121212121213,0.3000038180814235)
(0.23232323232323232,0.30266307996603603)
(0.25252525252525254,0.3591061140154175)
(0.2727272727272727,0.4653143572262747)
(0.29292929292929293,0.499310268156153)
(0.31313131313131315,0.5020308685817344)
(0.3333333333333333,0.5489174991265012)
(0.35353535353535354,0.5976889841833125)
(0.37373737373737376,0.6016798985343625)
(0.3939393939393939,0.6459507714192051)
(0.41414141414141414,0.7107919074404033)
(0.43434343434343436,0.7821223920207362)
(0.45454545454545453,0.8043336590910921)
(0.47474747474747475,0.8790100367131767)
(0.494949494949495,0.9845720118383806)
(0.5151515151515151,0.9998962129616848)
(0.5353535353535354,0.9999999831278427)
(0.5555555555555556,0.9999999999999466)
(0.5757575757575758,0.9999999999999999)
(0.5959595959595959,0.9999999999999999)
(0.6161616161616161,0.9999999999999999)
(0.6363636363636364,0.9999999999999999)
(0.6565656565656566,0.9999999999999999)
(0.6767676767676768,0.9999999999999999)
(0.696969696969697,0.9999999999999999)
(0.7171717171717171,0.9999999999999999)
(0.7373737373737373,0.9999999999999999)
(0.7575757575757576,0.9999999999999999)
(0.7777777777777778,0.9999999999999999)
(0.797979797979798,0.9999999999999999)
(0.8181818181818182,0.9999999999999999)
(0.8383838383838383,0.9999999999999999)
(0.8585858585858586,0.9999999999999999)
(0.8787878787878788,0.9999999999999999)
(0.898989898989899,0.9999999999999999)
(0.9191919191919192,0.9999999999999999)
(0.9393939393939394,0.9999999999999999)
(0.9595959595959596,0.9999999999999999)
(0.9797979797979798,0.9999999999999999)
(1.0,0.9999999999999999)
};

\addplot[black,mark=*,mark=none,smooth,dashed]  coordinates{
(-1.0,0.12631482498597957)
(-0.9797979797979798,0.13031660144039411)
(-0.9595959595959596,0.13440596542788588)
(-0.9393939393939394,0.138583332642073)
(-0.9191919191919192,0.1428490629835984)
(-0.898989898989899,0.14720345891355166)
(-0.8787878787878788,0.15164676385440912)
(-0.8585858585858586,0.15617916064253104)
(-0.8383838383838383,0.16080077003623464)
(-0.8181818181818182,0.16551164928343076)
(-0.797979797979798,0.17031179075276992)
(-0.7777777777777778,0.1752011206321894)
(-0.7575757575757576,0.18017949769868993)
(-0.7373737373737373,0.1852467121630937)
(-0.7171717171717171,0.1904024845934479)
(-0.696969696969697,0.19564646492063995)
(-0.6767676767676768,0.20097823152968042)
(-0.6565656565656566,0.20639729043998814)
(-0.6363636363636364,0.21190307457787944)
(-0.6161616161616161,0.21749494314432097)
(-0.5959595959595959,0.22317218108084994)
(-0.5757575757575758,0.2289339986364027)
(-0.5555555555555556,0.23477953103761584)
(-0.5353535353535354,0.24070783826498032)
(-0.5151515151515151,0.2467179049370345)
(-0.494949494949495,0.2528086403045785)
(-0.47474747474747475,0.25897887835667893)
(-0.45454545454545453,0.26522737804001423)
(-0.43434343434343436,0.2715528235928812)
(-0.41414141414141414,0.2779538249949486)
(-0.3939393939393939,0.2844289185336014)
(-0.37373737373737376,0.2909765674874715)
(-0.35353535353535354,0.29759516292749644)
(-0.3333333333333333,0.30428302463559176)
(-0.31313131313131315,0.31103840214075773)
(-0.29292929292929293,0.31785947587217966)
(-0.2727272727272727,0.3247443584286091)
(-0.25252525252525254,0.33169109596304797)
(-0.23232323232323232,0.33869766968148524)
(-0.21212121212121213,0.3457619974541639)
(-0.1919191919191919,0.3528819355375919)
(-0.1717171717171717,0.36005528040523405)
(-0.15151515151515152,0.3672797706845645)
(-0.13131313131313133,0.3745530891978909)
(-0.1111111111111111,0.38187286510410434)
(-0.09090909090909091,0.3892366761382572)
(-0.0707070707070707,0.3966420509456187)
(-0.050505050505050504,0.40408647150662025)
(-0.030303030303030304,0.41156737564886414)
(-0.010101010101010102,0.419082159642147)
(0.010101010101010102,0.4266281808722255)
(0.030303030303030304,0.43420276058885016)
(0.050505050505050504,0.4418031867233884)
(0.0707070707070707,0.44942671677117574)
(0.09090909090909091,0.4570705807335549)
(0.1111111111111111,0.4647319841144002)
(0.13131313131313133,0.4724081109657733)
(0.15151515151515152,0.48009612697721826)
(0.1717171717171717,0.4877931826030818)
(0.1919191919191919,0.4954964162221298)
(0.21212121212121213,0.503202957323643)
(0.23232323232323232,0.5109099297140878)
(0.25252525252525254,0.5186144547383983)
(0.2727272727272727,0.5263136545098508)
(0.29292929292929293,0.5340046551424857)
(0.31313131313131315,0.5416845899800057)
(0.3333333333333333,0.5493506028150816)
(0.35353535353535354,0.5569998510930145)
(0.37373737373737376,0.5646295090937241)
(0.3939393939393939,0.5722367710860878)
(0.41414141414141414,0.5798188544487153)
(0.43434343434343436,0.5873730027513142)
(0.45454545454545453,0.5948964887909085)
(0.47474747474747475,0.6023866175772653)
(0.494949494949495,0.6098407292620205)
(0.5151515151515151,0.6172562020061239)
(0.5353535353535354,0.62463045478038)
(0.5555555555555556,0.6319609500940239)
(0.5757575757575758,0.6392451966464507)
(0.5959595959595959,0.6464807518974082)
(0.6161616161616161,0.6536652245511615)
(0.6363636363636364,0.660796276950358)
(0.6565656565656566,0.6678716273755396)
(0.6767676767676768,0.6748890522464853)
(0.696969696969697,0.6818463882218121)
(0.7171717171717171,0.688741534193508)
(0.7373737373737373,0.6955724531733369)
(0.7575757575757576,0.7023371740683129)
(0.7777777777777778,0.7090337933427183)
(0.797979797979798,0.7156604765644132)
(0.8181818181818182,0.7222154598334629)
(0.8383838383838383,0.7286970510913977)
(0.8585858585858586,0.7351036313096961)
(0.8787878787878788,0.7414336555563815)
(0.898989898989899,0.7476856539398969)
(0.9191919191919192,0.7538582324297192)
(0.9393939393939394,0.7599500735534527)
(0.9595959595959596,0.7659599369704323)
(0.9797979797979798,0.7718866599221406)
(1.0,0.7777291575600217)

};

\end{axis}
\end{tikzpicture} 
\subcaption{The cumulative distribution function. }
\end{minipage}
\begin{minipage}[t]{0.29\textwidth}
\centering
\begin{tikzpicture}[scale=1]
\begin{axis}[scale only axis,width = 0.8\textwidth,xticklabel=\empty,yticklabel=\empty,
xmin=-3.5,
xmax=1,
ymax =3.2,
ymin = 0.8,
]

\node at (axis cs:-1,2.6) {$\r{\X}$};
\node at (axis cs:-1,1.4) {$\r{\X}$};
\node at (axis cs:-2.5,2.5) {$\Y$};
\node at (axis cs:0.5,2) {$\X$};
\draw[->,line width=0.5mm,blue]  (axis cs:-1.7,2) --  (axis cs:-0.9,2);
\addplot[black,mark=*,only marks]  coordinates{

(-0.5,0.02462939350711979+1.4)
(-0.5,-0.16432510072986817+1.4)
(-0.5,0.07262754059098453+1.4)
(-0.5,0.049453473473094034+1.4)
(-0.5,0.08393711061059439+1.4)
(-0.5,0.018291126402634417+1.4)
(-0.5,0.10593074895018362+1.4)
(-0.5,-0.05000728335205201+1.4)
(-0.5,0.1012427587765079+1.4)
(-0.5,-0.24177976822919817+1.4)

};

\addplot[black,mark=*,only marks]  coordinates{  (-0.5,-0.054816505045317654+2.6)
(-0.5,-0.35+2.6)
(-0.5,0.15012780336976594+2.6)
(-0.5,0.05000000003864502+2.6)
(-0.5,0.24987236182783845+2.6)
(-0.5,-0.14518349495139027+2.6)
(-0.5,0.43976710264652097+2.6)
(-0.5,-0.25+2.6)
(-0.5,0.36023273211393747+2.6)
(-0.5,-0.45+2.6)
};
\addplot[black,mark=*,only marks]  coordinates{            (0.2683063433795556-2.5,2)
(-0.24101800485454633-2.5,2)
(0.3949788762882986-2.5,2)
(0.3336047093763347-2.5,2)
(0.4251451142305751-2.5,2)
(0.2516772258276241-2.5,2)
(0.4843401169681668-2.5,2)
(0.0722047095629883-2.5,2)
(0.471656195475935-2.5,2)
(-0.4772175315745213-2.5,2) };

\addplot[black,mark=*,only marks]  coordinates{        

(0,0.12822247037247592+2)
(0,-0.16109623467035983+2)
(0,0.3655337576296037+2)
(0,0.25032783003708703+2)
(0,0.4103234012005146+2)
(0,0.10268383236085121+2)
(0,0.4672326106139879+2)
(0,-0.008392909753772437+2)
(0,0.458337058912901+2)
(0,-0.3820296082418867+2)
};

\draw  (axis cs:-2,1.8) -- (axis cs:-2,2.2);
\draw  (axis cs:-3,1.8) -- (axis cs:-3,2.2);
\draw  (axis cs:-0.7,2.1) -- (axis cs:-0.3,2.1);
\draw  (axis cs:-0.7,3.1) -- (axis cs:-0.3,3.1);
\draw  (axis cs:-0.7,1.9) -- (axis cs:-0.3,1.9);
\draw  (axis cs:-0.7,0.9) -- (axis cs:-0.3,0.9);
\draw  (axis cs:-0.2,1.5) -- (axis cs:0.2,1.5);
\draw  (axis cs:-0.2,2.5) -- (axis cs:0.2,2.5);
\end{axis}
\end{tikzpicture}

\subcaption{The transformed samples. }
\end{minipage}
\caption{Illustration of the problem of over- and underfitting.}
\label{fig:pic_fitting}
\end{figure} 

\subsection{Smoothing parameter selection}
There are some very simple and easy to compute mathematical formulas for estimating the smoothing parameter $\sigma$. They are
often called the \textit{rules-of-thumb} (ROT). One possibility is
\begin{equation}\label{eq:ROT}
\sigma_{\text{ROT}} = 1.06 \,\min \left\{\text{std}(\Y),\frac{\text{IQR}(\Y)}{1.34}\right\} \, M^{-1/5}, 
\end{equation}
where $\text{std}$ is the \textit{standard deviation} and $\text{IQR}$ is the \textit{inter-quartile range}, see \cite[Section 4.2.1]{Gra18}.
The assumption for that rule is that the unknown density belongs to the family of the normal distribution. In practice, we do not know, whether $\rho(y)$ is a normal distribution. If it is, then $\sigma_{\text{ROT}}$ gives the optimal smoothing parameter. If not, then $\sigma_{\text{ROT}}$ will give a parameter not too far from the optimum, if the distribution of the samples $\Y$ is not too different from the normal distribution. 

\subsubsection{Data on the real axis}\label{sec:KDE_R}
Another approach, which is more general and performs better than the ROT is the \textit{Direct Plug-In-Selector} (DPI), see \cite[Section 4.2.2]{Gra18}. 
To describe this approach, we have to introduce some notation. The second moment of the kernel $k$ is defined by
$$\mu_2(k) = \int_{\R} y^2 k(y) \d y.$$
Since the MISE is set as the error criterion to be minimized, our aim is to find 
$$\sigma_{\MISE} := \argmin_{\sigma>0} \MISE(\r{\rho}).$$
The dominating part of $\MISE$ is denoted by $\AMISE$, which stands for \textit{Asymptotic MISE},
$$\AMISE(\r{\rho}) = \frac 14 \mu_2^2(k) \norm{\rho''}^2_{L_2(\R)}\sigma^4 +\frac{\norm{k}^2_{L_2(\R)}}{M\sigma},$$
The minimizer $\sigma_{\AMISE}$ is given by
\begin{equation}\label{eq:sigma_AMISE}
\sigma_{\AMISE}=\left(\frac{\norm{k}^2_{L_2(\R )}}{\mu_2^2(k)\norm{\rho''}^2_{L_2(\R )}M}\right)^{1/5} = \left(\frac{\norm{k}^2_{L_2(\R )}}{\mu_2^2(k)\Psi_4M}\right)^{1/5},
\end{equation}
where 
$$\Psi_4 = \int_\R \rho^{(4)}(y)\rho(y)\d y,\quad  \text{or more generally } \Psi_r = \int_\R \rho^{(r)}(y)\rho(y)\d y,$$
where $r$ is an even number. The naming convention of $\Psi_r$ was introduced in~\cite[Section 3.5]{WaJo95}. 
The critical step is to estimate $\Psi_4$ in~\eqref{eq:sigma_AMISE}, as this is the only unknown value. Assuming that $\rho $ 
is some normal distribution, would lead to~\eqref{eq:ROT}. This is an example of a zero-stage plug-in selector, a terminology 
inspired by the fact that $\Psi_4$ was estimated by directly plugging in a parametric assumption. Another possibility is 
to estimate $\Psi_4$ non-parametrically then to plug it into $\sigma_{\AMISE}$. First note, that integration by parts gives
$$ \norm{\rho^{(r)}}_{L^2(\R)}^2 = (-1)^r \int_{\R} \rho^{(2r)}(y)\rho(y)\d y.$$ 
Therefore, a possible way to estimate $\Psi_r$ is
\begin{equation}\label{eq:Psi_r}
\overline{\Psi_r}=\frac{1}{M^2 g^{r+1}} \sum_{i=1}^M\sum_{j=1}^M k^{(r)}\left(\frac{y_i-y_j}{g}\right),
\end{equation}
where $g$ is the smoothing parameter of a kernel density estimation. Typically, two stages are considered to have a good trade-off between 
bias and variance. This is the method proposed by~\cite{ShJo91}, and does the following steps.
\begin{enumerate}
	\item Estimate $\Psi_8$ by $\overline{\Psi}_8=\left(\frac{105}{32\sqrt{\pi} (\overline{\text{std}}(\rho))^9}\right)$, where $\overline{\text{std}}(\rho)$ is an estimate for the standard derivation of $\rho$, which can 
	be $\text{std}(\X)$ or $\min \left\{\text{std}(\X),\frac{\text{IQR}(\X)}{1.34}\right\}$.
	\item Estimate $\Psi_6$ using $\overline{\Psi}_6$ from~\eqref{eq:Psi_r}, where 
	$g_1=\left(-\frac{2k^{(6)}(0)}{\mu_2(k)\overline{\Psi}_8M}\right)^{1/9}.$
	\item Estimate $\Psi_4$ using $\overline{\Psi}_4$ from~\eqref{eq:Psi_r}, where
	$g_2=\left(-\frac{2k^{(4)}(0)}{\mu_2(k)\overline{\Psi}_6M}\right)^{1/7}.$
	\item The selected smoothing parameter is 
	$\sigma_{\text{DPI}} := \left(\frac{\norm{k}_{L_2(\R)}^2}{\mu_2^2(k)\overline{\Psi}_4M}\right)^{1/5}.$
\end{enumerate}

\subsection{Non-periodic data }\label{sec:KDE_cube}
A general problem with kernel density estimation is that certain difficulties can arise at the boundaries
and near them. In many practical situations the values of a random variable are
bounded. For example, the age of a person obviously cannot be a negative number.
On the other hand, the normal kernel has unlimited support. Even if a kernel with finite support
is used, the estimated density can usually go beyond the permissible
domain. \\
For that reason we use on $\Omega =[0,1]$ compactly supported kernels $k$, but we allow the approximated density $\r{\rho}$ 
to be nonzero outside the interval $[0,1]$. Especially, we receive a nonzero density in the interval $\r{\Omega}=[\omega_1,\omega_2]$ with 
$-|\supp k|\,\tfrac{\sigma}{2}\leq\omega_1\leq 0$ and $1\leq\omega_2\leq 1+|\supp k|\,\tfrac{\sigma}{2}$. 
This allows us to create a function $\r{f}$, which smoothly extends the function $f$ to a function on whole $\r{\Omega}$,  
such that the transformed function $f\circ\r{\Rho}^{-1}$ becomes a periodic function. 
This idea of an extension of the function is similar to the studies in Section~\ref{sec:extension_cube}. The choice of the extension parameter $\eta$ is now hidden in the choice of the smoothing parameter selection $\sigma$, which determines the interval $\r{\Omega}$.
This is illustrated in Figure~\ref{fig:pic_boundary_cube}. 
The kernel density estimation of $\rho$ can be seen as a periodization of the function $f$, 
such that we can use approximation operators for functions on $\T$. In contrast to the tent transformation, see~\cite{nasdaladiss}, which passes trough the function forth and back, we have here no need to double the number of sample points.\\  
\begin{figure}[ht]
\centering
\begin{subfigure}[c]{0.49\textwidth}
\centering
\begin{tikzpicture}[scale=1]
\begin{axis}[scale only axis, width = 0.6\textwidth,legend entries={ $f$ ,$\r{f}$},
xtick={-0.1,0,1,1.1},xticklabels={$\omega_1$,$0$,$1$,$\omega_2$},yticklabel=\empty,
legend pos=outer north east,
]
\addplot[black,domain = 0:1] {x} ;
\addplot[black,dashed, smooth] coordinates {(-0.1,0.5)(0,0)(0.1,0.1)(0.2,0.2)(0.3,0.3)(0.4,0.4)(0.5,0.5)(0.6,0.6)(0.7,0.7)(0.8,0.8)(0.9,0.9)(1,1)(1.1,0.5)};	

\addplot[dotted, samples=50, smooth,domain=0:1,black] coordinates {(1,-0.3)(1,1.3)};
\addplot[dotted, samples=50, smooth,domain=0:1,black] coordinates {(0,-0.3)(0,1.3)};
\addplot[dotted, samples=50, smooth,domain=0:1,black] coordinates {(-0.1,-0.3)(-0.1,1.3)};
\addplot[dotted, samples=50, smooth,domain=0:1,black] coordinates {(1.1,-0.3)(1.1,1.3)};
\end{axis}

\end{tikzpicture} 
\subcaption{Function $f$ and extended function $\r{f}$ on $\r{\Omega}$. }
\end{subfigure}
\begin{subfigure}[c]{0.49\textwidth}
\centering
\begin{tikzpicture}[scale=1]
\begin{axis}[scale only axis,width = 0.6\textwidth,legend entries={$\rho$,$\r{\rho}$},
xtick={-0.1,0,1,1.1},xticklabels={$\omega_1$,$0$,$1$,$\omega_2$},yticklabel=\empty,
legend pos=outer north east,
]

\addplot[black,mark=none,domain = 0:1] {1} 	;
\addplot[black,dashed, smooth] coordinates {(-0.1,0)(0,0.5)(0.1,1)(0.2,1.2)(0.3,0.9)(0.4,0.97)(0.5,1.1)(0.6,1.2)(0.7,1)(0.8,0.9)(0.9,1.2)(1,0.5)(1.1,0)};
\addplot[dotted, samples=50, smooth,domain=0:1,black] coordinates {(1,-0.3)(1,1.3)};
\addplot[dotted, samples=50, smooth,domain=0:1,black] coordinates {(0,-0.3)(0,1.3)};
\addplot[dotted, samples=50, smooth,domain=0:1,black] coordinates {(-0.1,-0.3)(-0.1,1.3)};
\addplot[dotted, samples=50, smooth,domain=0:1,black] coordinates {(1.1,-0.3)(1.1,1.3)};
\end{axis}
\end{tikzpicture} 
\subcaption{The density $\rho$ and the estimated density $\r{\rho}$.}
\end{subfigure}
\caption{Periodization of $f$ using instead of the real density $\rho$ the estimated density $\r{\rho}$. }
\label{fig:pic_boundary_cube}
\end{figure}
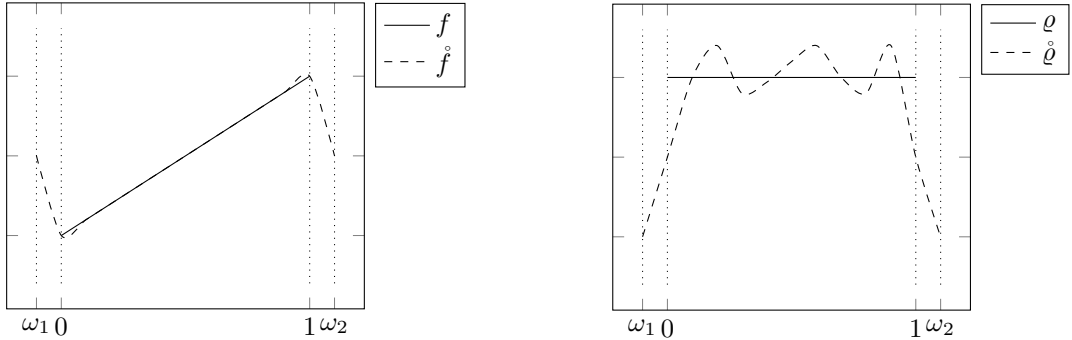 To select the smoothing parameter $\sigma$, one simple possibility is to again use the estimator $\sigma_{\text{ROT}}$ from~\eqref{eq:ROT}. 
Analogously to Theorem~\ref{thm:bound_KDE} we give an estimate for the error decay, namely,
\begin{align*}
\norm{\r{e}_f}^2_{L_2([0,1],\rho)} = \int_{0}^1 |\r{e}_f(y)|^2 \frac{\rho(y)}{\r{\rho}(y)} \r{\rho}(y)\dx y\leq \left(\max_{y\in [0,1]} \frac{\rho(y)}{\r{\rho}(y)} \right)\norm{\r{e}_f}^2_{L_2([0,1],\r{\rho})}.
\end{align*}
The extension $\r{f}$ is similar to the proposed method in Section~\ref{sec:extension_cube}. But here naturally, the extension is on both sides of the interval because of the kernel density estimation. Instead of the factor $\tfrac{1}{1-\eta}$ we have now the the term  $\max_{y\in [0,1]} \tfrac{\rho(y)}{\r{\rho}(y)} $. 

\subsection{Numerical experiments}In this section we endorse our theoretical findings by two numerical experiments on $\R$ and $[0,1]$. We compare the approximation results of unknown density with the results if the density $\rho$ is known. In both cases we use Chui-Wang wavelets of order $m=3$.\\

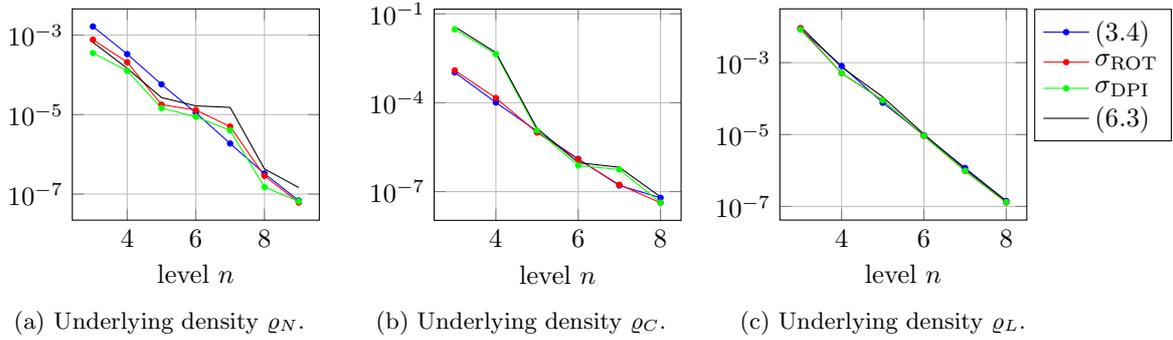
\begin{figure}[ht]
\begin{minipage}[t]{.29\linewidth}
\centering
\begin{tikzpicture}[scale=1]
\begin{semilogyaxis}[scale only axis, width = 0.7\textwidth,xlabel = {level $n$},
grid=major,
]
\addplot[blue,mark=*,mark size=1.0pt]  coordinates { 		
(3,0.0016466185575652217)
(4,0.00033328989188780337)
(5,5.769500626777274e-5)
(6,1.1181500415555705e-5)
(7,1.8691505122131892e-6)
(8,3.2895616497605937e-7)
(9,6.864341386377416e-8)
};

\addplot[red,mark=*,mark size=1.0pt]  coordinates {

(3,0.0007673242132792716)
(4,0.00020766246143368275)
(5,1.7759776803015017e-5)
(6,1.3025968947183485e-5)
(7,5.002664834948591e-6)
(8,2.877271195109052e-7)
(9,6.125588774264284e-8)
};

\addplot[green,mark=*,mark size=1.0pt]  coordinates {
(3,0.00035470355547357815)
(4,0.00012559733332303153)
(5,1.4621857846176841e-5)
(6,8.820825741499063e-6)
(7,4.001101697254171e-6)
(8,1.5025732299876296e-7)
(9,6.578995041517759e-8)
};

\addplot[black,mark=none]  coordinates {
(3,0.0006591842261025673)
(4,0.0001460188310975341)
(5,2.687035538281663e-5)
(6,1.6575889773159994e-5)
(7,1.5194572000148738e-5)
(8,4.301183563968523e-7)
(9,1.4601740963796513e-7)

};

\end{semilogyaxis}
\end{tikzpicture} \subcaption{Underlying density $\rho_N$. }
\end{minipage}
\begin{minipage}[t]{.29\linewidth}
\centering
\begin{tikzpicture}[scale=1]
\begin{semilogyaxis}[scale only axis,width = 0.7\textwidth,xlabel = {level $n$},
grid=major,
]
\addplot[blue,mark=*,mark size=1.0pt]  coordinates { 		

(3,0.0010584775789297327)
(4,0.000102062122353739)
(5,1.0981497409963376e-5)
(6,1.2714389460559163e-6)
(7,1.60951953600631e-7)
(8,6.207937834941898e-8)

};

\addplot[red,mark=*,mark size=1.0pt]  coordinates {

(3,0.0012395196127276465)
(4,0.00014470995850951658)
(5,9.890275690367647e-6)
(6,1.22328385078784e-6)
(7,1.7060469936396376e-7)
(8,4.213893659715505e-8)

};

\addplot[green,mark=*,mark size=1.0pt]  coordinates {

(3,0.030051558223749632)
(4,0.004386705815527351)
(5,1.1553395212345e-5)
(6,7.620938046727559e-7)
(7,5.569246853031671e-7)
(8,4.2476904775263925e-8)

};
\addplot[black,mark=none]  coordinates {
(3,0.03675969665771165)
(4,0.0049225418148206804)
(5,1.3733643254085836e-5)
(6,9.49079409448594e-7)
(7,6.61347867849382e-7)
(8,6.657571819170098e-8)

};

\end{semilogyaxis}
\end{tikzpicture} \subcaption{Underlying density $\rho_C$. }
\end{minipage}
\begin{minipage}[t]{.29\linewidth}
\centering
\begin{tikzpicture}[scale=1]
\begin{semilogyaxis}[scale only axis,width = 0.7\textwidth,xlabel = {level $n$},
legend entries={\eqref{eq:S_nXR} ,$\sigma_{\text{ROT}}$,$\sigma_{\text{DPI}}$,\eqref{eq:bound_KDE2} },
grid=major,
legend style={at={(axis cs:12.1,10^2)},anchor=north west,scale =0.5},
legend cell align={left},
legend pos=outer north east,
]
\addplot[blue,mark=*,mark size=1.0pt]  coordinates { 		

(3,0.008704228877848658)
(4,0.0008090153343181091)
(5,7.700153673425294e-5)
(6,9.613731227831342e-6)
(7,1.1776836864402535e-6)
(8,1.4206076392318394e-7)

};

\addplot[red,mark=*,mark size=1.0pt]  coordinates {

(3,0.00915746124368886)
(4,0.0005107864176535586)
(5,8.995502042603405e-5)
(6,9.167402811214762e-6)
(7,9.785274108848377e-7)
(8,1.2887401746802769e-7)

};

\addplot[green,mark=*,mark size=1.0pt]  coordinates {

(3,0.008295362395300534)
(4,0.0005137371499151079)
(5,8.95564619294154e-5)
(6,9.285448082516797e-6)
(7,9.801459041230313e-7)
(8,1.3162562485705166e-7)

};

\addplot[black,mark=none]  coordinates {
(3,0.010179953815571297)
(4,0.0007553937855584819)
(5,0.00011145215793252741)
(6,1.0260366912448281e-5)
(7,1.097102529977316e-6)
(8,1.4096783373117679e-7)

};

\end{semilogyaxis}
\end{tikzpicture} \subcaption{Underlying density $\rho_L$. }
\end{minipage}
\caption{RMSE of the approximation on $\R$ of the test function~\eqref{eq:f_gauss} using kernel density estimation.}
\label{fig:num_gauss2}
\end{figure}

\textbf{The Gauss kernel on the real axis}\\
For the case where $\Omega =\R$, we choose $k(y) = \rho_N$ to be the standard normal distribution~\eqref{eq:rho_N}. 
The expressions used for estimating the smoothing parameter $\sigma_{\text{DPI}}$ are
$$\norm{k}_{L_2(\R)}^2= \frac{1}{2\,\sqrt{\pi}},\quad \mu_2(k) =1.$$
To study the performance of our algorithm, we use as test function again the function in~\eqref{eq:f_gauss}. 
We investigate the densities $\rho_N$, \eqref{eq:rho_N}, $\rho_C$, \eqref{eq:rho_C} and $\rho_L$, \eqref{eq:rho_L}.\\
Doing the same procedure as described in Section~\ref{sec:Numerics1} and we study the 
resulting RMSE~\eqref{eq:RMSE}. We used the two different proposed parameter selection methods and plotted the results in Figure~\ref{fig:num_gauss2}. 
It is reasonable to compare the results with the approximation error of the previous section, where we assume that the density $\rho$ is known, see~\eqref{eq:S_nXR} in Figure~\ref{fig:num_gauss2}. One can 
see that our approximation approach without knowing the density works well for all three examples, since for 
both investigated smoothing parameter selectors we end up with nearly the same error, which we get with knowing the density.
Furthermore, we had a look at the bound in Theorem~\ref{thm:bound_KDE}. We choose
$$\Eps = (-\infty,-\max |\Y_{\text{test}}|]\cap [\max |\Y_{\text{test}}|,\infty),$$
since this is the interval where we do not expect data. Then we calculated the right-hand side of~\eqref{eq:bound_KDE2} 
numerically for the choice $\sigma_{\text{DPI}}$. In Figure~\ref{fig:num_gauss2} we see that this is indeed a good upper bound for 
the approximation error from choosing $\sigma = \sigma_{\text{DPI}}$. \\

\textbf{The polynomial kernel on the cube}\\
Let the density be the beta distribution $\rho_{B,\alpha}$  from~\eqref{eq:rho_B} with shape parameter $\alpha\in \{\sfrac 12,1,2\}$. 
Let us study the test function $f(y) =\e^{y}$.
\begin{figure}[ht]
\begin{minipage}[t]{.29\linewidth}
\centering
\begin{tikzpicture}[scale=1]
\begin{semilogyaxis}[scale only axis, width = 0.7\textwidth,xlabel = {level $n$},
grid=major,
]
\addplot[blue,mark=*,mark size=1.0pt]  coordinates { 		

(3,0.002660999003417345)
(4,0.0009386321923443813)
(5,0.0006631824762506578)
(6,0.00034918183244285407)
(7,0.000196279908227528)
(8,6.444665186039913e-5)
(9,5.540737859184017e-5)

};

\addplot[red,mark=*,mark size=1.0pt]  coordinates {

(3,0.0946952780967285)
(4,0.048211420425529616)
(5,0.036775432385576866)
(6,0.028958495552849257)
(7,0.01045107479534404)
(8,0.0057691530527724435)
(9,0.0017727740961156688)

};

\addplot[black,mark=none]  coordinates {

(3,0.17951490333734782)
(4,0.12601756604277883)
(5,0.08871147373625804)
(6,0.06587549450487665)
(7,0.037340548657272245)
(8,0.02075966522800273)
(9,0.00782805143263955)

};

\end{semilogyaxis}
\end{tikzpicture}
\subcaption{Beta distribution with parameter $\alpha =2$. }
\end{minipage}
\begin{minipage}[t]{.29\linewidth}
\centering
\begin{tikzpicture}[scale=1]
\begin{semilogyaxis}[scale only axis,width = 0.7\textwidth,xlabel = {level $n$},
grid=major,
]
\addplot[blue,mark=*,mark size=1.0pt]  coordinates { 		
(3,3.949073316399214e-6)
(4,3.8753374706080655e-7)
(5,4.66377523690624e-8)
(6,5.4049933688129084e-9)
(7,6.399862849633528e-10)
(8,7.971056619158929e-11)
(9,9.894252122341505e-12)

};

\addplot[red,mark=*,mark size=1.0pt]  coordinates {
(3,0.04273237098898993)
(4,0.008924993495578843)
(5,2.5845092300748198e-5)
(6,6.156034187086481e-6)
(7,7.499656302944174e-7)
(8,1.5794891122407873e-7)
(9,2.7348571930724727e-8)

};

\addplot[black,mark=none]  coordinates {

(3,0.10723376259561074)
(4,0.023548432234421098)
(5,6.72141612759032e-5)
(6,1.4747631438311672e-5)
(7,1.835520724758931e-6)
(8,4.241695717859355e-7)
(9,6.715021560191758e-8)

};

\addplot[black,mark=none,dashed] coordinates {(3, 0.01* 2^-3*3 ) (9, 0.01 * 2^-3*9 )};

\end{semilogyaxis}
\end{tikzpicture} \subcaption{Uniformly distributed samples on $[0,1]$, i.e.\,$\alpha=1$. }
\end{minipage}
\begin{minipage}[t]{.29\linewidth}
\centering
\begin{tikzpicture}[scale=1]
\begin{semilogyaxis}[scale only axis,width = 0.7\textwidth,xlabel = {level $n$},
legend entries={\eqref{eq:S_nXR} ,$\sigma_{\text{ROT}}$,\eqref{eq:bound_KDE2}, $2^{-3n}$ },
grid=major,
legend cell align={left},
legend pos=outer north east,
]

\addplot[blue,mark=*,mark size=1.0pt]  coordinates { 		

(3,5.899921446550911e-5)
(4,6.149607854984911e-6)
(5,6.651732055677562e-7)
(6,7.984545740994539e-8)
(7,9.650569529217553e-9)
(8,1.1894379793330843e-9)
(9,1.4822157415535511e-10)

};

\addplot[red,mark=*,mark size=1.0pt]  coordinates {

(3,0.003445050423008265)
(4,3.2013660452462746e-5)
(5,3.961117025762617e-6)
(6,5.767933849945137e-7)
(7,8.812626801360644e-8)
(8,1.145050597440606e-8)
(9,1.471420709871818e-9)

};

\addplot[black,mark=none]  coordinates {

(3,0.006915447850318444)
(4,9.893003506865859e-5)
(5,9.641769309982119e-6)
(6,1.3039867776452978e-6)
(7,2.239009046564121e-7)
(8,2.5307471653514937e-8)
(9,3.7528437565297675e-9)

};

\addplot[black,mark=none,dashed] coordinates {(3, 0.01* 2^-3*3 ) (9, 0.01 * 2^-3*9 )};

\end{semilogyaxis}
\end{tikzpicture} \subcaption{Beta distribution with parameter $\alpha =\sfrac 12$. }
\end{minipage}
\hfill
\caption{RMSE of the approximation on the interval $[0,1]$ of the test function $f(y)=\e^{y}$ using kernel density estimation.}
\label{fig:num_cube2}
\end{figure}
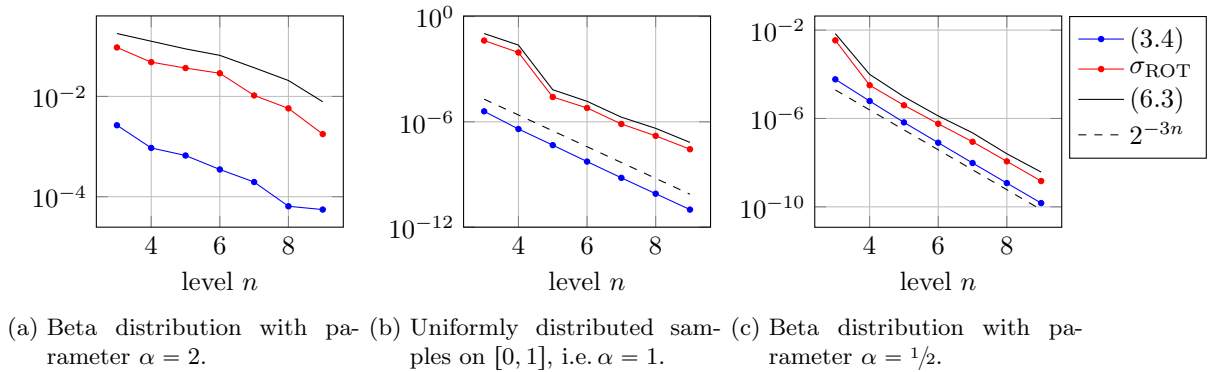 
We use a B-Spline kernel $k(y) = B_3(y)$, see~\eqref{eq:BSpline}. In this case the integral $K(y)$ can be calculated easily. 
The resulting RMSE are plotted in Figure~\ref{fig:num_cube2}. We compare with the case where the density in known,~\eqref{eq:S_nXR}. In the case where $\alpha=\sfrac 12$ the density 
$\rho_{B,\sfrac 12}$ is large on the boundary, see Figure~\ref{fig:example2}a), which means that we have more points at the boundary. 
In this case we receive the error rate $2^{-3n}$, which is specified by the order of the wavelets.
This behavior occurs even in the case $\alpha =1$, 
which means that the samples are uniformly distributed on $[0,1]$. 
In the case where $\alpha = 2$, the benefit of our approximation approach is not as big as in the other cases, since the density $\rho$ 
tends to zero at the boundary and we do not have much samples at the boundary. This means that the function $\r{f}$ does not have a support bigger 
than $[0,1]$ and the smoothing effect of $\r{f}\circ\r{\Rho}^{-1}$ at the boundary does not apply. The approximation is slightly better in the case where the density is known, but nevertheless this density does not inherit the decay rate from the other examples due to the lack of sample points at the boundary.

 \section{High-dimensional approximation }\label{sec:high-dim}
The main aim of this paper is the fast and effective approximation of high-dimensional functions. We study the setting where the variables $y_1,\ldots,y_d$ 
are \textit{independent}, which means that density $\rho(y)$ is a product of the one-dimensional densities~\eqref{eq:rho_d}.
Therefore, we transform every variable of the given samples separately after estimating one-dimensional densities. 
Additionally, we utilize the ANOVA decomposition from Section~\ref{sec:ANOVA_trafo} to deal with the curse of dimensionality. 

For the function $f\in L_2(\Omega,\rho)$ we have the ANOVA decomposition~\eqref{eq:ANOVA_terms_transform}.
The number of ANOVA terms of a function is equal to $2^d$ and therefore grows exponentially in the dimension $d$. This reflects the curse of dimensionality 
in a certain way and poses a problem for the approximation of a function. \changed{In high-dimensional settings, the underlying function can very often effectively represented as a sum of lower-order functions. In other words, the function can be expressed as a combination of component functions, where only $\nu \gg d$ variables out of the total $d$ variables are active in each component, \cite{DePeVo10, KuSlWaWo09}. Recent methods as ANOVAapprox~\cite{LiPoUl21,PoSc19a} (and the successful application to different datasets in \cite{PoSc21}), SALSA~\cite{SALSA}, SRFE~\cite{Ha21}, SHRIMP~\cite{Xie22}}
To this end we introduce the notion of effective dimension, see \cite{CaMoOw97}.  
\begin{Definition}For $0<\epsilon\leq 1$ the \textit{effective dimension} of $f$, in the \textit{superposition sense}, is the smallest integer $\nu\leq d$, such that 
$$\sum_{|\vec u|\leq \nu}\sigma^2_\rho(f_\vec u)\geq \epsilon \sigma_\rho^2 (f).$$ 
\end{Definition}
A function with low effective dimension allows a good approximation using only ANOVA-terms up to order $\nu$. To approximate $f_{\vec u}$, we have to use the transformation $\Rho_{\vec u}$, or if the density $\rho_{\vec u}$ is unknown, we have to estimate the one-dimensional densities $\rho_{i}(y_i)$ for $i\in \vec u$, and use a transformation 
$\r{\Rho}_{\vec u}(\vec y_{\vec u})\approx \Rho_{\vec u}(\vec y_{\vec u})$ to transform the 
samples $\Y_{\vec u}$ to $\r{\X_{\vec u}} = \r{\Rho}_{\vec u}(\Y_{\vec u})$.
Since we deal with independent input variables, we transform every variable separately, i.e.
\begin{equation*}\r{\Rho}_{\vec u}(\vec y) = \left(\r{\Rho}_i(y_i)\right)_{i\in \vec u},
\end{equation*}
where we get $\r{\Rho}_i$ from one-dimensional transformations~\eqref{eq:Rho_rho}.
In the truncated hyperbolic wavelet matrix $\vec A$ we insert only the indices $(\vec j,\vec k)$ belonging to the low-dimensional terms, i.e.
$$\vec A = \left(\psi_{\vec j,\vec k}^\per(\vec x_{\vec u})\right)_{\vec x_{\vec u}\in \r{\X}_{\vec u},(\vec j,\vec k)\in I_n^{\vec u}}.$$
For chosen $\nu \ll d$ we do this for all $\vec u$ with $|\vec u|\leq\nu$, i.e. 
$U_{\nu}: = \{\vec u\in [d]\mid |\vec u|\leq \nu\}$. This is summarized in Algorithm~\ref{alg:1}. 
The notation $\Y_{\{i\}}$ means analogously to the notation $\vec y_{\vec u}$ that we only consider the components $y_i$ of 
the samples in $\Y$. Similar to~\cite{LiPoUl21,PoSc19a} we calculate the variances of the approximations of the ANOVA-terms $f_{\vec u}$ and omit 
in a second approximation step the ones with low variance in order to increase the accuracy with a higher wavelet level $n$ for the important ones. Hence, in a second approximation step we use only the ANOVA-terms $\vec u\in U\subset \mathcal P([d]) $ and get the approximant $S_n^{\Y,U} f$, \eqref{eq:SNY}, also for an arbitrary ANOVA index set $U$. \\

\begin{figure}[ht]
\begin{center}
\begin{tikzcd}[row sep=scriptsize, column sep=scriptsize]
\text{known } \rho_i \arrow[rr]&&\text{transformation } \Rho_i \text{ from \eqref{eq:Rho}}\arrow[rr]&& \X_i = \Rho(\Y_i)\\
\Y_i \arrow[d] \arrow[u] \arrow[rrrru] \arrow[rrrrd]&&&&\\
 \text{unknown } \rho_i   \arrow[rr]&& \text{transformation } \r{\Rho}_i \text{ from \eqref{eq:Rho_rho}}\arrow[rr]&&\r{\X}_i = \r{\Rho}(\Y_i)
\end{tikzcd}
\end{center}
\caption{For every variable $y_i$ of the sample set $\Y$ we have these possibilities.} 
\label{fig:skizze2}
\end{figure}
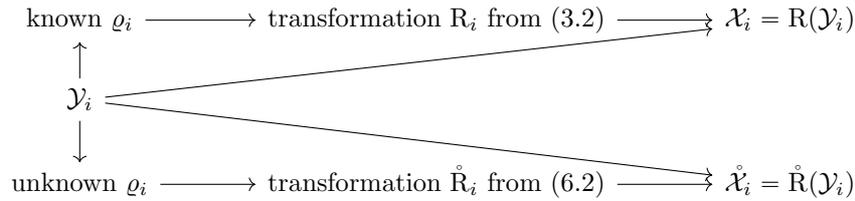

\begin{algorithm}[ht]
\caption{Transformed ANOVA hyperbolic wavelet regression}
	\vspace{2mm}
	\begin{tabular}{ l l l }
		\textbf{Input:} & $ d $ & dimension\\
		& $ \nu $ & superposition dimension\\
		&	$\Y = (y_i)_{i=1}^M\in \Omega$ & sampling nodes \\
		& $\vec f = (f(y_i))_{i=1}^M$ & function values at sampling nodes 
	\end{tabular}
	\begin{algorithmic}[1]
			\STATE{Choose $n$ such that for $N=\sum_{|\vec u|\leq\nu}|I_n^{\vec u}|$ holds $M>N\log N$.}
			\FOR{$i = 1,\ldots,d$ }
				\IF{$\rho_{i}(y_i)$ is known}
				\STATE{Calculate the transformation ${\Rho}_{i}$ by~\eqref{eq:Rho}.}
				\STATE{Transform the samples $\Y_{i}$ to $\X_{i}={\Rho}_{i}(\Y_{i})$.}
				\STATE{Do the steps 13-15 with $\Rho$ instead of $\r{\Rho}$.}
				\ELSE
					\STATE{Estimate $\r{\rho_{i}}$ with~\eqref{eq:KDE}.}
					\STATE{Calculate the transformation $\r{\Rho}_{i}$ by \eqref{eq:Rho_rho}.}
					\STATE{Transform the samples $\Y_{i}$ to $\r{\X}_{i}=\r{\Rho}_{i}(\Y_{i}) $.}
				\ENDIF
      \ENDFOR
			\STATE{Construct the sparse matrix 
			$$\vec A=[\vec A_{\vec u}]_{|\vec u |\leq \nu}\in \C^{M\times N}, \quad \vec A_{\vec u} = (\psi_{\vec j,\vec k}^\per(\vec x_{\vec u}))_{\vec x_{\vec u}\in \r{\X_{\vec u}},(\vec j,\vec k)\in I_n^{\vec u}}.$$}
	    \STATE{Solve the overdetermined linear system
			$\vec A \left(a_{\vec j,\vec k}\right)_{\vec j,\vec k} =\vec f$ 
			via an LSQR-algorithm. This gives us the approximation
			$$S_n^{\r{\X},\nu} (f\circ \r{\Rho})(\vec x) :=\sum_{|\vec u|\leq \nu}\sum_{(\vec j,\vec k)\in I_n^{\vec u}}a_{\vec j,\vec k}\psi_{\vec j,\vec k}^\per(\vec x)$$
			}
			\STATE{Transform the approximation back to $\Omega$ using $\r{\Rho}_{\vec u}^{-1}$ for $|\vec u|\leq \nu$. }
	\end{algorithmic}
	\begin{tabular}{ l l l }
		\textbf{Output:} 
		 &  $\left(a_{\vec j,\vec k}\right)_{\vec j,\vec k}\in \C^N$ coefficients of the approximant 
	\end{tabular}
	\begin{equation}\label{eq:SNY}
		S_n^{\Y,U} f (\vec y):=\sum_{\vec u \in U }\sum_{(\vec j,\vec k)\in I_n^{\vec u}}a_{\vec j,\vec k}\psi_{\vec j,\vec k}^\per(\r{\Rho}_{\vec u}^{-1}(\vec y_{\vec u})),
		\end{equation}
		where $U = U_\nu$.
	\label{alg:1}
\end{algorithm}  
To summarize our algorithm, we have for every variable $y_i$ two possibilities, which is summarized in Figure~\ref{fig:skizze2}.
Algorithm~\ref{alg:1} works well if the underlying density is a tensor density. 

\subsection{Approximating the global sensitivity indices}A direct calculation of sensitivity indices $S(\vec u,f)$ from~\eqref{eq:gsi_rho} would require integral evaluations in~\eqref{eq:ANOVA_rho}, followed by numerous integral
evaluations of sensitivity indices~\eqref{eq:sigma_rho}. For high-dimensional systems, such an approach is impractical and possibly prohibitive.
Therefore, alternative routes must be charted to estimate the sensitivity indices both accurately and efficiently. Our approach is 
to approximate the function $f$ by $S_n^\Y f$ and afterwards calculating the \changed{GSI}s of the approximation. It was also shown in~\cite{GrKuSl10} that ANOVA terms inherit the smoothness of a function, i.e.\,if $f\in H^s_\mix(\T^d)$, then $f_{\vec u}\in H^s_{\mix}(\T^{|\vec u|})$ or even smoother. This fact was also shown in~\cite[Theorem 3.10]{PoSc19a} by using a Fourier based approach. Following this lines we immediately have this result for Besov spaces, i.e.\,if $f\in \bB_{2,\infty}^s(\T^d)$ then $f_{\vec u}\in \bB^s_{\mix}(\T^{|\vec u|})$. \\

The intuition is that a good approximation 
of the function means also a good approximation of the ANOVA-terms, and hence a good approximation of the variances. Calculating the variances of the ANOVA-terms for functions on $\T^d$ is easy because of the connection~\eqref{eq:I_n_decomp}. Therefore we approximate the variances $\sigma^2_{\rho}(f_{\vec u})$ by the following estimated variances,
\begin{align*}
\tilde{\sigma}_{\rho}^2(f_{\vec u}) &:= \int_{\T^{|\vec u|}} |(S_n^{\X}(f\circ \Rho_{\vec u}^{-1}))_{\vec u}|^2 \d \vec x_{\vec u} = \int_{\Omega_{\vec u}} |(S_n^{\X}(f\circ \Rho_{\vec u}^{-1}))_{\vec u}\circ\Rho_{\vec u}|^2 \,\frac{\dx }{\dx \vec y}\Rho_{\vec u}(\vec y_{\vec u})\dx \vec y_{\vec u},\\
\tilde{\sigma}_{\r{\rho}}^2(f_{\vec u}) &:= \int_{\T^{|\vec u|}} |(S_n^{\r{\X}}(f\circ \r{\Rho}^{-1}))_{\vec u}|^2 \d \vec x_{\vec u} = \int_{\Omega_{\vec u}} |(S_n^{\r{\X}}(f\circ \r{\Rho}_{\vec u}^{-1}))_{\vec u}\circ\r{\Rho}_{\vec u}|^2 \,\frac{\dx }{\dx \vec y}\r{\Rho}_{\vec u}(\vec y_{\vec u})\dx \vec y_{\vec u},
\end{align*}
see \cite[Section 4]{LiPoUl21} for computing details with hyperbolic wavelet regression. In the following we will study the error between the estimated variances $\tilde{\sigma}^2_{\rho}(f_{\vec u})$ and the variances $\sigma^2_{\rho}(f_{\vec u})$.
First, let us consider the case where the density is known.
\begin{theorem}\label{thm:gsi1}
Let $\varnothing\neq \vec u\in U$ and $\vec v = \{i\in \vec u\mid \Omega_i = [0,1]\}$ and $f_{\vec u}\in \bB_{2,\infty}^s(\Omega_{\vec u},\rho_{\vec u}) $. Denote furthermore $e_2 := \norm{f_{\vec u}- (S_n^{\X}(f\circ \Rho_{\vec u}^{-1}))\circ \Rho}$ and $a = 1-(1-\eta)^{\sfrac{|\vec v|}{2}-1}$. Then
\begin{equation*}
|\sigma^2_\rho(f_{\vec u})- \tilde{\sigma}^2_\rho(f_{\vec u})| \leq \left(e_2+a\norm{f_{\vec u}}_{L_2(\Omega_{\vec u},\rho_{\vec u})}\right) \,\left( 2+ e_2+a\norm{f_{\vec u}}_{L_2(\Omega_{\vec u},\rho_{\vec u})}\right)\, \norm{f_{\vec u}}_{L_2(\Omega_{\vec u},\rho_{\vec u})}.
\end{equation*}
\end{theorem}
\begin{proof}
Let us denote in this proof $g=(S_n^{\X}(f\circ \Rho_{\vec u}^{-1}))\circ \Rho$. Then we have
$$\tilde{\sigma}_{\rho}^2(f_{\vec u}) = \int_{\Omega_{\vec u}} |g(\vec y_{\vec u})|^2\frac{\dx }{\dx \vec y}\Rho_{\vec u}(\vec y_{\vec u}) \dx \vec y_{\vec u} = \int_{\Omega_{\vec u}} |g(\vec y_{\vec u})|^2\left(1-\eta\right)^{-|\vec v|}\rho_{\vec u}(\vec y_{\vec u}) \dx \vec y_{\vec u}.$$
For simplicity we denote in this proof $\tilde{g} = (1-\eta)^{-\sfrac{|\vec v|}{2}} \, g$.
Then we estimate the difference of the variances of the ANOVA terms given in~\eqref{eq:ANOVA_rho} by the reverse triangle inequality and Cauchy-Schwarz inequality,
\begin{align*}|\sigma^2_\rho(f_{\vec u})-\tilde{\sigma}^2_\rho(f_{\vec u})| &=\left|\int_{\Omega_{\vec u} } \left(|f_{\vec u}(\vec y_{\vec u})|^2-|\tilde{g}(\vec y_{\vec u} )|^2\right) \,\rho_{\vec u}(\vec y_{\vec u})\d \vec y_{\vec u} \right| \notag\\
&=\left|\int_{\Omega_{\vec u} } \left(|f_{\vec u}(\vec y_{\vec u})|-|\tilde{g}(\vec y_{\vec u} )|\right)\,\left(|f_{\vec u}(\vec y_{\vec u})|+|\tilde{g}(\vec y_{\vec u} )|\right)\,\rho_{\vec u}(\vec y_{\vec u})\d \vec y_{\vec u} \right|\notag\\
&\leq \norm{f_{\vec u}-\tilde{g}}_{L_2(\Omega_{\vec u},\rho_{\vec u})} \, \norm{f_{\vec u}+\tilde{g}}_{L_2(\Omega_{\vec u},\rho_{\vec u})}\notag\\
&\leq \norm{f_{\vec u}-\tilde{g}}_{L_2(\Omega_{\vec u},\rho_{\vec u})} \, \norm{f_{\vec u}}_{L_2(\Omega_{\vec u},\rho_{\vec u})}\left(2+\norm{f_{\vec u}-\tilde{g}}_{L_2(\Omega_{\vec u},\rho_{\vec u})}\right).
\end{align*}
Only in the case where $\vec v\neq \varnothing$ we have the additional factor, which depends on $\eta$:
$$
\norm{f_{\vec u}-\tilde{g}}_{L_2(\Omega_{\vec u},\rho_{\vec u})}\leq \frac{1}{(1-\eta)^{\sfrac{|\vec v|}{2}}} \norm{f_{\vec u}-g}_{L_2(\Omega_{\vec u},\rho_{\vec u})} + \left(1-(1-\eta)^{\sfrac{|\vec v|}{2}-1}\right)\norm{f_{\vec u}}_{L_2(\Omega_{\vec u},\rho_{\vec u})} 
$$
In the case of periodic functions, the second term is zero.
Putting all inequalities together gives the desired result. 
\end{proof}
The error between the function $f_{\vec u}$ and its approximation $g$ can be estimated as follows. We are now concerned with a $|\vec u|$-dimensional function. Therefore estimates~\eqref{eq:P_n_2} to \eqref{eq:S_n_3} hold for $d=|\vec u|$ for the transformed function. The connection~\eqref{eq:L2_equality} or rather Theorem~\ref{thm:decay_trafo} transforms the results to $\Omega_{\vec u}$. Therefore, we have
\begin{equation*}
\norm{f_{\vec u}-g}_{L_2(\Omega_{\vec u},\rho_{\vec u})}\lesssim  2^{-2ns }n^{|\vec u|-1}\norm{f_{\vec u}}_{\bB_{2,\infty}^s},
\end{equation*}
with high probability.
The logarithmic term in the approximation error appears only for ANOVA terms with $|\vec u|\geq 2$. In the case where the 
density $\rho$ is unknown we get an additional term, which depends on the error 
between the estimated density $\r{\rho}$ and the actual density $\rho$, similar as in Theorem~\ref{thm:bound_KDE}.\\
We introduce the subset 
$$\Eps := \bigtimes\limits_{i=1}^d \left(\Omega_i\backslash [\min \X_{\{i\}},\max \X_{\{i\}} ]\right),$$
for which the estimate in the next theorem follows.
\begin{theorem}
Let $\varnothing\neq \vec u\in U$ and $f_{\vec u}\in \bB_{2,\infty}^s(\Omega_{\vec u},\rho_{\vec u}) $. Let furthermore $g :=S_n^{\r{\X},U}f$ be an approximation received from our procedure in Algorithm~\ref{alg:1} with unknown density. Then
\begin{align*}
|\sigma^2_\rho(f_{\vec u})- \tilde{\sigma}^2_\rho(f_{\vec u})|&\leq \underbrace{|\sigma^2_\rho(f_{\vec u})- \sigma^2_\rho(g_{\vec u})|}_{A} + \underbrace{\left|\int_{\Eps} |g_{\vec u}(\vec y)|^2(\rho(\vec y)-\r{\rho}(\vec y))\d \vec y\right|}_{B}\\
& \quad+ \underbrace{\left|\int_{\Omega_{\vec u} \backslash \Eps} |g_{\vec u}(\vec y)|^2 \rho(\vec y)\d \vec y\right| \, \sup_{\vec y\in \Omega_{\vec u} \backslash \Eps} \left(1-\frac{\r{\rho}(\vec y)}{\rho(\vec y)}\right)}_{C}.
\end{align*}
\end{theorem}
Let us briefly explain the appearing terms. Term $A$ is the error from Theorem~\ref{thm:gsi1}, terms $C$ depends on the quality of the approximation $\r{\rho}\approx \rho$ and term $B$ describes the variance of $g_{\vec u}$ in the part where we have no samples, i.e.\,where we extend the original function $f_{\vec u}$. Of course, if the original function is non-periodic, we use an extension and study the variances, so we have to accept the additional term. 
\begin{proof}
We do the following splitting
\begin{equation*}
|\sigma^2_\rho(f_{\vec u})- \sigma^2_{\r{\rho}}(g_{\vec u})|\leq |\sigma^2_\rho(f_{\vec u})- \sigma^2_\rho(g_{\vec u})| +  |\sigma^2_\rho(g_{\vec u})- \sigma^2_{\r{\rho}}(g_{\vec u})|.
\end{equation*}
Then by splitting the domain of the second integral the assertion follows.

\end{proof}
The quintessence of this subsection is that the approximation of the \changed{GSI} of the function $f$ by the \changed{GSI} 
of the approximant is a reasonable approach to get insides about the variances of the ANOVA terms.
In a second approximation step we reduce the index set to the ANOVA indices in the set 
$$ U = \{\vec u\in U_{\nu}\mid S(\vec u,S_n^\Y f)>\epsilon \},$$
for some threshold $\epsilon >0$. This allows on the other hand to increase the maximal wavelet level for the 
important ANOVA terms and therefore decrease the approximation error, while \changed{ensuring} logarithmic oversampling.

\subsection{\changed{A synthetic} numerical example}\label{sec:num_high}
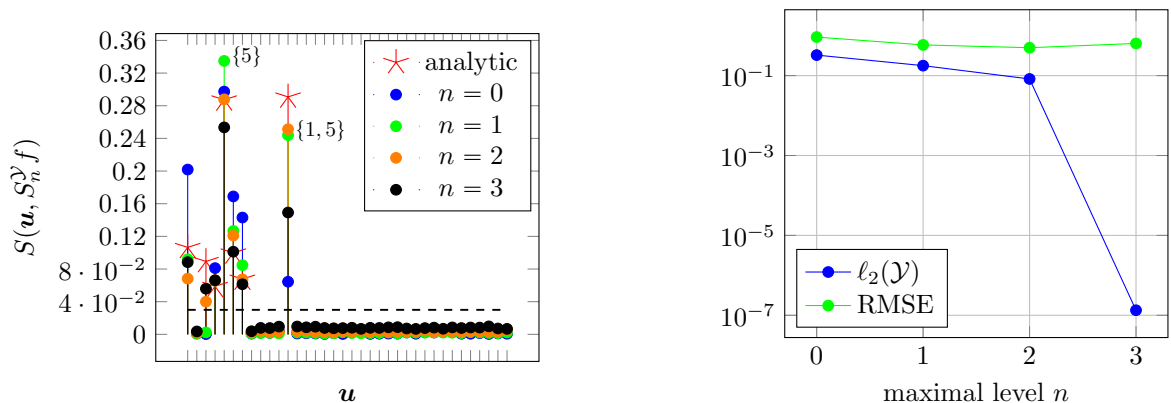
\begin{figure}[htb]
\begin{minipage}[t]{0.45\textwidth}
\centering \begin{tikzpicture}[scale=1]
\begin{axis}[scale only axis,width = 0.7\textwidth,
ylabel={$S(\vec u,S_n^\Y f)$},
			xlabel={$\vec u$},
			xtick={1,...,36},
			xticklabel=\empty,
			ytick distance=0.04,
			xtick distance=10,
			legend entries={analytic,$n=0$,$n=1$,$n=2$,$n=3$},
legend style={scale =0.7},
			]
			
			\addplot+[ycomb, red,mark=star, mark size =5pt,mark options={red}] plot coordinates
			{

(1,0.1065)
(3,0.0893)
(4,0.0588)
(5,0.2872)
(6,0.0998)
(7,0.0677)
(12,0.2908)

};
			\addplot+[ycomb, blue,mark=*, mark options={blue}] plot coordinates  {
(1,0.20185261681660402)
(2,0.0010621892400989232)
(3,0.00030840051724795076)
(4,0.0811317714717183)
(5,0.2973783662567179)
(6,0.16891609051387502)
(7,0.14308036600667012)
(8,0.0009241182720922288)
(9,0.001892050790021587)
(10,0.0018042297341934529)
(11,0.0038664577935982464)
(12,0.06443863122811791)
(13,0.0010031291972906566)
(14,0.001274277412006862)
(15,0.0009911065291333717)
(16,0.000445701842342459)
(17,0.0014141635198800368)
(18,0.00041568564747853997)
(19,0.002643135317399652)
(20,0.003178613540262975)
(21,0.000539294749203888)
(22,0.000771923425932143)
(23,0.0009024937652668005)
(24,0.0008713048050472376)
(25,0.0022973683713043225)
(26,0.0007300418031358001)
(27,0.002541314997896301)
(28,0.0022166648117421016)
(29,0.0032063405619181195)
(30,0.0027931955312338877)
(31,0.0007234633620241199)
(32,0.0007223584086842382)
(33,0.0014397268504037566)
(34,0.0004324027436983703)
(35,0.001057690180981764)
(36,0.0007333139847770298)

};

			\addplot+[ycomb, green,mark=*, mark options={green}] plot coordinates  {
(1,0.09192763891150246)
(2,0.0006372905365779366)
(3,0.00234171777665402)
(4,0.06584940853309933)
(5,0.3349056758956196)
(6,0.1266269361766236)
(7,0.08473708294396545)
(8,0.0009161269867354382)
(9,0.0015970156373992055)
(10,0.0017224926908846242)
(11,0.0010583211332251662)
(12,0.24396276090162142)
(13,0.0021459387198220145)
(14,0.002018106107461745)
(15,0.0015479368442728255)
(16,0.0009842881904205908)
(17,0.001158511539140351)
(18,0.002527312615090182)
(19,0.0015159423587749305)
(20,0.0012690919178997312)
(21,0.0015065262877247405)
(22,0.00160088456712782)
(23,0.002239618611393597)
(24,0.0010190374348954673)
(25,0.0020936768587674756)
(26,0.0019899745986401616)
(27,0.002079054513522335)
(28,0.0023521786720228028)
(29,0.002107400106628336)
(30,0.0019417973954664574)
(31,0.001521340351869637)
(32,0.0026050949555128907)
(33,0.0021461035233148696)
(34,0.0016293095806970668)
(35,0.0014796573895138364)
(36,0.0022387487361118725)

};
			\addplot+[ycomb, orange,mark=*, mark options={orange}] plot coordinates  {
(1,0.0682531989309356)
(2,0.0019071555108854694)
(3,0.040012016758232447)
(4,0.06617033857600366)
(5,0.2878724924726393)
(6,0.12098616544734031)
(7,0.06730499649143677)
(8,0.002118199036984878)
(9,0.004150364296299761)
(10,0.0028779286392966374)
(11,0.003256839229050903)
(12,0.25127140551454585)
(13,0.0040239252321905275)
(14,0.00431217291560764)
(15,0.003251417418761797)
(16,0.0027544467755198664)
(17,0.0022859535573446126)
(18,0.003815270331649489)
(19,0.0030644448926863146)
(20,0.004186621157549701)
(21,0.002598700552557204)
(22,0.004306221925817571)
(23,0.002494134816889332)
(24,0.003726042658498863)
(25,0.0034038161558703054)
(26,0.0033756082610097633)
(27,0.004553705373523317)
(28,0.002676882060685362)
(29,0.0035823047744568772)
(30,0.0037240823669593666)
(31,0.0030037136421155646)
(32,0.0035897118610402263)
(33,0.0035011486779258235)
(34,0.0032689140262553915)
(35,0.002897674283247811)
(36,0.00542198537818577)

};
\addplot+[ycomb, black,mark=*, mark options={black}] plot coordinates  {
(1,0.08838984721145937)
(2,0.0035353635959733997)
(3,0.055885498860113535)
(4,0.06625104555878056)
(5,0.25349334989666905)
(6,0.10120781305601394)
(7,0.06135132699346565)
(8,0.0038318332720535483)
(9,0.007887467334716741)
(10,0.007746372289709529)
(11,0.009485703649498353)
(12,0.14920196010848855)
(13,0.009602122688141663)
(14,0.00864136921982371)
(15,0.009373405420685538)
(16,0.007842903438571136)
(17,0.0074344057934787435)
(18,0.007679256492375027)
(19,0.00808124054488589)
(20,0.006714394705498264)
(21,0.007983199853292083)
(22,0.007842706364288847)
(23,0.008610810146905741)
(24,0.008814127751488463)
(25,0.007152544775779537)
(26,0.00666758535353251)
(27,0.007667767818454333)
(28,0.008246554754230449)
(29,0.00684733831347162)
(30,0.008505567615845928)
(31,0.007716473559365581)
(32,0.008372153714740056)
(33,0.008288396090516013)
(34,0.009770991898635527)
(35,0.0071792801421613995)
(36,0.006697821716889821)

};
\addplot[ black,mark=none, mark options={black},dashed,line width=0.25mm] plot coordinates  {(1,0.03) (36,0.03)};

\node at (axis cs:1,0.4) [anchor=west,scale =0.8] {$\{1\}$};
\node at (axis cs:5,0.34) [anchor=west,scale =0.8] {$\{5\}$};
\node at (axis cs:12,0.25) [anchor=west,scale =0.8] {$\{1,5\}$};
		\end{axis}
\end{tikzpicture}\subcaption{The approximated \changed{GSI}s $S(\vec u,S_n^\Y f)$ for ${|\vec u|\leq 2}$ and different level $n$ compared to the \changed{GSI}s $S(\vec u, f)$.  }
\label{fig:gsis}
\end{minipage}
\hfill
\begin{minipage}[t]{0.45\textwidth}
\centering \begin{tikzpicture}[scale=1]
\begin{semilogyaxis}[scale only axis,width = 0.7\textwidth,xlabel = {maximal level $n$},
legend entries={$\ell_2(\Y)$,RMSE},
grid=major,
legend style={scale =0.7},
legend cell align={left},
legend pos= south west
]

\addplot[blue,mark=*,]  coordinates {

(0,0.3265599971900143)
(1,0.17788958734087637)
(2,0.08255619726394425)
(3,1.3361428334825033e-7)

};

\addplot[green,mark=*]  coordinates {
(0,0.923873423198099)
(1,0.5874770701336299)
(2,0.4997157754066017)
(3,0.6394584383984888)

};
\end{semilogyaxis}
\end{tikzpicture} 
\subcaption{Approximation errors, the $\ell_2$-norm on the given samples and the RMSE for different wavelet levels $n$.}
\label{fig:errors}
\end{minipage}
\caption{Approximation of the function $f$ from~\eqref{eq:f_8} from $M=1000$ samples. }
\end{figure}

As a conclusion of this paper we want to apply Algorithm~\ref{alg:1} to a high-dimensional test function. For that reason let
\begin{equation}\label{eq:f_8}
f\colon \R^5\times [0,1]^3\rightarrow \R,\quad f(\vec y) =\tfrac 15 y_1^2+\tfrac 12\cos(2\pi y_3) +\e^{-y_4^2} + y_5^{\sfrac 12} +30\,(y_6^3\,(1-y_6^2)) +\tfrac 12|4\,y_7-2|+5\e^{-y_1^2-y_5^2}
\end{equation} 
be an $8$-dimensional function where $y_5>0$. We assume the given data $\Y$ to be sampled from the distribution
\begin{equation*}\rho\colon \R^5\times [0,1]^3\rightarrow \R_+,\quad \rho(\vec y) = \prod_{i=1}^8 \rho_i(y_i),
\end{equation*} 
where we use different distributions, already studied in this paper, 
\begin{align*}
\rho_1(y_1) &= \rho_N(y_1), \quad \rho_2(y_2) = \rho_L(y_2), \quad \rho_3(y_3) = \rho_C(y_3), \quad \rho_5(y_5) = \frac 12 \e^{-\tfrac {y_5}{2}}, \quad \rho_6(y_6)=\rho_8(y_8) =1,\\
\rho_4(y_4) &= \frac{1}{\sqrt{11.52\,\pi}}\e^{-\frac{(y_4+2)^2}{2.88}}+\frac{1}{\sqrt{50\,\pi}}\e^{-\frac{(y_4-3)^2}{12.5}}, \quad
\rho_7(y_7) =\changed{ \frac{1}{\pi}} \,{y_7}^{-\sfrac 12} (1-y_7)^{-\sfrac 12}.
\end{align*} 
We draw $M=1000$ samples and use the corresponding function values $\vec f = f(\Y)$. 
These samples projected to the directions $y_1$ and $y_2$ are plotted in the introduction in Figure~\ref{fig:samples} together with the transformed samples $\r{\Rho}(\Y)$, also projected to the directions $y_1$ and $y_2$. We use as superposition dimension $\nu=2$, 
which is a suitable guess if we have a look at the function equation, which suggests only one ANOVA-term of order $2$. With the 
choice $\nu=3$ we would conclude in a first step that we do not need the three-dimensional terms. Furthermore we use Chui Wang wavelets of order $m=2$.\\ 
We consider the setting where we do not know the underlying density, so we use for the variables $y_1,\ldots,y_5$ the kernel 
density estimation for data on $\R$ from Section~\ref{sec:KDE_R} using the Gaussian kernel introduced there and for the remaining 
variables $y_6,y_7,y_8$ the kernel density estimation for data on $[0,1]$ from Section~\ref{sec:KDE_cube} using the B-Spline kernel introduced there. 
For different wavelet levels $n$ we plot in Figure~\ref{fig:gsis} the approximated \changed{GSI}'s $S(\vec u,S_n^\Y f)$, i.e.\,the $8$ \changed{GSI}'s of order $1$ and then the $28$ \changed{GSI}'s of order $2$. Since we know the function explicitly, we compare this to the analytically calculated \changed{GSI}'s $S(\vec u, f)$. One can see that we could indeed figure out even with a low wavelet level $n=2$ the ANOVA terms with high variances. So we filter out the unnecessary variables $y_2$ and $y_8$ and all two-dimensional terms except the term with $\vec u=\{1,5\}$. 
Furthermore, we plotted in Figure~\ref{fig:errors} the error $\norm{f-S_n^\Y f}_{\ell_2(\Y)}=\left(\sum_{\vec y\in \Y} |f(\vec y)-S_n^\Y f(\vec y)|^2\right)^{\sfrac 12}$ and the RMSE~\eqref{eq:RMSE} for a test set $\Y_{\text{test}}$ sampled according to $\rho$ with $|\Y_{\text{test}}|= 3M$. The low $\ell_2(\Y)$- error indicates that $n=3$ is already overfitting, i.e.\,using to many parameters for the $1000$ samples. \\
In a second approximation step we omit the ANOVA-terms with low variance and use only the ANOVA indices $ U = \{\varnothing, \{1\},\{3\},\{4\},\{5\},\{6\},\{7\},\{1,5\}\}$ for the approximation. This procedure is similar to the proposed method suggested in~\cite{LiPoUl21, PoSc21}. It is also reasonable to choose different maximal levels $n$ for the one- and the two-dimensional terms, since for different dimension, these index-sets have different sizes. For the choice $n=5$ for the one-dimensional terms and the choice $n =3 $ for the two-dimensional term we are able to reduce the RMSE on the test-set additionally from $0.4997$ to $0.2248$. Our procedure reduces the RMSE significantly, hence we are able to approximate an $8$-dimensional function \changed{using} only $1000$ samples very well. \changed{Using only a min-max transformation of the data, it is not possible to detect the non-zero ANOVA-terms.}\\

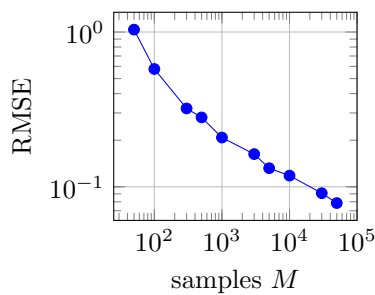
\begin{figure}[htb]
\centering \begin{tikzpicture}[scale=1]
\begin{loglogaxis}[scale only axis,width = 0.2\textwidth,xlabel = {samples $M$},
ylabel = {RMSE},
grid=major,
legend cell align={left},
legend pos= south west
]

\addplot[blue,mark=*,]  coordinates {

(50,1.036650276346503)
(100,0.5769362282352161)
(300,0.32077989614930696)
(500,0.2807097114222667)
(1000,0.20800088400117833)
(3000, 0.1626812590886379)
(5000,0.13186221499927248)
(10000, 0.11809370793506997)
(30000,0.09089382540635382)
(50000,0.0785240746878328)

};

\end{loglogaxis}
\end{tikzpicture} 
\caption{Approximation of the function $f$ from~\eqref{eq:f_8} for different numbers of samples $M$. }
\label{fig:num8_2}
\end{figure}

In a second experiment we only use the ANOVA indices in $ U$ for the approximation and do our procedure for different sample sizes $M$, \changed{with all other parameters kept the same}. The results are plotted in Figure~\ref{fig:num8_2}. We used the maximal wavelet level $n$ for which the RMSE is minimal for the given fixed data set.

\subsection{Real-world data}
\changed{
The proposed Algorithm~\ref{alg:1} was used in~\cite{WeKr23} to estimate the vertical ground reaction force from time series
of plantar pressure from instrumented insoles. The study
included data from $16$ persons moving at different speeds on a two-belt treadmill equipped with force
plates. In total about $M\approx 1.2 \cdot 10^6$ data points were used and the data was modelled as an $8$-dimensional function with ANOVA terms up to order $\nu =2$. The approach successfully reached relative RMSEs of up to $10.6\,\%$, which is comparable to other studies in the literature with the advantage of being interpretable and having much more data available in the study used in~\cite{WeKr23}.} \\

\changed{In the following we compare the performance of models obtained by Algorithm~\ref{alg:1} with
other state-of-the-art algorithms when applied to seven real-world
datasets from the the UCI repository (\url{http://archive.ics.uci.edu/ml}). For comparison we test against Random Forest Regression (RF) and Gaussian Processes (GP), both implemented in ScikitLearn.jl, which implements the popular scikit-learn interface and algorithms in Julia. Furthermore, we compare with sparse additive models \cite{SALSA, HARFE, Xie22} and follow the the experimental setup in~\cite{SALSA}, where the training data is normalized so that the
input and output values have zero mean and unit variance along each dimensions. Each dataset is divided in half to form the training and testing sets, we use exactly the same splitting as in~\cite{SALSA} for the datasets used there and do the same procedure for the datasets with more samples. Note that they test only one single random splitting. In our experiments we use a cross-validation of the original training data-set as training data and $20\%$
as validation data to select the best parameters $\nu, n, \lambda, \epsilon$. An overview of the datasets is presented in Table~\ref{tab:setup}.
Furthermore, we give details of our trained model: The used basis functions on the torus (chui-m are Chui-Wang-wavelets of order $m$ and per means trigonometrical polynomials) as well as the total number $N$ of trained coefficients in the final model. \\
For the transformation $\Rho$ we use the DPI, described in Section~\ref{sec:KDE_R} applied to the Gaussian kernel. }


\begin{table}[hbt]\centering
\begin{tabular}{|l|ccc|cc|}
    \hline
    Dataset &$d$   & $M_{\text{train}} $& $M_{\text{test}}$ & basis & $N$\\
    \hline
		Propulsion			&$15$ & $200$ 		& $200$ &chui-2 & $271$\\	
		Galaxy 					&$20$ & $2000$ 		& $2000$&per &$983$\\
		Airfoil					&$41$ & $750$ 		& $750$ &chui-2 &$339$\\	
		Forestfire  		&$10$ & $211$ 		& $167$&per  &$100$	\\	
		Boston Housing 	&$12$ & $256$ 		& $250$&chui-2  &$166$	\\	
		Protein 				&$9$ & $22685$ 		& $22685$&per &	$	2082$\\
		Elevators				&$17$&	$8300$		&$8399$&chui-4&$2332$\\
				\hline

  \end{tabular}
\caption{Overview of seven datasets: Dimension $d$ and number of datapoints in training and testing set.
The experimental setup and datasets for each test follows from~\cite{SALSA}.
We give also the used basis functions on the torus, the total number $N$ of trained coefficients in the final model of Algorithm~\ref{alg:1}. }
\label{tab:setup}
\end{table}

\changed{The approximation results and comparisons are shown in Table~\ref{tab:results}. The results of SALSA are obtained from~\cite{SALSA}, the results from 
HARFE and SRFE are obtained from~\cite{HARFE},
and we included the results of the SHRIMP algorithm~\cite{Xie22}.
Since the results for the random algorithms depend on the draw of the random features, 
in contrast to the given results for SHRIMP, HARFE and SALSA, we did the approximation validation $50$ times and present the mean in 
Table~\ref{tab:results}.  Note that our parameter $\nu$ coincides with the parameter $q$ in the random feature literature. 
Furthermore, for the datasets with too many samples, i.e. Protein and Elevators, the random feature algorithms are not able to calculate an approximation, since the involved random matrices are getting too big.}

\begin{table}[hbt]\centering
		%
%
	%
	
	\begin{tabular}{|l|lllllll|}
    \hline
													& Propulsion 				&  Galaxy   			 & Airfoil 		&  Forestfire		& Housing  & Protein & Elevators\\
    \hline
Alg.~\ref{alg:1} ($\nu$)	& $\textbf{0.0001126} (2)$		& $0.00344 (1)$			&$\textbf{0.1530}(2)$	&$0.3460(1)$		&$0.3779(1)$ &$\textbf{0.4095}(3)$ &$\textbf{0.2488}(1)$\\
\hline
RF												& $0.005928 $				& $0.1092$							&$0.6358$ 		&$0.3372$ 			&$0.3339$	&$0.4225$&0.2908\\
GP												& $0.009031 $				& $0.02765$							&$1.0091$ 		&$0.4729$ 			&$0.4231$	&$0.4414$&0.4324\\
\hline
HARFE ($\nu$) 						&$0.000140(2)$			&$\textbf{0.000110}(2)$ 			&$0.5350(2)$ 	&$\textbf{0.3122} (2)$		&$\textbf{0.2994}(2)$ &-&-\\
SHRIMP ($\nu$) 						&$0.000147(1)$			&$0.000190(2)$			&$0.3616(2)$	&$0.3501(1)$		&$0.4551(7)$&-&-\\
SALSA 										& $0.00918$					&$0.000135$					&$0.5470 $ 			&$0.3635$ 		&$0.3607$ &-&-\\
SRFE 											& $0.0154 $					&$0.0012$						&$0.5702$ 		&$0.4067$ 			&$0.6395$		&-&-\\
			\hline
  \end{tabular}
\caption{MSE on real datasets using various approximation techniques. Details to the corresponding algorithms can be found in HARFE \cite{HARFE}, SHRIMP \cite{Xie22}, SALSA \cite{SALSA}, SRFE~\cite{Ha21} or are available via ScikitLearn.jl. The best results for every dataset are highlighted. }
\label{tab:results}
\end{table}
\changed{We want to highlight two cases: the cases where the dataset has many samples $M$ compared to the dimension $d$, and second, the dataset has not much data samples available. 
In the first case (datasets Protein and Elevators) Algorithm~\ref{alg:1} performs better than the other machine learning algorithms. Since the random feature models set up a matrix as large as the number of unknowns, they can not handle such big datasets. 
Our Algorithm~\ref{alg:1} even provides similar or smaller approximation errors compared to the random features models in the case of a low sample size, with the advantage of being interpretable, that means it is easy to calculate the GSIs of the involved ANOVA terms for the final approximation model. In the application the user can use this to derive conclusions.
Furthermore, we confirm the thesis, that real-world data can be described by functions with low effective dimension. It should also be noted that in the dataset Airfoil additional $36$ dimensions with random noise are added and Algorithm~\ref{alg:1} easily finds the non-importance of these dimensions and reduces the data to the important $5$ dimensions. }


 \section{Conclusion \changed{and outlook}}
In this paper, we introduced a new method for function approximation from given fixed samples from an arbitrary density. 
This method combines previous work on least-squares approximation on the torus $\T^d$ and the truncation of the ANOVA decomposition with a variable transformation and a kernel density estimation. We are able to transfer the error decay rates and the fast algorithms from the torus to the domain $\Omega$. The new extension method, which benefits from the Chui-Wang wavelets, even allows the approximation of non-periodic functions. 
As shown in
our numerical experiments, this procedure is beneficial in function approximation. The code is available in the Julia package ANOVAapprox on GitHub, see  \url{https://github.com/NFFT/ANOVAapprox}.
\changed{We assume for our algorithm that the input variables are independent, which is not necessarily the case in applications. In future work we want to adapt our algorithm also to the setting of dependent input variables.}

\section*{Acknowledgments}
\changed{The authors thank the referees for useful comments about this work.} The authors would like to thank Tino Ullrich for several fruitful discussions. Next to that, the authors acknowledge funding by Deutsche Forschungsgemeinschaft (DFG, German Research Foundation) - Project-ID 416228727 - SFB 1410. Furthermore, the \changed{second} author acknowledges the support by the BMBF grant 01|S20053A (project SA$\ell$E).

\section*{Statements and Declarations}
\textbf{Conflict of interest}
The authors declare no competing interests.

\bibliographystyle{abbrv}

\begin{thebibliography}{10}

\bibitem{AdBrWe22}
B.~Adcock, S.~Brugiapaglia, and C.~G. Webster.
\newblock {\em Sparse polynomial approximation of high-dimensional functions}.
\newblock Computational science \& engineering. SIAM, Philadelphia,
  Pennsylvania, 2022.

\bibitem{AdHu19}
B.~Adcock and D.~Huybrechs.
\newblock Frames and numerical approximation.
\newblock {\em SIAM Review}, 61(3):443--473, 2019.

\bibitem{AdHu20}
B.~Adcock and D.~Huybrechs.
\newblock Approximating smooth, multivariate functions on irregular domains.
\newblock {\em Forum Math. Sigma}, 8:E26, 2020.

\bibitem{Be34}
E.~T. Bell.
\newblock Exponential polynomials.
\newblock {\em Ann. Math}, 35(2):258--277, 1934.

\bibitem{Boyd10}
J.~P. Boyd.
\newblock Six strategies for defeating the runge phenomenon in {G}aussian
  radial basis functions on a finite interval.
\newblock {\em Comput. Math. Appl.}, 60(12):3108--3122, 2010.

\bibitem{CaMoOw97}
R.~Caflisch, W.~Morokoff, and A.~Owen.
\newblock Valuation of mortgage-backed securities using {B}rownian bridges to
  reduce effective dimension.
\newblock {\em J. Comput. Finance}, 1(1):27--46, 1997.

\bibitem{chui92}
C.~K. Chui.
\newblock {\em An Introduction to Wavelets}.
\newblock Academic Press, Boston, 1992.

\bibitem{CoDaLe13}
A.~Cohen, M.~A. Davenport, and D.~Leviatan.
\newblock {On the stability and accuracy of least-squares approximations}.
\newblock {\em {Found. Comput. Math.}}, 13:819--834, 2013.

\bibitem{CoMi17}
A.~Cohen and G.~Migliorati.
\newblock Optimal weighted least-squares methods.
\newblock {\em SMAI J. Comput. Math.}, 3:181--203, 2017.

\bibitem{DuTeUl16}
D.~{D{\~u}ng}, V.~N. {Temlyakov}, and T.~{Ullrich}.
\newblock {\em Hyperbolic Cross Approximation}.
\newblock Advanced Courses in Mathematics -- CRM Barcelona. Birkh\"auser, Cham,
  2018.

\bibitem{DaKuUr99}
W.~Dahmen, A.~Kunoth, and K.~Urban.
\newblock Biorthogonal spline wavelets on the interval—stability and moment
  conditions.
\newblock {\em Appl. Comput. Harmon. Anal.}, 6(2):132--196, 1999.

\bibitem{DePeVo10}
R.~DeVore, G.~Petrova, and P.~Wojtaszczyk.
\newblock Approximation of functions of few variables in high dimensions.
\newblock {\em Constr. Approx.}, 33(1):125--143, 2010.

\bibitem{CoDo21}
M.~Dolbeault and A.~Cohen.
\newblock Optimal pointwise sampling for $l^2$ approximation.
\newblock {\em J. Complex.}, 68:101602, 2022.

\bibitem{GiKuSl22}
A.~D. Gilbert, F.~Y.Kuo, and I.~H. Sloan.
\newblock Equivalence between {S}obolev spaces of first-order dominating mixed
  smoothness and unanchored {ANOVA} spaces on {$\R^d$}.
\newblock {\em Math. Comp.}, 91:1837--1869, 2022.

\bibitem{Gra18}
A.~Gramacki.
\newblock {\em Nonparametric Kernel Density Estimation and Its Computational
  Aspects}, volume~37.
\newblock Springer International, 2018.

\bibitem{GrKuSl10}
M.~Griebel, F.~Y. Kuo, and I.~H. Sloan.
\newblock The smoothing effect of the {ANOVA} decomposition.
\newblock {\em J. Complex.}, 26(5):523--551, 2010.

\bibitem{Ha21}
A. Hashemi, H. Schaeffer, R. Shi, U. Topcu, G. Tran and R. Ward.
\newblock Generalization bounds for sparse random feature expansions.
\newblock {\em 	Appl. Comput. Harmon. Anal.},62:310--330, 2023.

\bibitem{Holtz11}
M.~Holtz.
\newblock {\em Sparse grid quadrature in high dimensions with applications in
  finance and insurance}, volume~77 of {\em Lecture Notes in Computational
  Science and Engineering}.
\newblock Springer-Verlag, Berlin, 2011.

\bibitem{Ho07}
G. Hooker
\newblock Generalized functional anova diagnostics for high-dimensional functions of dependent variables.
\newblock {\em J. Comput. Graph. Stat.}, 16(3):709--732, 2007.

\bibitem{Hu10}
D.~Huybrechs.
\newblock On the {F}ourier extension of nonperiodic functions.
\newblock {\em SIAM J. Numer. Anal.}, 47(6):4326--4355, 2010.

\bibitem{Jia09}
R.-Q. Jia.
\newblock Spline wavelets on the interval with homogeneous boundary conditions.
\newblock {\em Adv. Comput. Math.}, 30:177–200, 2009.

\bibitem{KaUlVo19}
L.~K\"{a}mmerer, T.~Ullrich, and T.~Volkmer.
\newblock Worst case recovery guarantees for least squares approximation using
  random samples.
\newblock {\em Constr. Approx.}, 54:295--352, 2021.

\bibitem{SALSA}
K. Kandasamy, Y. Yu.
\newblock Additive approximations in high dimensional nonparametric regression via
the salsa.
\newblock{International Conference on Machine Learning}, pp. 69–78. PMLR, 2016.

\bibitem{KuSiUl14}
T.~K\"uhn, W.~Sickel, and T.~Ullrich.
\newblock Approximation numbers of {S}obolev embeddings --- sharp constants and
  tractability.
\newblock {\em J. Complex.}, 30:95--116, 2014.

\bibitem{KuSlWaWo09}
F.~Y. Kuo, I.~H. Sloan, G.~W. Wasilkowski, and H.~Wo{\'{z}}niakowski.
\newblock On decompositions of multivariate functions.
\newblock {\em Math. Comp.}, 79(270):953--966, 2010.

\bibitem{LiPoUl21}
L.~Lippert, D.~Potts, and T.~Ullrich.
\newblock Fast hyperbolic wavelet regression meets {ANOVA}.
\newblock {\em Numer. Math.}, 154, 155-207, 2023.

\bibitem{LiOw06}
R.~Liu and A.~B. Owen.
\newblock Estimating mean dimensionality of analysis of variance
  decompositions.
\newblock {\em J. Amer. Statist. Assoc.}, 101(474):712--721, 2006.

\bibitem{nasdaladiss}
R.~Nasdala.
\newblock Efficient multivariate approximation with transformed rank-1
  lattices.
\newblock Dissertation, Fakult{\"a}t f{\"u}r Mathematik, Technische
  Universit{\"a}t Chemnitz, 2021.

\bibitem{NW08}
E.~Novak and H.~Wo\'zniakowski.
\newblock {\em Tractability of Multivariate Problems Volume I: Linear
  Information}.
\newblock Eur. Math. Society, EMS Tracts in Mathematics Vol 6, 2008.

\bibitem{NuSu21}
D.~Nuyens and Y.~Suzuki.
\newblock Scaled lattice rules for integration on {$\R^d$} achieving
  higher-order convergence with error analysis in terms of orthogonal
  projections onto periodic spaces.
\newblock {\em Math. Comp.}, 92:307--347, 2023.

\bibitem{PoSc19a}
D.~Potts and M.~Schmischke.
\newblock Approximation of high-dimensional periodic functions with
  {F}ourier-based methods.
\newblock {\em {SIAM} J. Numer. Anal.}, 59(5):2393--2429, 2021.

\bibitem{PoSc21}
D.~Potts and M.~Schmischke.
\newblock Interpretable approximation of high-dimensional data.
\newblock {\em SIAM J. Math. Data Sci.}, 3(4):1301--1323, 2021.

\bibitem{PoSc22}
D.~Potts and M.~Schmischke.
\newblock Interpretable transformed {ANOVA} approximation on the example of the
  prevention of forest fires.
\newblock {\em Front. Appl. Math. Stat.}, 8, 2022.

\bibitem{Ra14}
S.~Rahman.
\newblock Approximation errors in truncated dimensional decompositions.
\newblock {\em Math. Comput.}, 83(290):2799--2819, 2014.

\bibitem{Ra142}
S.~Rahman.
\newblock A Generalized ANOVA Dimensional Decomposition for Dependent Probability Measures.
\newblock {\em SIAM-ASA J. Uncertain.}, 2(1):670--697, 2014.

\bibitem{HARFE}
E. Saha, H. Schaeffer and G. Tran.
\newblock HARFE: hard-ridge random feature expansion. 
\newblock{\em Sampl. Theory Signal Process. Data Anal.}, 21, 27, 2023.


\bibitem{Diss_Schmischke}
M.~Schmischke.
\newblock Dissertation: Interpretable Approximation of High-Dimensional Data based on the ANOVA Decomposition.
\newblock {\em Universitaetsverlag Chemnitz}, 2022.  

\bibitem{ShJo91}
S.~J. Sheather and M.~C. Jones.
\newblock A reliable data-based bandwidth selection method for kernel density
  estimation.
\newblock {\em J. R. Stat. Soc., Series B, Methodological}, 53:683--690, 1991.

\bibitem{So01}
I.~M. Sobol.
\newblock Global sensitivity indices for nonlinear mathematical models and
  their {M}onte {C}arlo estimates.
\newblock {\em Math. Comput. Simul.}, 55(1-3):271--280, 2001.

\bibitem{TriebelIII}
H.~Triebel.
\newblock {\em Theory of Function Spaces III}.
\newblock Birkhäuser Basel, 1 edition, 01 2006.

\bibitem{WaJo95}
M.~Wand and M.~Jonas.
\newblock {\em Kernel Smoothing}, volume~60.
\newblock London ; New York : Chapman \& Hall, 1995.

\bibitem{WeKr23}
L.~Weidensager, D. Krumm, D. Potts and S. Odenwald.
\newblock Estimating vertical ground reaction forces from plantar
pressure using interpretable high-dimensional approximation.
\newblock{\em Sports Eng. (accepted)}, 2023. 

\bibitem{Wu2011}
C.~F.~J. Wu and M.~S. Hamada.
\newblock {\em Experiments - Planning, Analysis, and Optimization}.
\newblock John Wiley \& Sons, New York, 2011.

\bibitem{Xie22}
Y. Xie, B. Shi, H. Schaeffer and R. Ward.
\newblock Shrimp: Sparser random feature models via iterative magnitude pruning. 
\newblock{\em Math. Sci. Mach. Learn.}, PMLR 190, 303–318, 2022.

\end{thebibliography}

 \appendix{}
\section{Besov-Nikolskij-Sobolev spaces of mixed smoothness on the $d$-torus}
\label{sec:A1}
Here we summarize some relevant results from \cite[Chapt. 3]{DuTeUl16}. 
In particular we give the standard definition of the used function spaces. Let us first define Besov-Nikolskij spaces of mixed smoothness. We 
will use the classical definition via mixed moduli of smoothness. Therefore
first recall the basic concepts. For univariate functions $f:\T \to \C$ the
$m$-th difference operator $\Delta_h^{m}$ is defined by
\begin{equation*}
\Delta_h^{m}(f,x) := \sum_{j =0}^{m} (-1)^{m - j} \binom{m}{j} f(x +
jh),\quad x\in \T, h\in [0,1]\,.
\end{equation*}

Let $\vec u$ be any subset of $\{1,...,d\}$. For multivariate functions $f:\T^d\to
\C$ and $\vec h\in [0,1]^d$ the mixed $(m,\vec u)$-th difference operator $\Delta_{\vec h}^{m,\vec u}$
is defined by 
\begin{equation*}
\Delta_{\vec h}^{m,\vec u} := \
\prod_{i \in \vec u} \Delta_{h_i,i}^m\quad\mbox{and}\quad \Delta_{\vec h}^{m,\varnothing} =  \operatorname{Id},
\end{equation*}
where $\operatorname{Id}f = f$ and $\Delta_{h_i,i}^m$ is the univariate
operator applied to the $i$-th coordinate of $f$ with the other variables kept
fixed. 

\begin{definition}Let $s > 0$ and $1 \le p \le \infty$. Fixing an integer $m > s$, we define the space $\bB^s_{p,\infty}(\T^d)$ as the set of all $f\in L_p(\T^d)$ such that for any $\vec u \subset \{1,...,d\}$
\[
\big\|\Delta_{\vec h}^{m,\vec u}(f,\cdot)\big\|_{L_p(\T^d)}
\ \le \  
C\, \prod_{i \in \vec u} |h_i|^s
\]
for some positive constant $C$ and introduce the norm in this space
\[
\| \, f \, \|_{\bB^s_{p,\infty}} :=
\sum_{\vec u \subseteq \{1,...,d\}}
\, | \, f \, |_{\bB^s_{p,\infty}(\vec u)},
\]
where 
\[
| \, f \, |_{\bB^s_{p,\infty}(\vec u)} :=
\sup_{0 < |h_i| \le 1, \ i \in \vec u} \,  
\left(\prod_{i \in \vec u} |h_i|^{-s} \right) \,\big\| \, \Delta_{\vec h}^{m,\vec u}(f,\cdot) \, \big\|_{L_p(\T^d)}\,.
\]
\end{definition}

We define functions in a Sobolev space with dominating mixed derivatives,
\begin{equation*}
H^m_\mix(\T^d):=\left\{f:\T^d\to \C\mid \norm{f}_{H^m_\mix(\T^d)}<\infty\right\},
\end{equation*}
where the norm is defined by
\begin{equation}\label{eq:Hmnorm}
\norm{f}_{H^m_\mix(\T^d)}=\sum_{0\leq \norm{\vec k}_\infty \leq m}\left\|\Dx^{\vec k}f\right\|_{L_2(\T^d)},
\end{equation}
with the partial derivatives $\Dx^{\vec k} f = \tfrac{\partial^{k_1+\ldots+k_d} }{\partial x_1^{k_1}\cdots\partial x_d^{k_d} }$.
It clearly holds for $d=1$
\begin{equation*}
      H^m_{\mix}(\T) = H^m_2(\T) =: H^m(\T)\,.
\end{equation*}
The case $p=2$ and $\Omega=\T^d$ allows for a straight-forward extension to fractional smoothness parameters.

\begin{definition}
Let $s>0$. Then we define
\begin{equation*}H^s_{\mix}(\T^d):=\left\{f:\T^d\to \C\mid \norm{f}_{H^s_{\mix}(\T^d)}<\infty\right\},
\end{equation*}
where the norm is defined by
\begin{equation}\label{Fourier_char}
\norm{f}_{H^s_{\mix}(\T^d)}^2=\sum_{\vec k\in\Z^d}|c_{\vec k}(f)|^2\prod_{i=1}^d(1+|k_i|^2)^{s}.
\end{equation}
\end{definition}
This norm  is equivalent to the norm in \eqref{eq:Hmnorm} for $s\in \N$, see \cite{KuSiUl14}. We will 
consider the case where $s>\tfrac 12$, since in this case we have that $H^s_{\mix}(\T^d)	\hookrightarrow C(\T^d)$, which is necessary to sample the function.
There is a further useful equivalent norm which is based on a decomposition of $f$ in dyadic blocks. The dyadic blocks~\eqref{eq:def_J} and the 
decomposition~\eqref{eq:decomp_dyadic}
immediately give
\begin{equation*}\|f\|_{H^s_{\mix}}^2 \asymp \sum\limits_{\vec j \in \N_0^d} 2^{2|\vec j|_1s} \|f_{\vec j}\|_{L_2(\T^d)}^2.    
\end{equation*}

Interestingly, there is also a Fourier analytic characterization of the above defined Besov-Nikolskij spaces $\bB^s_{p,\infty}(\T^d)$ which even works for $1<p<\infty$. Instead of taking the $\ell_2$-norm of the weighted sequence $(2^{|\vec j|_1s}\|f_{\vec j}\|_{L_p(\T^d)})_{\vec j \in \N_0^d}$ we take the $\ell_{\infty}$-norm,   
\begin{equation}\label{dyadic_besov}
    \|f\|_{\bB^s_{p,\infty}} \asymp \sup\limits_{\vec j \in \N_0^d} 2^{|\vec j|_1 s }\|f_{\vec j}\|_{L_p(\T^d)}\,.
\end{equation}

\section{Cardinal B-Splines and Chui-Wang Wavelets}
We mostly work in this paper with Spline-wavelets, which have useful properties for our purposes. Therefore, we introduce here the cardinal B-Splines and the corresponding Chui-Wang-Wavelets. 
The cardinal B-Spline $B_m\colon \R\to \R$ of order $m$ is a piecewise polynomial function recursively defined by
\begin{equation}\label{eq:BSpline}
B_1(x)= 
\begin{cases}
1, & -\tfrac 12<x<\tfrac 12\\
0, &\text{otherwise},
\end{cases} \quad \text{and } \quad B_m(x) = \int_{x-\sfrac 12}^{x+\sfrac 12} B_{m-1}(t)\d t.
\end{equation}
The function $B_m(x)$ is a piecewise polynomial function of order $m-1$. Furthermore, the support 
of $B_m(x)$ is $(-\tfrac m2,\tfrac m2)$ and they are normalized by $\int_{-\sfrac m2}^{\sfrac m2}B_m(x)\dx x =1$. 
\begin{definition}
If we use the cardinal B-spline of order $m$ as scaling function $\phi = B_m$, the corresponding wavelet functions are the  
Chui-Wang wavelets \cite{chui92}, which are given by
\begin{equation*}
\psi(x) = \sum_{n}q_n B_m(2x-n-\tfrac m2),
\end{equation*}
where 
\begin{equation*}
q_n=\frac{(-1)^n}{2^{m-1}}\sum_{k=0}^m\binom{m}{k}B_{2m}(n+1-k-m).
\end{equation*}
\end{definition}
As in~\cite{LiPoUl21} we introduce the function 
\begin{equation}\label{eq:Psi_m}
\Psi_m(x):=\int_{-\infty}^x \frac{\psi(t)(x-t)^{m-1}}{(m-1)!}\d t,
\end{equation}
which is supported on $[0,2m-1]$ and fulfills $\Dx^m\Psi_m(x)=\psi(x)$ and is bounded, i.e.\,there exists a constant $C_m$, such that
\begin{equation}\label{eq:Psi_bounded}
|\Psi_m(x)|\leq C_m \text{ for all } x\in \R.
\end{equation}
Furthermore, this function has a nice property, which we use in Lemma~\ref{thm:nice_nonperiodic}. 
\begin{lemma}\label{lem:Psi_m=0}
For $k\in \Z$ we have for the function~\eqref{eq:Psi_m} that
$$\Psi_m(k)=0.$$
\end{lemma}
\begin{proof}
The result is a consequence of results from~\cite[Chapter 6]{chui92}. The wavelet function $\psi$ can also be written as
$$\psi(x) = L_{2m}^{(m)}(2x-1),$$
where $L_m(x)$ is the fundamental cardinal spline, which is a piece-wise polynomial of degree $m$ with the interpolation property $L_m(j)=\delta_{j,0}$ for $j\in \Z$.
Them we have $\Psi_m(x)= \tfrac{1}{2^m}L_{2m}(2x-1)$, which is zero for $x\in \Z$. 
\end{proof}

\end{document}